\documentclass[preprint,12pt]{elsarticle} 

\usepackage{url}
\usepackage{graphicx}
\usepackage{epsfig} 
\usepackage{hyperref}
\usepackage{natbib}
\usepackage[usenames, dvipsnames]{color}
\usepackage{amsmath,amssymb,amsfonts}
\usepackage{mathtools}
\usepackage{bm}
\usepackage{stmaryrd}
\usepackage{tikz}
\usepackage{marvosym}
\usepackage{pgfplotstable}
\usepackage{pgfplots}
\usepackage{courier}
\usepackage{geometry}
\usepackage{cleveref}
\usepackage{filecontents}
\usepackage[utf8]{inputenc}
\usepackage{booktabs}
\usepackage{tabularx}
\usepackage{multirow}
\usepackage{grffile}
\usepackage{subcaption}
\usepackage{algpseudocode}
\usepackage{algorithm}
\usepackage{xspace}
\usepackage{cancel}
\usepackage{mathrsfs}
\usepackage{longtable}
\usepackage{upgreek}

\usetikzlibrary{spy}
\usetikzlibrary{patterns}
\usetikzlibrary{calc,shadings}
\usetikzlibrary{snakes,arrows}
\usepgfplotslibrary{patchplots}
\usetikzlibrary{plotmarks}
\usetikzlibrary{calc,positioning}

\tikzset{
dimen/.style={<->,>=latex,thin,every rectangle node/.style={midway,above,font=\sffamily}}
}

\makeatletter
\def\els@aparagraph[#1]#2{\elsparagraph[#1]{#2}}
\def\els@bparagraph#1{\elsparagraph*{#1}}
\makeatother

\newcommand{\eref}[1]{Eq.~\ref{#1}}
\newcommand{\fref}[1]{Fig.~\ref{#1}}

\newcommand{\bl}[1]{\textcolor{black}{#1}}

\algnewcommand{\algorithmicgoto}{\textbf{goto}}%
\algnewcommand{\Goto}[1]{\algorithmicgoto~\ref{#1}}%

\def\addlegendimage{\csname pgfplots@addlegendimage\endcsname}

\pgfkeys{/pgfplots/number in legend/.style={%
        /pgfplots/legend image code/.code={%
            \node at (0.125,-0.0225){#1}; 
        },%
    },
}
\pgfplotsset{
every legend to name picture/.style={west}
}

\pgfplotsset{compat=newest}
  \usetikzlibrary{plotmarks}
  \usetikzlibrary{snakes,arrows}
  \usepgfplotslibrary{patchplots}

\begin{document}

\begin{frontmatter}

\title{Isogeometric analysis of insoluble surfactant spreading on a thin film}

\author[buw,uc1,inspire]{David Medina}
\author[buw]{Navid Valizadeh}
\author[uc1,uc2]{Esteban Samaniego}
\author[espol,usfq,inspire]{Alex X. Jerves}
\author[tdtu1,tdtu2,buw]{Timon Rabczuk\corref{Timon}}
\cortext[Timon]{Corresponding author at: Division of Computational Mechanics, Ton Duc Thang University, Ho Chi Minh City, Viet Nam.\\ {\it Email address:} timon.rabczuk@tdtu.edu.vn (Timon Rabczuk)}

\address[tdtu1]{Division of Computational Mechanics, Ton Duc Thang University, Ho Chi Minh City, Viet Nam}
\address[tdtu2]{Faculty of Civil Engineering, Ton Duc Thang University, Ho Chi Minh City, Viet Nam}
\address[buw]{Institute of Structural Mechanics, Bauhaus-Universit\"{a}t Weimar, Marienstrasse 15, 99423 Weimar, Germany}
\address[uc1]{Facultad de Ingeniería, Universidad de Cuenca, Ecuador}
\address[uc2]{Departamento de Recursos H{\'\i}dricos y Ciencias Ambientales, Universidad de Cuenca, Cuenca 010151, Ecuador}
\address[espol]{\bl{Escuela Superior Polit\'ecnica del Litoral, ESPOL, Facultad de Ciencias Naturales y Matem\'aticas, Campus Gustavo Galindo Km. 30.5 v\'{i}a Perimetral, P.O. Box 09-01-5863,  Guayaquil, Ecuador}}
\address[usfq]{\bl{College of Science \& Engineering, Universidad San Francisco de Quito, Quito, Pichincha 1712841, Ecuador}}
\address[inspire]{\bl{Fundaci\'{o}n {INSPIRE}, INS$\pi$RE, Quito, Pichincha, Ecuador}}

\begin{abstract}
In this paper we tackle the problem of surfactant spreading on a thin liquid film in the framework of isogeometric analysis. We consider a mathematical model that describes this phenomenon as an initial boundary value problem (IBVP) that includes two coupled fourth order partial differential equations (PDEs), one for the film height and one for the surfactant concentration. In order to solve this problem numerically, it is customary to transform it into a mixed problem that includes at most second order PDEs. However, the higher-order continuity of the approximation functions in Isogeometric Analysis (IGA) allows us to deal with the weak form of the fourth order PDEs directly, without the need of resorting to mixed methods. We demonstrate numerically that the IGA solution is able to reproduce results obtained before with mixed approaches. Complex phenomena such as Marangoni-driven fingering instabilities triggered by perturbations are easily captured.  
\end{abstract}

\begin{keyword}
 Isogeometric Analysis \sep Insoluble surfactant spreading \sep Thin film \sep High-order PDEs \sep Fingering instability
\end{keyword}

\end{frontmatter}

\section{Introduction}
Surfactant spreading on a thin liquid film is an interesting phenomenon both for
its biomedical and industrial applications and also for the mathematical and numerical challenges that arise from modeling it. Applications include coating flows, drainage flows in emulsions and foams, and surfactant replacement therapy \cite{warner_craster_matar_2004,jensen_grotberg_1992,Matar2001,YEO2001233}. 
The mathematical model describing the dynamics of surfactant spreading is usually derived by using lubrication theories and results in a system of coupled time-dependent nonlinear fourth-order partial differential equations. Solving such a strongly coupled problem which includes high-order spatial differential operators is not an easy task and requires robust numerical approaches. In this work, we present an isogeomtric finite element formulation \cite{Hughes20054135,Cottrell20091} for solving this problem.\\ 
Once a drop of surfactant is deposited on a thin liquid film over a smooth planar substrate, the initial difference between the lower surface tension of surfactant and the higher surface tension of surfactant-free thin film triggers the spreading of surfactant on the surface of thin film from the regions which are fully covered with surfactant to surfactant-free regions. This spreading process is mainly driven by the surface tension gradient (which is induced by the gradient in surfactant concentration) and causes a so-called Marangoni flow in the thin film \cite{jensen_grotberg_1992,warner_craster_matar_2004,Craster2009}. The Marangoni flow may cause large deformations in the thin film provided that the diffusion of the surfactant on the free surface of thin film is a slow process \cite{borgas_grotberg_1988,gaver_grotberg_1990,jensen_grotberg_1992}.
In such case, the Marangoni flow can be described by a ramped region with a sharp leading front \cite{warner_craster_matar_2004,jensen_grotberg_1992}. The sharp leading front is generated where the surfactant-covered film meets the surfactant-free undisturbed film. The peak value of the front is about two times the height of undisturbed film \cite{jensen_grotberg_1992}. Due to the mass conservation of liquid film, this elevation in film height is compensated by a film thinning region forming at the drop edge, just behind the ramped region. The fingering instability (see e.g. \cite{warner_craster_matar_2004,Craster2009} for more details) appears in this film thinning region.
As it will be shown later, a rough solid substrate may change the shape of ramped region and may also initiate the fingering instabilities.

Theoretically, the dynamics of surfactant spreading on thin liquid films can be modeled by a free boundary problem of Navier-Stokes equations for the incompressible flow of thin film and a surfactant mass balanace equation on the free surface. However, taking into account that the height of the film is considerably smaller than its other dimensions, one can use a lubrication theory with certain assumptions (see \cite{jensen_grotberg_1992,warner_craster_matar_2004} and the review paper \cite{Craster2009}) to approximate the problem by a system of coupled nonlinear fourth-order PDEs for the film height and the surfactant concentration. 

Several numerical methods have been proposed for solving surfactant spreading on thin films. \citet{jensen_grotberg_1992} carried out one-dimensional computations using second-order finite differences in space and Gear's method in time. \citet{warner_craster_matar_2004} designed a numerical scheme by using a finite element collocation in space and Gear's method in time for one-dimensional computations and \bl{a finite difference alternating direction implicit (ADI) scheme} for two-dimensional computations. \citet{Wong2011} split the equations into two parts, hyperbolic and parabolic systems, and solved the hyperbolic system using a finite volume method and the parabolic system using an ADI scheme. For a simplified form of the equations, in the absence of high-order regularization terms, \citet{PETERSON20125157} devised a numerical scheme which tracks the leading edge of the surfactant and compared its results with a hybrid method which uses a Godunov's method for the height equation and a an implicit finite difference method for the surfactant equation. We also refer to \cite{Swanson2015,Momoniat2013,GrLeRu02} for other finite difference and finite volume schemes which have been applied to the problem. \bl{While most of the  aforementioned methods enjoy properties such as high accuracy and/or ease of implementation, they are either limited to 1D problems, or require considerably more effort for modeling surfactant spreading over irregular geometries (this is generally an issue with finite difference and finite volume methods), or in some cases impose severe restrictions on the maximum allowable time step size \cite{Wong2011}. As an alternative, finite element methods can easily handle surfactant spreading problems on complex geometries and have shown better convergence behavior and larger time step sizes for fine meshes \cite{LIU2019429}.}
In the context of finite element methods, \citet{Barrett2003} presented a mixed finite element method which split the two fourth-order coupled nonlinear equations into three second-order equations by approximating the surfactant concentration, the film height, and the fluid pressure, thereby having three degrees of freedom per node.   
A different mixed finite element method was presented by \citet{LIU2019429} which approximates the surfactant concentration, the film height, the fluid pressure, and the surface and depth-averaged velocity fields, and results in seven degrees of freedom per node and a five-field system of coupled nonlinear equations. 
The authors showed that the fingering instabilities can be triggered by introducing perturbations to the roughness of solid substrate.
They also mentioned that, to their knowledge, their work was the first to incorporate the effect of substrate roughness explicitly into the mathematical model of surfactant spreading on thin liquid films. However, we remark that a more complete model, which considers \bl{the effects of substrate roughness, substrate incline, and intermolecular interactions,} had been already presented in an excellent review paper by \citet{Craster2009}, among others. 
Moreover, \cite{LIU2019429} addressed the effect of surfactant concentration on fluid pressure. We note that in practice the term representing the effect of surfactant concentration on fluid pressure is negligible with respect to other terms present in approximating the fluid pressure \cite{Craster2009} and therefore it is reasonable that the effect of surface tension variations on capillarity has been neglected in previous studies \cite{jensen_grotberg_1992,Barrett2003,warner_craster_matar_2004,Craster2009,Wong2011}. 


As mentioned earlier, to discretize fourth-order PDEs using the mixed finite element method with globally $C^0$-continuous basis functions, each of the equations should be transformed into a system
of PDEs which at the most involve second-order spatial derivative operators; this allows for the use of lower order finite elements. 
However, a more natural finite element approach to discretizing fourth-order PDEs is a primal variational formulation by dealing directly with the weak form of fourth order PDEs using globally $C^1$-continuos basis functions.   
In this contribution, we take this latter approach by using a spatial discretization based on isogeometric analysis (IGA) \cite{Hughes20054135,Cottrell20091} and exploiting the high continuity of high-order spline basis functions. For temporal discretization, we use the generalized-$\alpha$ time integration method. Numerical experiments show that the proposed approach is able to match the results of previous methods used to describe the physical phenomena. Specifically, we are able to study the so-called fingering phenomena and its sensitivity to initial conditions for the height of the thin liquid film. We also study the effect of the roughness of the solid substrate in the solution. We use sufficiently refined meshes which rules out the possibility of numerical instabilities due to the degeneracy of the fourth-order film height equation \cite{Bertozzi1999}.

The rest of this paper is organized as follows: Sec.~\ref{sec:mathModel} presents a detailed derivation of the mathematical model for insoluble surfactant spreading on thin liquid films deposited on rough solid substrates. We start from a free boundary problem of Navier-Stokes equations for the incompressible flow and a surfactant mass balanace equation on the free surface and show how the final system of two coupled fourth-order PDEs for the film height and the surfactant concentration are obtained using lubrication theory by considering certain approximations. Sec.~\ref{sec:NumericalFormulation} is devoted to the numerical formulation of the problem in an isogeometric framework, where details of the weak forms, semi-discrete formulation and the time discretization are given. In Sec.~\ref{sec:NumericalExamples}, we first verify the results of our isogeometric formulation by making comparisons with similarity solutions and other numerical schemes. Thereafter, we consider several examples of fingering instabilities triggered either by perturbations in the initial film height or the roughness of solid substrate. Finally, conclusions are drawn in Sec.~\ref{sec:conclusions}.   
\section{Mathematical model} \label{sec:mathModel}
In this section, we present the mathematical model for describing the spreading of a drop of liquid, uniformly covered with an insoluble surfactant, on an initially-undisturbed thin film of the same liquid as the drop. The thin film is resting on a rough solid substrate and is in contact with air at its free surface. \\
We consider the thin liquid film to be an incompressible Newtonian fluid with a constant dynamic viscosity $\mu$ and a constant density $\rho$. Let $\mathcal{L}$ and $\mathcal{H}$ be the characteristic length and height of the film, respectively. The characteristic height of the film is assumed to be significantly smaller than its characteristic length and therefore we can define a fineness ratio as $\epsilon=\frac{\mathcal{H}}{\mathcal{L}}\ll1$. Based on this assumption, we can simplify the thin film's Navier--Stokes and continuity equations, surfactant mass balance equation, and related boundary conditions by only considering leading-order terms in $\epsilon$. The thin film is considered to be rested on a rough solid substrate represented by a roughness function $f$. \\
The mathematical model for surfactant spreading on thin films resting on smooth surfaces (see \cite{jensen_grotberg_1992,warner_craster_matar_2004,Wong2011,Craster2009}) was extended to thin films on rough surfaces in \cite{Craster2009,LIU2019429}, among others. 
The model is based on the assumption that the magnitude of roughness of solid substrate is small enough (\bl{much smaller than the film height}) so that it does not violate the main assumptions of lubrication theory \cite{Craster2009,LIU2019429}. Our derivation is mainly based on this mathematical model. In what follows, we derive the governing equations. The final equations are similar to the ones in \cite{LIU2019429} except of a difference due to assumptions we made for approximating the pressure field in liquid film; similar to \cite{jensen_grotberg_1992,warner_craster_matar_2004,Craster2009,Wong2011} we neglect the effects of surface tension variations on capillarity and assume that it is only the term corresponding to the constant minimum surface tension which should be present in the approximation of the film's pressure. 
This is \bl{a valid assumption} as in practice the term representing the effects of surface tension variations on capillarity is much less than the other term (see e.g. \cite{Craster2009}). 

\subsection{Governing equations}
Let us assume a Cartesian coordinate system, $\left(x,y,z\right)$, where $x$, $y$, and $z$ correspond to the streamwise, transverse and vertical coordinates. The thin liquid film is bounded from below by a rough solid substrate, located at $z=f(x,y)$, and from above by a free surface which is located at $z=h(x,y,t)$. Here, $f(x,y)$ is the roughness function of the substrate, $h(x,y,t)$ represents the height of liquid film and $t$ denotes time.
The governing equations of the thin liquid film are Navier-Stokes and continuity equations for incompressible flows which are given by:
\begin{align} 
\frac{\partial\bm{u}}{\partial t}+\bm{u}\cdot\nabla\bm{u}+\frac{1}{\rho}\nabla p -\nu \nabla^2\bm{u}-\bm{g}&=\bm{0}, \label{LBLME}\\
\nabla\cdot\bm{u} &=0. \label{LBME}
\end{align}
where $\rho$ is the fluid density, $\bm{u}$ is the fluid velocity vector, $p$ is the fluid pressure, $\nu=\frac{\mu}{\rho}$ is the kinematic viscosity, $\mu$ is the dynamic viscosity, and $\bm{g}$ is the gravitational body force per unit mass. The components of velocity vector are given by $\bm{u}=\left(u,v,w\right)$ and $\bm{g}=-g\mathbf{e}_z$, where $g$ denotes the acceleration of gravity and $\mathbf{e}_z=\left(0,0,1\right)$ is the standard unit basis vector in the $z$ direction. \\
On the free surface of the thin liquid film, we have normal and tangential stress boundary conditions. Noting that the free surface is essentially a gas(air)--liquid interface, and assuming the air to be inviscid, the normal stress boundary condition is described as (see \cite{Bush2010} for full details)
\begin{align} \label{NSBC}
\bm{n}\cdot\hat{\bm{T}}\cdot\bm{n}-\bm{n}\cdot\bm{T}\cdot\bm{n} = \sigma (\nabla_s\cdot\bm{n}) \qquad \text{at $z=h$}
\end{align}   
where $\bm{n}$ denotes the outward unit normal vector to the free surface, $\nabla_s$ is the surface gradient operator, and $\sigma$ is the surface tension. In addition,
$\hat{\bm{T}}$ and $\bm{T}$ are the stress tensors of the air and the film, respectively, which are defined as, 
\begin{align} 
\hat{\bm{T}}&= -\hat{p}\bm{I} \nonumber \\
\bm{T} &= -p\bm{I} + \mu \left(\nabla \bm{u}+\nabla^T\bm{u} \right) \nonumber
\end{align}  
Also, considering $\bm{t}$ to be any unit tangent vector to
the interface, the tangential stress boundary condition is given by (see \cite{Bush2010})
\begin{align} \label{TSBC}
\bm{n}\cdot\hat{\bm{T}}\cdot\bm{t}-\bm{n}\cdot\bm{T}\cdot\bm{t} =- \nabla_s\sigma\cdot\bm{t} \qquad \text{at $z=h$}
\end{align} 
We consider a kinematic boundary condition on the free surface, $z=h(x,y,t)$. On the free surface, we have 
\begin{align} 
z-h(x,y,t)=0 \nonumber
\end{align} 
which holds at all times. Taking the material time derivative of this equation yields
\begin{align} 
&\frac{D \left(z-h(x,y,t)\right)}{Dt}=0, \nonumber\\
& \frac{\partial \left(z-h(x,y,t)\right)}{\partial t}+\bm{u}_s\cdot\nabla \left(z-h(x,y,t)\right) = 0. \label{levelset-Dt}
\end{align} 
where $\bm{u}_s=\left(u_s,v_s,w_s\right)=\left(u\big\rvert_{z=h(x,y,t)},v\big\rvert_{z=h(x,y,t)},w\big\rvert_{z=h(x,y,t)}\right)$ is the fluid velocity $\bm{u}$ evaluated at the free surface.\\
On the bottom surface of the film, $z=f(x,y)$, we assume the no-slip and no-penetration boundary conditions for $\bm{v} = (u,v)$ which is the projection of the velocity vector $\bm{u}$ onto
$xy$-plane,
\begin{align} \label{no-slip}
\bm{v}=\bm{0} \qquad \text{at} \qquad z=f
\end{align} 
Moreover, on the free surface, we have the mass balance equation for surfactant concentration $c$ which is expressed by ( see  \cite{Wong1996,cermelli_fried_gurtin_2005,Pereira2007} for details of derivation),
\begin{align} \label{mass_balance}
\frac{\partial c}{\partial t} + \left(\bm{u}\cdot\bm{n}\right)\bm{n}\cdot \nabla c +\nabla_s\cdot(c\bm{u})=D\Delta_s c
\end{align}
where $\Delta_s$ is the surface Laplacian which is defined as $\Delta_s=\nabla_s^2=\nabla_s\cdot\nabla_s$, and $D$ is the surface diffusion coefficient of surfactant. 

\subsection{Scaling and assumptions}
	Considering $\mathcal{L}$, $\mathcal{H}$, and $\epsilon=\frac{\mathcal{H}}{\mathcal{L}}$ as characteristic length and height of the film and fineness ratio, respectively, and defining $\mathcal{U}$ and $\mathcal{P}$ as the characteristic velocity and pressure of the fluid, respectively, we use the following scalings to nondimensionalize the governing equations, 
	\begin{equation} \label{dimensionless}
	\begin{split}
	&x^*=\frac{x}{\mathcal{L}}, \qquad y^*=\frac{y}{\mathcal{L}}, \qquad z^*=\frac{z}{\mathcal{H}}=\frac{z}{\epsilon \mathcal{L}},\qquad t^*=\frac{\mathcal{U}t}{\mathcal{L}}, \qquad h^*=\frac{h}{\mathcal{H}},\qquad f^*=\frac{f}{\mathcal{H}},\\
	& u^*=\frac{u}{\mathcal{U}}, \qquad v^*=\frac{v}{\mathcal{U}},\qquad w^*=\frac{w}{\epsilon~\mathcal{U}}, \qquad p^*=\frac{p}{\mathcal{P}}, \qquad \hat{p}^* = \frac{\hat{p}}{\mathcal{P}}.\\
	\end{split}
	\end{equation} 
	where superscript stars denote dimensionless quantities.
	Moreover, we scale both the surfactant concentration and the surface
	tension by the following relations:
	\begin{equation} \label{dimensionless2}
	\begin{split}
	& \qquad c= c_m c^*, \\ 
	&\sigma = \sigma_m + S \sigma^*(c^*).
	\end{split} 
	\end{equation} 
	where $c_m$ is a characteristic concentration defined as the concentration of surfactant at saturation. When surfactant concentration attains $c_m$ the surface tension takes a constant minimum value $\sigma_m$. $S=\sigma_0 -\sigma_m$ is a spreading coefficient and signifies the difference between the maximum value of surface tension $\sigma_0$ and its minimum value $\sigma_m$. \\
	In addition, we define $\bm{v}_s=\bm{v}(x,y,z=h)$ as the projected velocity field at the free surface, $\bar{\bm{v}}=\frac{1}{h-f}\int_{f}^{h} \bm{v}(x,y,z) ~dz$ as the depth-averaged velocity of the liquid film, $\nabla_p=(\frac{\partial}{\partial x},\frac{\partial}{\partial y})$ as the planar gradient operator, and $\Delta_p=\frac{\partial^2}{\partial x^2}+\frac{\partial^2}{\partial y^2}$ as the planar Laplace operator. These quantities are scaled as follows,  
	\begin{equation} \label{dimensionless3}
	\begin{split}
	& \nabla_p = \frac{1}{\mathcal{L}}\nabla^*_p, \quad \Delta_p = \frac{1}{\mathcal{L}^2}\Delta^*_p,\quad \nabla_p\sigma=\frac{S}{\mathcal{L}}\nabla^*_p\sigma^*,\\
	&\bm{v} = \mathcal{U}\bm{v}^*,\qquad \bm{v}_s = \mathcal{U}\bm{v}^*_s,\qquad \bar{\bm{v}} = \mathcal{U}\bar{\bm{v}}^*
	\end{split}
	\end{equation} 	
	We further define several dimensionless parameters which will appear through nondimensionalization of governing equations:
	\begin{equation} \label{dimensionless4}
	\begin{split}
	& \mathrm{Re}=\frac{\rho\mathcal{U}\mathcal{L}}{\mu}, \qquad \Lambda = \frac{\mu\mathcal{U}\mathcal{L}}{\mathcal{H}^2\mathcal{P}}, \qquad \Upsilon=\frac{\mu\mathcal{U}\mathcal{L}}{\rho g \mathcal{H}^3}, \qquad \mathrm{Ca}=\frac{\sigma_m}{\mu \mathcal{U}}, \qquad \mathrm{M}=\frac{\epsilon S}{\mu \mathcal{U}},\\
	&\mathscr{C} = \frac{\epsilon^2 \sigma_m}{S} = \epsilon^3 \mathrm{Ca},\qquad
	\mathscr{G}=\frac{\rho g \mathcal{L}^2}{\sigma_m}=\frac{1}{\Upsilon\mathscr{C}}, \qquad \mathrm{Pe} = \frac{S \mathcal{H}}{\mu D}=\frac{\mathcal{U} \mathcal{L}}{D}, 
	\end{split}
	\end{equation} 
where $\mathrm{Re}$ is the Reynolds number, $\Lambda$ is the so-called bearing number \cite{KUNDU2016409} which is defined as the ratio of viscous forces to pressure forces, $\Upsilon$ is a dimensionless number which is indicative of the ratio of viscous forces to gravity forces, $\mathrm{Ca}$ is the capillary number which is reflecting the relative significance of capillary to viscous forces, $\mathrm{M}$ is Marangoni number which is reflecting the relative significance of Marangoni forces (driven by surface tension gradients) to viscous forces,  
$\mathscr{C}$ is a capillarity parameter, $\mathscr{G}$ is a gravitational parameter, and $\mathrm{Pe}$ is a surface Peclet number which is indicative of the relative importance of Marangoni-driven spreading to surface diffusive spreading \cite{warner_craster_matar_2004}. \\
We will also make several assumptions: As we are modeling thin film flows for which Marangoni effects are dominant (see also surface-tension-gradient-dominated flows in \cite{Craster2009}), Marangoni number is set to unity, $\mathrm{M}=1$, which gives the characteristic (Marangoni) velocity as $\mathcal{U} = \frac{S\mathcal{H}}{\mu\mathcal{L}}$. The pressure scaling is chosen as $\mathcal{P}=\frac{\mu\mathcal{U}\mathcal{L}}{\mathcal{H}^2}$ which renders $\Lambda=1$. We assume that $\epsilon\ll1$, Reynolds number is finite, $\epsilon^2\mathrm{Re}\ll1$, and $\Upsilon$ is near unity. Similar to \cite{Craster2009}, we further assume that $\frac{\sigma_m}{S}\gg |\sigma^*|$, and that the concentration field (and as a result the surface tension) is independent of $z$ coordinate.

\subsection{Reduced set of dimensionless governing equations based on lubrication theory approximations}
Introducing the scalings into governing equations and taking into account the assumptions described in the previous section, one can find a set of approximate equations which are valid for $\epsilon\ll1$. Accordingly, to leading order in $\epsilon$, \cref{LBLME,LBME,NSBC,TSBC,levelset-Dt,no-slip,mass_balance} can be described in a dimensionless form as,
\paragraph{Navier--Stokes and continuity equations}
\begin{align} 
&\nabla^*_p p^* = \frac{\partial^2 \bm{v}^*}{\partial z^{{*}^2}}, \label{RE1-2}\\
&\frac{\partial p^*}{\partial z^*} +\frac{1}{\Upsilon}=0, \label{RE3}\\
&\frac{\partial u^*}{\partial x^*}+\frac{\partial v^*}{\partial y^*}+\frac{\partial w^*}{\partial z^*}=0, \label{RE4}
\end{align} 
\paragraph{Normal stress boundary condition}
\begin{align} 
& p^* = \hat{p}^*-\mathscr{C}\Delta^*_p h^*, \qquad \text{at $z^*=h^*$} \label{RE5}
\end{align}  
\paragraph{Tangential stress boundary condition}
\begin{align} 
& \frac{\partial \bm{v}^*}{\partial z^*} = \nabla^*_p\sigma^*, \qquad \text{at $z^*=h^*$} \label{RE6}
\end{align}  
\paragraph{No-slip, no-penetration boundary condition}
\begin{align} 
& \bm{v}^* = \bm{0}, \qquad \text{at $z^*=f^*$} \label{RE7}
\end{align} 
\paragraph{Kinematic boundary condition}
\begin{align} 
& \frac{\partial h^*}{\partial t^*} + u^*_s \, \frac{\partial h^*}{\partial x^*} + v^*_s \, \frac{\partial h^*}{\partial y^*} - w^*_s = 0,  \label{RE8}
\end{align}
\paragraph{Surfactant mass balance equation}
\begin{align} 
& \frac{\partial c^*}{\partial t^*} + \nabla^*_p\cdot\left(c^* \bm{v}^*_s\right)= \frac{1}{\mathrm{Pe}}\Delta^*_p c^*.  \label{RE9}
\end{align}

\subsection{Deriving the dimensionless evolution equations}
Integrating \eref{RE3} and applying the normal stress boundary condition (\eref{RE5}), the dimensionless pressure field in the liquid film is given by,
\begin{align} \label{pressure}
& p^* = \hat{p}^*+\frac{1}{\Upsilon} (h^*-z^*) -\mathscr{C} \Delta^*_p h^*
\end{align}
Substituting \cref{pressure} into \cref{RE1-2} yields,
\begin{align} \label{LBLME5}
\frac{\partial^2 \bm{v}^*}{\partial z^{{*}^2}} &= \nabla^*_p ( \frac{1}{\Upsilon} \, h^* -\mathscr{C} \, \Delta^*_p h^*)
\end{align}
Integrating \cref{LBLME5} twice and applying the tangential stress boundary condition (\cref{RE6}) at the free surface, and the no-slip and no-penetration boundary condition (\cref{RE7}) at the bottom surface of the film, the dimensionless velocity field $\bm{v}^*$ is obtained as,
\begin{align} \label{projectedVelocity}
\bm{v}^* &= \left(\frac{1}{\Upsilon} \, \nabla^*_p h^* - \mathscr{C} \, \nabla^*_p\Delta^*_p h^*\right)\left(\frac{1}{2}(z^*-f^*)^2-(h^*-f^*)(z^*-f^*)\right)+(z^*-f^*)\nabla^*_p \sigma^*
\end{align}     
Now, we can define the dimensionless velocity field at the free surface of the thin liquid film as 
\begin{align} \label{FreeSurfaceVelocity}
\bm{v}^*_s = \bm{v}^*(x^*,y^*,z^*=h^*) = \frac{(h^*-f^*)^2}{2}\left( \mathscr{C} \, \nabla^*_p\Delta^*_p h^* -\frac{1}{\Upsilon} \, \nabla^*_p h^* \right) +(h^*-f^*)\nabla^*_p \sigma^*
\end{align}  
Next, we integrate \cref{RE4} with respect to $z^*$ from $z^*=f^*$ to $z^*=h^*$. Using \cref{RE7}, we obtain $w^*_s$ as,
\begin{align} \label{w_s}
w^*_s = -\int_{f^*}^{h^*} \frac{\partial u^*}{\partial x^*} ~dz^* -\int_{f^*}^{h^*} \frac{\partial v^*}{\partial y^*} ~dz^* 
\end{align}
Using Leibniz's integral rule and \cref{RE7}, \cref{w_s} can be written as,
\begin{align} \label{w_s-Leibniz}
w^*_s = &-\frac{\partial}{\partial x^*}\left(\int_{f^*}^{h^*} u^* ~dz^*\right) + \frac{\partial h^*}{\partial x^*} u^*\big\rvert_{z^*=h^*} - \cancelto{0}{\frac{\partial f^*}{\partial x^*} u^*\big\rvert_{z^*=f^*}} \nonumber\\ &-\frac{\partial}{\partial y^*}\left(\int_{f^*}^{h^*} v^* ~dz^*\right) + \frac{\partial h^*}{\partial y^*} v^*\big\rvert_{z^*=h^*} - \cancelto{0}{\frac{\partial f^*}{\partial y^*} v^*\big\rvert_{z^*=f^*}}
\end{align}  
hence, 
\begin{align} \label{w_s-Leibniz-2}
w^*_s = -\frac{\partial}{\partial x^*}\left(\int_{f^*}^{h^*} u^* ~dz^*\right) + \frac{\partial h^*}{\partial x^*} u^*\big\rvert_{z^*=h^*} -\frac{\partial}{\partial y^*}\left(\int_{f^*}^{h^*} v^* ~dz^*\right) + \frac{\partial h^*}{\partial y^*} v^*\big\rvert_{z^*=h^*}
\end{align}
Substituting \cref{w_s-Leibniz-2} into \cref{RE8} gives,
\begin{align} \label{evolution-h}
\frac{\partial h^*}{\partial t^*}+\frac{\partial}{\partial x^*}\left(\int_{f^*}^{h^*} u^* ~dz^*\right)+\frac{\partial}{\partial y^*}\left(\int_{f^*}^{h^*} v^* ~dz^*\right)=0.
\end{align}
Using the projected velocity vector (\cref{projectedVelocity}) for computing the integrals in the second and third terms of the above equation, the evolution equation for the film height is obtained as, 
\begin{align} \label{evolution-h-2}
\frac{\partial h^*}{\partial t^*}+\nabla^*_p\cdot\left((h^*-f^*)\bar{\bm{v}}^*\right)=0.
\end{align}
where
\begin{align} \label{evolution-h-3}
\bar{\bm{v}}^*= \frac{1}{h^*-f^*}\int_{f^*}^{h^*} \bm{v}^*(x^*,y^*,z^*) ~dz^* = \frac{(h^*-f^*)^2}{3}\left(\mathscr{C}\nabla^*_p\Delta^*_p h^*-\frac{1}{\Upsilon} \nabla^*_p h^*\right)+\frac{(h^*-f^*)}{2}\nabla^*_p\sigma^*
\end{align}
Finally, \cref{RE9,FreeSurfaceVelocity,evolution-h-2,evolution-h-3} together with a surfactant equation of state (which will be discussed in Sec.~\ref{equation-of-state}) constitute a coupled system of fourth-order PDEs. These equations are subject to suitable boundary and initial conditions.  
\subsection{Surfactant equation of state} \label{equation-of-state}
The surfactant equation of state, $\sigma^*(c^*)$ in \cref{dimensionless2}, describes how the surface tension of the thin film is influenced by the surfactant concentration. We consider the following three variants of dimensionless equations of state:
\begin{align} \label{eq-surface-tension}
& \sigma^*(c^*) = \frac{\alpha+1}{\left(1+\Theta(\alpha)c^*\right)^3}-\alpha,\qquad \Theta(\alpha)=\left(\frac{\alpha+1}{\alpha}\right)^{1/3}-1.\\
& \sigma^*(c^*) = 1-c^*,\\
&\sigma^*(c^*)= 
\begin{dcases}
\left(1-c^*\right)^3,& c^*\leq 1\\
0,              & c^*>1
\end{dcases}
\end{align}
where $\alpha=\frac{\sigma_m}{S}$ is a positive constant.
The first equation was originally proposed in \cite{SHELUDKO1967391}. The second equation is a linear equation of state which is achieved by considering $\alpha\to\infty$ in the first equation. This equation or its non-linear version has been used in several studies (see e.g. \cite{gaver_grotberg_1990,jensen_grotberg_1992,Barrett2003,warner_craster_matar_2004}) for modeling a monolayer of surfactant. The last equation of state was proposed in \cite{Wong2011} for the case of multiple layers of surfactant and was recently used in \cite{LIU2019429}.  
\subsection{Strong form of dimensionless partial differential equations governing the evolution of surfactant concentration and film height. Initial and boundary conditions}
Let $\Omega^* \subset \mathbb{R}^2$ be an open set that represents the computational domain, $\partial \Omega^*$ the boundary of $\Omega^*$, and $\bm{m}^*$ the outward unit normal vector to $\partial \Omega^*$. The strong form of the initial/boundary-value problem for the surfactant concentration and the thin film height in dimensionless coordinates is stated as: Given a time interval of interest $\left[0,T^*\right]$, find $c^*: \overline{\Omega}^*\times\left(0,T^*\right)\mapsto \mathbb{R}$ and $h^*: \overline{\Omega}^*\times\left(0,T^*\right)\mapsto \mathbb{R}$ such that
\begin{align} 
&\frac{\partial c^*}{\partial t^*}+\nabla^*_p\cdot\left(c^*\bm{v}_s^*\right)=\frac{1}{\mathrm{Pe}}\Delta^*_p c^* &\text{in}\quad &\Omega^*\times\left(0,T^*\right),\label{evolution-eqs-c}\\
&\bm{v}^*_s= \mathscr{C}\frac{(h^*-f^*)^2}{2}\left(\nabla^*_p\Delta^*_p h^*-\mathscr{G} \nabla^*_p h^*\right)+(h^*-f^*)\nabla^*_p\sigma^*,\\
&\frac{\partial h^*}{\partial t^*}+\nabla^*_p\cdot\left((h^*-f^*)\bar{\bm{v}}^*\right)=0   &\text{in}\quad &\Omega^*\times\left(0,T^*\right),\label{evolution-eqs-h}\\
&\bar{\bm{v}}^*= \mathscr{C}\frac{(h^*-f^*)^2}{3}\left(\nabla^*_p\Delta^*_p h^*-\mathscr{G} \nabla^*_p h^*\right)+\frac{(h^*-f^*)}{2}\nabla^*_p\sigma^*,\\
& c^*\bm{v}^*_s \cdot \bm{m}^* = 0  &\text{on}\quad &\partial\Omega^*\times\left(0,T^*\right), \label{bc1}\\
& (h^*-f^*) \bar{\bm{v}}^* \cdot \bm{m}^* = 0 &\text{on}\quad &\partial\Omega^*\times\left(0,T^*\right),\label{bc2}\\
& \nabla^*_p h^* \cdot \bm{m}^* = 0  &\text{on}\quad &\partial\Omega^*\times\left(0,T^*\right),\label{bc3}\\
& c^*(\bm{x}^*,0) = c^*_0  &\text{on}\quad &\bm{x}^*\in\Omega^*,\\
& h^*(\bm{x}^*,0) = h^*_0  &\text{on}\quad &\bm{x}^*\in\Omega^*.
\end{align}
In what follows, we will use the dimensionless form of the equations; the superscript stars will be omitted for the sake of notational simplicity. Moreover, as we will only use planar gradient and Laplace operators, we will omit their $p$ subscript henceforth.

\section{Numerical formulation} \label{sec:NumericalFormulation}
\subsection{Continuous problem in the weak form}
Let $\mathcal{S}= \left \{ c \mid c (\cdot, t) \in H^{2}(\Omega) \right \}$ be the trial solution space and
$\mathcal{V} = \left \{ w \mid w \in H^{2}(\Omega) \right \}$ the weighting function space, where $H^{2}(\Omega)$ represents a Sobolev space of square-integrable functions with square-integrable first and second derivatives in $\Omega$. For the time being, let us also assume periodic boundary conditions in all directions and strongly impose these conditions on the finite element spaces. The weak form of the evolution \cref{evolution-eqs-c,evolution-eqs-h} is derived by multiplying these equations by weighting functions and applying integration by parts repeatedly. The variational formulation is therefore stated as follows: find $\mathbf{U} = \left \{c,h \right \} \in \mathcal{S}^{2}$ such that for all $ \mathbf{W} = \left \{ w,q \right \} \in \mathcal{V}^{2}$  
\begin{align}	
\mathcal{B} \left(\mathbf{W},\mathbf{U}\right) &= 0, \label{eq:weak} \\
\left(w,c(\bm{x},0)\right)_{\Omega}&=\left(w,c_{0}\right)_{\Omega},   \label{eq:weakq}\\
\left(q,h(\bm{x},0)\right)_{\Omega}&=\left(q,h_{0}\right)_{\Omega}. \label{eq:weakw}
\end{align}
where
\begin{equation} \label{eq:weak2}
\begin{split}
\mathcal{B}\left(\mathbf{W},\mathbf{U}\right)  =  
&\left( w, \frac{\partial c}{\partial t}\right)_{\Omega} +
\left(\nabla w, \mathscr{C} \, \mathscr{G} \, \mathcal{M}_{2}(h_p) \, c \, \nabla h \right)_{\Omega}
-\left(\nabla w, \mathcal{M}_1(h_p) \, \sigma'(c) \, c \, \nabla c \right)_{\Omega}\\
&+\left(\nabla w, \mathscr{C} \, \mathcal{M}_2(h_p) \, \Delta h \, \nabla c \right)_{\Omega} 
+ \left ( \Delta w, \mathscr{C} \, \mathcal{M}_2(h_p) \, c \, \Delta h \right)_{\Omega}\\
&+\left ( \nabla w, \mathscr{C} \mathcal{M}'_2(h_p) \, \nabla (h-f) \, c \, \Delta h \right  )_{\Omega}
- \left ( w, \frac{1}{\mathrm{Pe}}  \Delta c \right)_{\Omega} \\
&+\left( q, \frac{\partial h}{\partial t}\right)_{\Omega} 
+\left( \nabla q, \mathscr{C} \, \mathcal{M}'_3(h_p) \, \nabla(h-f) \,\Delta h \right)_{\Omega} 
+\left( \nabla q, \mathscr{C} \, \mathscr{G} \,\mathcal{M}_3(h_p) \, \nabla h \right)_{\Omega} \\
&-\left( \nabla q, \mathcal{M}_2(h_p) \, \sigma'(c) \, \nabla c \right)_{\Omega}
+\left( \Delta q, \mathscr{C} \, \mathcal{M}_3(h_p) \, \Delta h \right)_{\Omega}
\end{split} \end{equation}
Here, $(\cdot,\cdot)_{\Omega}$ denotes the $L^{2}$ inner product with respect to the domain $\Omega$, $h_p$ denotes the physical height which is defined as $h_p=h-f$, $\mathcal{M}_1(h_p)=h_p$, $\mathcal{M}_2(h_p)=h_p^2/2$, $\mathcal{M}_3(h_p)=h_p^3/3$, and the prime denotes differentiation with respect to the argument.
\bl{\paragraph*{Nitsche's method for imposing Dirichlet boundary condition.} Following the variational formulations of Nitsche's method for fourth-order PDEs \cite{Embar2010,Harari2015,ZHAO2017177}, we impose the Dirichlet boundary condition (\eref{bc3}) by adding the following boundary terms to the residual (\eref{eq:weak2})
\begin{equation} \label{eq:weak-dirichlet}
\begin{split}
&-\left( \nabla w\cdot\bm{m},\mathscr{C} \,  \mathcal{M}_2(h_p) \, c \, \Delta h \right)_{\partial\Omega}
-\left( \nabla q\cdot\bm{m}, \mathscr{C} \, \mathcal{M}_3(h_p) \, \Delta h \right)_{\partial\Omega}
-\left( \Delta q, \mathscr{C} \, \mathcal{M}_3(h_p) \, \nabla h \cdot\bm{m} \right)_{\partial\Omega}\\
&+\left( \nabla q\cdot\bm{m}, \alpha_h \, \nabla h\cdot\bm{m} \right)_{\partial\Omega}
\end{split} \end{equation}
where $\alpha_h$ is a stabilization parameter. Here, the assumption is that no admissibility constraints are imposed on trial and weighting function spaces. The consistency of the Nitsche's--based variational formulation with the strong form of boundary value problem (\cref{evolution-eqs-c,evolution-eqs-h,bc1,bc2,bc3}) can be shown by applying repeated integration by parts to the augmented residual. This leads to  
\begin{align}	
&\left( w,\frac{\partial c}{\partial t}+\nabla\cdot\left(c\bm{v}_s\right)-\frac{1}{\mathrm{Pe}}\Delta c \right)_{\Omega} -\left(  w,c\bm{v}_s\cdot \bm{m} \right)_{\partial\Omega} = 0, \label{vc:c} \\
&\left( q,\frac{\partial h}{\partial t}+\nabla\cdot\left((h-f)\bar{\bm{v}}\right) \right)_{\Omega}
-\left(  q,(h-f) \bar{\bm{v}} \cdot \bm{m} \right)_{\partial\Omega} \nonumber\\
&-\left( \left[\mathscr{C} \, \mathcal{M}_3(h_p) \, \Delta q - \alpha_h (\nabla q\cdot\bm{m})\right],\nabla h\cdot \bm{m} \right)_{\partial\Omega} = 0.   \label{vc:h}
\end{align}
The above equations must hold for arbitrary weighting functions $w$ and $q$. Invoking this assumption yields the Euler-Lagrange equations which recovers the strong form of the boundary value problem (\cref{evolution-eqs-c,evolution-eqs-h,bc1,bc2,bc3}). This completes the proof of variational consistency for the Nitsche's--based variational formulation.
}

\subsection{Semidiscrete formulation}
Following the Galerkin formulation, we let $\mathcal{S}^{h}\subset \mathcal{S}$ and $\mathcal{V}^{h} \subset \mathcal{V}$ be the finite-dimensional trial and weighting function spaces such that $\mathcal{S}^{h}=\mathcal{V}^{h}=\text{span}\{N_A\}_{A=1}^{n_b}$ where $N_A$ denotes the spline basis function, and $n_b$ is the dimension of the discrete function space which implies $N_A$'s are linearly independent. To define the functions $N_A$ over the domain in the physical space, the spline basis functions with $C^1$ global continuity are defined on a parametric space and then mapped into the physical space $\Omega$ using the isoparametric concept. Having a non-singular mapping, the discrete solution possess $C^1$ global continuity in the physical space which ensures the satisfaction of the requirement $\mathcal{S}^{h}\subset H^{2}(\Omega)$. \\
The semidiscrete Galerkin form of the problem is stated as: Find $\mathbf{U}^{h}=\left \{ c^{h}, h^{h}\right \} \in (\mathcal{S}^{h})^{2}$ such that for all $\mathbf{W}^{h}=\left \{ w^{h}, q^{h}\right \} \in (\mathcal{V}^{h})^{2}$
\begin{align}	
\mathcal{B} \left(\mathbf{W}^h,\mathbf{U}^h\right) &= 0, \label{eq:semidiscrete} \\
\left(w^h,c^h(\bm{x},0)\right)_{\Omega}&=\left(w^h,c_{0}\right)_{\Omega},   \label{eq:semidiscrete-w}\\
\left(q^h,h^h(\bm{x},0)\right)_{\Omega}&=\left(q^h,h_{0}\right)_{\Omega}. \label{eq:semidiscrete-q}
\end{align}
where
\begin{align*} 
& c^h\left(\bm{x},t\right)=\sum_{A=1}^{n_{bf}}N_A(\bm{x})c_A(t), \qquad w^h\left(\bm{x}\right)=\sum_{A=1}^{n_{bf}}N_A(\bm{x})w_A,\\
& h^h\left(\bm{x},t\right)=\sum_{A=1}^{n_{bf}}N_A(\bm{x}) h_A(t), \qquad
q^h\left(\bm{x}\right)=\sum_{A=1}^{n_{bf}}N_A(\bm{x})q_A.
\end{align*} 
\subsection{Time discretization}
The above semidiscrete Galerkin formulation is discretized in time using the generalized $\alpha$-method \cite{Chung1993371}. We discretize the time interval of interest $\left(0,T\right)$ into $N$ subintervals $I_n=\left(t_n,t_{n+1}\right)$, $n=0,\dots,N-1$, where $0=t_0<t_1<\cdots<t_{N-1}<t_N=T$. The time step size is defined as $\Delta t_n = t_{n+1}-t_n$. We define the time discrete approximations of $c^h(\cdot,t_n)$ and $h^h(\cdot,t_n)$ as $c^h_n$ and $h^h_n$, respectively. 
The global vectors of control variables associated to $c^{h}_n$ and $h^{h}_n$ are written, respectively, as $\mathbf{C}_n$ and $\mathbf{H}_n$, while their first time derivatives are denoted as $\dot{\mathbf{C}}_n$ and $\dot{\mathbf{H}}_n$. The global solution vector at time step $n$ is therefore defined as $\mathbf{S}_n=\{\mathbf{C}_n,\mathbf{H}_n\}^T$. Let us also define the residual vector as
\begin{align} 
&\bm{R}=\begin{Bmatrix}\bm{R}^{c} \\ \bm{R}^{h}\end{Bmatrix},\\
& \bm{R}^{c} = \{R^{c}_A\}, \qquad \bm{R}^{h} = \{R^{h}_A\},\\
& R^{c}_A =\mathcal{B}\left(\{N_A,0\},\{c^h,h^h\}\right), \qquad R^{h}_A=\mathcal{B}\left(\{0,N_A\},\{c^h,h^h\}\right).
\end{align} 
The generalized-$\alpha$ time integration scheme is then stated as follows: given $\mathbf{S}_n$, $\dot{\mathbf{S}}_n$, and $\Delta t_n$, find $\mathbf{S}_{n+1}$, $\dot{\mathbf{S}}_{n+1}$,$\mathbf{S}_{n+\alpha_f}$, and $\dot{\mathbf{S}}_{n+\alpha_m}$ such that
\begin{align} \label{eq:R-Coupled}
&\bm{R}(\dot{\mathbf{S}}_{n+\alpha_m},\mathbf{S}_{n+\alpha_f})=\bm{0}.
\end{align}
where
\begin{align}
& \dot{\mathbf{S}}_{n+\alpha_m} = \dot{\mathbf{S}}_{n} + \alpha_m(\dot{\mathbf{S}}_{n+1}-\dot{\mathbf{S}}_{n}),\\
& \mathbf{S}_{n+\alpha_f} = \mathbf{S}_{n} + \alpha_f(\mathbf{S}_{n+1}-\mathbf{S}_{n}),\\
& \mathbf{S}_{n+1} = \mathbf{S}_{n} +\Delta t_n ((1-\gamma)\dot{\mathbf{S}}_{n} +\gamma\dot{\mathbf{S}}_{n+1}). \label{newmark-Coupled}
\end{align}
The real-valued parameters $\alpha_f$, $\alpha_m$, and $\gamma$ in the above equations are determined such that the accuracy and stability of the time integration scheme is guaranteed. It was shown in \cite{Jansen2000305} that
for a first-order linear ordinary differential equation system, the generalized-$\alpha$ method attains second-order accuracy in time if
\begin{equation} \label{eq:gamma}
\gamma = \frac{1}{2}+\alpha_m-\alpha_f
\end{equation} 
and the method is unconditionally stable if 
\begin{equation}
\alpha_m\geqslant\alpha_f\geqslant\frac{1}{2}
\end{equation} 
The parameters $\alpha_m$ and $\alpha_f$ are defined as \cite{Jansen2000305}:
\begin{equation}\label{eq:alpha_m}
\alpha_m = \frac{1}{2}\left(\frac{3-\rho_\infty}{1+\rho_\infty}\right), \qquad \alpha_f = \frac{1}{1+\rho_\infty}
\end{equation}
where $\rho_\infty\in[0,1]$ is a parameter which controls the high-frequency damping and is defined as the spectral radius of the amplification matrix as $\Delta t \longrightarrow\infty$. In this work, we choose $\rho_\infty=0.5$ which has been shown to be an effective choice through several previous studies including turbulent computations \cite{Bazilevs2007173}, the Cahn-Hilliard phase-field model \cite{GOMEZ20084333} and other high-order phase-field models \cite{VALIZADEH2019599}.\\
The nonlinear system of equations of (\ref{eq:R-Coupled}) are solved using the Newton-Raphson method according to the following two-phase predictor-multicorrector algorithm:
\paragraph{Predictor phase} Set
\begin{align} 
& \mathbf{S}_{n+1}^{(0)}=\mathbf{S}_{n},\\
& \dot{\mathbf{S}}_{n+1}^{(0)}=\frac{\gamma-1}{\gamma}\dot{\mathbf{S}}_{n}.
\end{align}
The superscript 0 denotes the iteration index of the nonlinear solver. Note that the predictor is consistent with \eref{newmark-Coupled} which is the Newmark formula.
\paragraph{Multicorrector phase} Repeat the following steps for $i=1,2,\dots,i_{max}$, or until convergence is reached.
\begin{enumerate}
\item Evaluate iterates at the intermediate time levels:
\begin{align} 
& \dot{\mathbf{S}}_{n+\alpha_m}^{(i)} = \dot{\mathbf{S}}_{n} + \alpha_m(\dot{\mathbf{S}}_{n+1}^{(i-1)}-\dot{\mathbf{S}}_{n}),\\
& \mathbf{S}_{n+\alpha_f}^{(i)} = \mathbf{S}_{n} + \alpha_f(\mathbf{S}_{n+1}^{(i-1)}-\mathbf{S}_{n}).
\end{align}
\item Use the intermediate iterates to compute $\Delta\dot{\mathbf{S}}_{n+1}^{(i)}$ from the following linearized equation:
\begin{align} \label{eq:Linearized-Coupled}
& \bm{R}_{n+1}^{(i+1)}=\bm{R}_{n+1}^{(i)}+\frac{\partial\bm{R}_{n+1}^{(i)}}{\partial\dot{\mathbf{S}}_{n+1}^{(i)}}\Delta\dot{\mathbf{S}}_{n+1}^{(i)}=\bm{0}
\end{align}
Note that this equation is obtained by expanding the residual about a previous $i^{\text{th}}$-iterate solution for the residual at $(i+1)^{\text{th}}$-iterate.
\item Use $\Delta\dot{\mathbf{S}}_{n+1}^{(i)}$ to update the solution.
 \begin{align} 
& \dot{\mathbf{S}}_{n+1}^{(i)} = \dot{\mathbf{S}}_{n+1}^{(i-1)} + \Delta\dot{\mathbf{S}}_{n+1}^{(i)},\\
& \mathbf{S}_{n+1}^{(i)} = \mathbf{S}_{n+1}^{(i-1)} + \gamma \Delta t_n \Delta\dot{\mathbf{S}}_{n+1}^{(i)}.
 \end{align}
 Note that the last equation is obtained by subtracting \eref{newmark-Coupled} evaluated at iteration $(i-1)$ from \eref{newmark-Coupled} evaluated at iteration $(i)$.
\end{enumerate} 
The tangent matrix in \eref{eq:Linearized-Coupled} is given by
\begin{align} \label{eq:tangent-coupled-1}
 \bm{K}=&\frac{\partial\bm{R}(\dot{\mathbf{S}}_{n+\alpha_m},\mathbf{S}_{n+\alpha_f})}{\partial\dot{\mathbf{S}}_{n+\alpha_m}}\frac{\partial\dot{\mathbf{S}}_{n+\alpha_m}}{\partial\dot{\mathbf{S}}_{n+1}}
+\frac{\partial\bm{R}(\dot{\mathbf{S}}_{n+\alpha_m},\mathbf{S}_{n+\alpha_f})}{\partial\mathbf{S}_{n+\alpha_f}}\frac{\partial\mathbf{S}_{n+\alpha_f}}{\partial\dot{\mathbf{S}}_{n+1}} \nonumber\\
=&\alpha_m\frac{\partial\bm{R}(\dot{\mathbf{S}}_{n+\alpha_m},\mathbf{S}_{n+\alpha_f})}{\partial\dot{\mathbf{S}}_{n+\alpha_m}}
+\alpha_f\gamma\Delta t_n \frac{\partial\bm{R}(\dot{\mathbf{S}}_{n+\alpha_m},\mathbf{S}_{n+\alpha_f})}{\partial\mathbf{S}_{n+\alpha_f}}
 \end{align}
where the iteration index $i$ is omitted for clarity, and
 \begin{align} \label{eq:tangent-coupled-2}
& \frac{\partial\bm{R}(\dot{\mathbf{S}}_{n+\alpha_m},\mathbf{S}_{n+\alpha_f})}{\partial\dot{\mathbf{S}}_{n+\alpha_m}} = \begin{bmatrix} \frac{\partial\bm{R}^{c}(\dot{\mathbf{S}}_{n+\alpha_m},\mathbf{S}_{n+\alpha_f})}{\partial\dot{\mathbf{C}}_{n+\alpha_m}} & \frac{\partial\bm{R}^{c}(\dot{\mathbf{S}}_{n+\alpha_m},\mathbf{S}_{n+\alpha_f})}{\partial\dot{\mathbf{H}}_{n+\alpha_m}}\\
  \frac{\partial\bm{R}^{h}(\dot{\mathbf{S}}_{n+\alpha_m},\mathbf{S}_{n+\alpha_f})}{\partial\dot{\mathbf{C}}_{n+\alpha_m}} & \frac{\partial\bm{R}^{h}(\dot{\mathbf{S}}_{n+\alpha_m},\mathbf{S}_{n+\alpha_f})}{\partial\dot{\mathbf{H}}_{n+\alpha_m}} \end{bmatrix}, \nonumber\\
&\frac{\partial\bm{R}(\dot{\mathbf{S}}_{n+\alpha_m},\mathbf{S}_{n+\alpha_f})}{\partial \mathbf{S}_{n+\alpha_f}}  = 
\begin{bmatrix}
\frac{\partial\bm{R}^{c}(\dot{\mathbf{S}}_{n+\alpha_m},\mathbf{S}_{n+\alpha_f})}{\partial\mathbf{C}_{n+\alpha_f}} & \frac{\partial\bm{R}^{c}(\dot{\mathbf{S}}_{n+\alpha_m},\mathbf{S}_{n+\alpha_f})}{\partial\mathbf{H}_{n+\alpha_f}} \\
\frac{\partial\bm{R}^{h}(\dot{\mathbf{S}}_{n+\alpha_m},\mathbf{S}_{n+\alpha_f})}{\partial\mathbf{C}_{n+\alpha_f}} & \frac{\partial\bm{R}^{h}(\dot{\mathbf{S}}_{n+\alpha_m},\mathbf{S}_{n+\alpha_f})}{\partial\mathbf{H}_{n+\alpha_f}} 
\end{bmatrix}. 
  \end{align} 
  
\bl{\subsection{Adaptive time-stepping algorithm}  
To save the computational time while maintaining a good level of time accuracy, we use an adaptive time-stepping algorithm which is provided in the PETSc library \cite{petsc-web-page,petsc-user-ref}. This algorithm has been proven to be an effective option for time adaptivity in the context of isogeometric analysis with generalized-$\alpha$ time integration \cite{GOMEZ20084333,BARTEZZAGHI2015446,LIU2013321,VALIZADEH2019599}.
Borrowing ideas from embedded Runge-Kutta methods \cite{Gustafsson:1994,HAIRER1987}, the algorithm estimates a local truncation error by comparing the solutions computed by generalized-$\alpha$ method and the backward Euler method. Based on this error, the algorithm decides on the next time step size. Let $i$ denotes the iterate of the adaptive time-stepping algorithm. At the time step $t_n$, 
given $\dot{\mathbf{S}}_n$, $\mathbf{S}_n$, and $\Delta t_{n-1}$, and assuming $\Delta t_{n,0}=\Delta t_{n-1}$, we repeat the following steps for $i=1,\dots,i_{\text{max}}$  
\begin{enumerate}
	\item Compute $\mathbf{S}_{n+1,(i-1)}^{GA}$ using the generalized-$\alpha$ method and $\Delta t_{n,(i-1)}$.
	\item Compute $\mathbf{S}_{n+1,(i-1)}^{BE}$ using the backward Euler method and $\Delta t_{n,(i-1)}$.
	\item Calculate the weighted local truncation error according to (See \cite{HAIRER1987}, Eq. 4.11) \\
	 $$e_{n+1,(i-1)}=\sqrt{\frac{1}{N}\sum_{a=1}^{N}\left(\frac{|\mathbf{S}_{n+1,(i-1)}^{GA}(a)-\mathbf{S}_{n+1,(i-1)}^{BE}(a)|}{\textit{tol}_{n+1,(i-1)}(a)}\right)^2}$$
	 where $N$ is the dimension of the solution vector. The tolerance level in the denominator is computed by   
	 $$\textit{tol}_{n+1,(i-1)}(a) = \textit{tol}_{A}(a) + \max\left(|\mathbf{S}_{n+1,(i-1)}^{GA}(a)|,|\mathbf{S}_{n+1,(i-1)}^{BE}(a)|\right) \textit{tol}_{R}(a)
	 $$ 
	 where $\textit{tol}_{A}(a)$ and $\textit{tol}_{R}(a)$ are the desired absolute and relative tolerances, respectively.
	 \item Update the time step size based on the following formula\\
	$$\Delta t_{n,(i)}= \Delta t_{n,(i-1)} \min\left(\alpha_{max},\max\left(\alpha_{min},\beta \, e_{n+1,(i-1)}^{-1/2}\right)\right)$$
	 where $\beta<1$ is a safety factor, and $\alpha_{min}$ and $\alpha_{max}$ are used to ensure that the change in the time step size is kept within a certain factor. In our implementation, we use the following default values: $\beta=0.9$, $\alpha_{min}=0.1$, $\alpha_{max}=10$, $\textit{tol}_{A}(a)=10^{-4}$, $\textit{tol}_{R}(a)=10^{-4}$.
	 \item If {$e_{n+1,(i-1)}>1$}, go to step 1 and repeat the steps; otherwise, set  $\dot{\mathbf{S}}_{n+1}=\dot{\mathbf{S}}_{n+1,(i-1)}^{GA}$, $\mathbf{S}_{n+1}=\mathbf{S}_{n+1,(i-1)}^{GA}$, and $\Delta t_{n}=\Delta t_{n,i}$.
\end{enumerate}
}	 
  	
	\section{Numerical examples} \label{sec:NumericalExamples}
	Unless stated otherwise, we consider the surfactant equation of state to be $\sigma(c)=1-c$ in the following numerical examples. An adaptive time-stepping scheme (as explained in \cite{VALIZADEH2019599}) is utilized. All the computations are done in the dimensionless coordinates. Uniform cubic B-spline elements with $C^2$ global continuity are employed in all cases. The values of dimensionless parameters are chosen \bl{based on typical experimental parameter values} such that the use of lubrication theory is permissible (for more details, see \cite{warner_craster_matar_2004,Wong2011,jensen_grotberg_1992,SINZ2011519}). \bl{On the boundaries of the computational domain, we either impose periodic boundary conditions or boundary conditions according to \cref{bc1,bc2,bc3}. In the latter case, we utilize the Nitsche's method as derived in \eref{eq:weak-dirichlet}. The stability parameter of Nitsche's method is chosen as $5/h_e$ where $h_e$ is the characteristic length of each boundary element; this is an empirical choice following the works of \citet{ZHAO2017177} and \citet{WELLS2006860}.}
	The computer code for our isogeometric formulation has been implemented in PetIGA \cite{PetIGA}. PetIGA is a NURBS-based isogeometric analysis framework that allows for  parallel computing, and is built on top of the scientific library PETSc \cite{petsc-web-page}. 
	\subsection{Verification example: Comparing the surfactant spreading rates with  similarity solutions}
\begin{figure}
\centering
 	\begin{subfigure}[t]{0.5\textwidth}
       \pgfplotsset{legend style={at={(0.98,0.05)},
       		anchor=south east,font=\tiny}}
      \begin{tikzpicture}
            \begin{loglogaxis}[xlabel=$t$,ylabel=$r_s$,
            xmin=1e-1, xmax=5.e2,
           	ymin=1e0, ymax=1e1,
            	    xticklabel style={/pgf/number format/sci},
            	    grid=both,
            	tick align=inside,
            	tickpos=left] 
            	\addplot[blue,mark=triangle*] table [x=h,y=$Cubic$]{
            	h	$Cubic$
            	1.0		2.29
            	2.0		2.63
            	4.0		3.07
           		8.0		3.58
        		16.0	4.18
        		32.0	4.94
        		64.0    5.82
            	128.0   6.91
            	            };            
              \addplot[black] table [x=h,y=$DropSlope$]{
                       h	$DropSlope$
                       1.0	    2.2
                       128.0    7.4      
                       };
               \addlegendentry{IGA}
               \addlegendentry{Slope m=1/4}
            \end{loglogaxis}
            \end{tikzpicture}
      \caption{}
         \label{fig:SpreadingRate-Drop-Nitsche}
\end{subfigure} 
\begin{subfigure}[t]{0.5\textwidth}
       \pgfplotsset{legend style={at={(0.98,0.05)},
       		anchor=south east,font=\tiny}}
      \begin{tikzpicture}
            \begin{loglogaxis}[xlabel=$t$,ylabel=$x_s$,
            xmin=1e-1, xmax=5.e2,
           	ymin=1e0, ymax=2e1,
            	    xticklabel style={/pgf/number format/sci},
            	    grid=both,
            	tick align=inside,
            	tickpos=left] 
            	\addplot[blue,mark=triangle*] table [x=h,y=$Cubic$]{
            	h	$Cubic$
            	1.0		2.42
            	2.0		2.88
            	4.0		3.54
           		8.0		4.43
        		16.0	5.55
        		32.0	6.99
        		64.0    8.8
            	128     11.09
            	            };            
              \addplot[black] table [x=h,y=$DropSlope$]{
                       h	$DropSlope$
                       1.0	    2.4
                       128.0    12.095      
                       };
               \addlegendentry{IGA}
               \addlegendentry{Slope m=1/3}
            \end{loglogaxis}
            \end{tikzpicture}
      \caption{}
         \label{fig:SpreadingRate-Strip}
\end{subfigure}
\caption{\bl{The leading edge of the surfactant at different times computed by isogeometric analysis for two different initial distributions of surfactants: (a) axisymmetric drop for which the surfactant spreading scales as $r_s\sim t^{1/4}$, and (b) planar strip for which the surfactant spreading scales as $x_s\sim t^{1/3}$.}}	
 \label{fig:SpreadingRate}	 
\end{figure}
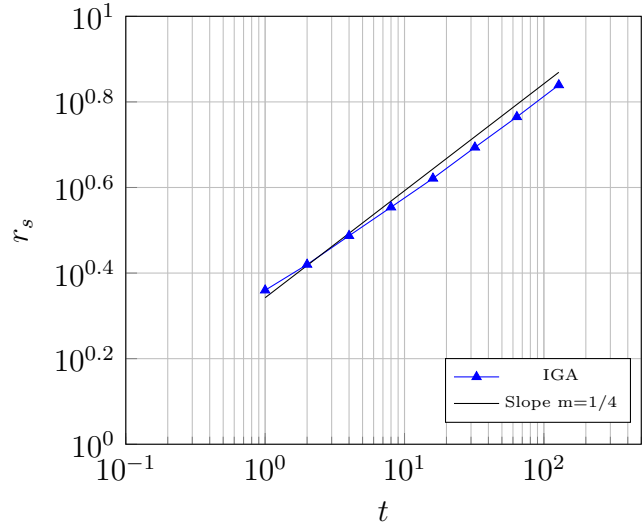
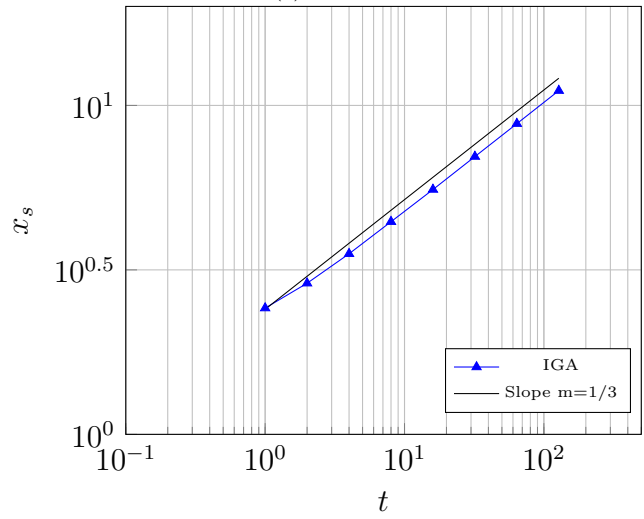

\bl{For the first set of numerical examples, we consider a simplified form of the \cref{evolution-eqs-c,evolution-eqs-h} for which there exist similarity solutions.} \bl{Assuming $\mathscr{C}=\mathscr{G}=f(x,y)=0$ and $\frac{1}{\mathrm{Pe}}\to 0$, the simplified equations are given by
\begin{align} 
&\frac{\partial c}{\partial t}+\nabla\cdot\left(c \, h \,\nabla\sigma\right)=0, &\text{in}\quad &\Omega\times\left(0,T\right),\label{simplified-evolution-eqs-c}\\
&\frac{\partial h}{\partial t}+\nabla\cdot\left(\frac{h^2}{2} \, \nabla\sigma\right)=0,   &\text{in}\quad &\Omega\times\left(0,T\right).\label{simplified-evolution-eqs-h}
\end{align}} 
These equations possess similarity solutions which show that the leading edge of surfactant spreads as $t^{1/4}$ for an initial axisymmetric drop of surfactant-covered liquid \cite{jensen_grotberg_1992}. The same study shows that when the initial distribution of surfactant is a planar strip, the spreading scales as $t^{1/3}$. These similarity solutions were obtained assuming that the total mass of surfactant is fixed. We model these problems numerically using isogeometric analysis with $256\times256$ cubic B-spline elements and an adaptive time-stepping scheme. The initial conditions for the surfactant and the film height are considered to be
\begin{align*}
& h(x,y,0)=1, \\   
& c(x,y,0)= \frac{1}{2} (1-\tanh(10(r-1))) , \qquad  \text{where } r=\sqrt{x^{2}+y^{2}}, \qquad
\text{(axisymmetric drop)}\\
& c(x,y,0)= \frac{1}{2} (1-\tanh(10(x-1))) , \qquad
\text{(planar strip)}
\end{align*}
The computational domain is $[-8,8]\times[-8,8]$ for the drop and $[0,16]\times[0,16]$ for the planar strip. \bl{For the planar strip problem, we consider periodic boundary conditions in the $y$ direction}. \fref{fig:SpreadingRate} presents the variation of the leading edge of surfactant front with respect to time for the axisymmetric drop and the planar strip. In both cases, our isogeomteric solution agrees well with the spreading rates predicted by the similarity solutions. The leading edge of surfactant is determined such that ahead of its position $c=0$. We note that the solution to the simplified equation for the film height has a shock-like structure which calls for the use of stabilization techniques such as Streamline-Upwinded Petrov-Galerkin method (SUPG) \cite{BROOKS1982199} and discontinuity-capturing techniques \cite{Bazilevs2007} (see also \cite{LIU2019429} for the use of SUPG in a mixed finite element framework for studying surfactant spreading on thin liquid films and \cite{BAZILEVS2020112638} for using SUPG plus a residual-based discontinuity-capturing technique for isogeometric analysis of Eulerian hydrodynamics). We will not pursue this direction as the solution to the full PDEs, especially when the contributions from capillary, gravity and surface diffusion terms are reasonably significant, is considerably smoothened and therefore usually no stabilization is needed when these terms are included (Cf. \cite{LIU2019429}).       
\subsection{Verification example: Spreading of a surfactant drop on a thin liquid film}
\begin{figure}
 	\centering
 	\begin{subfigure}{0.48\textwidth}
 	\centering
 	\includegraphics[width=\linewidth,trim={0cm 0cm 0cm 0cm},clip]{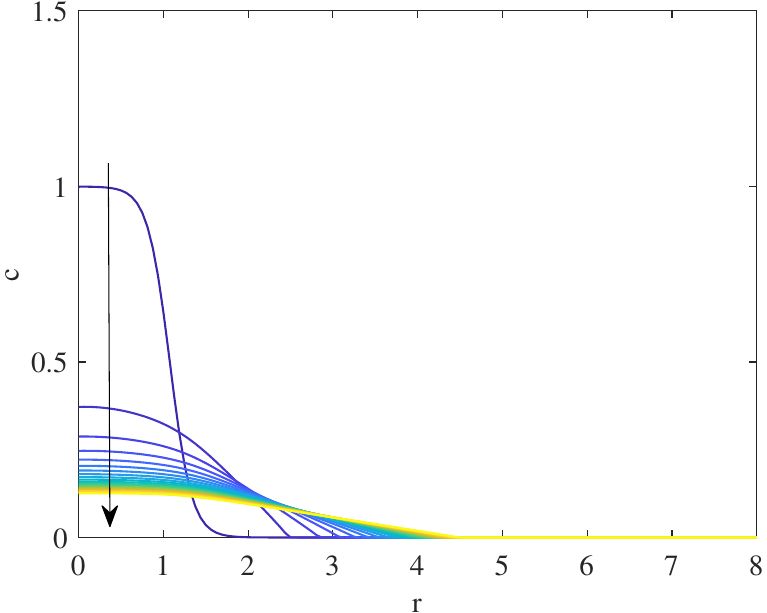}
 	\caption{}
 	\end{subfigure}
 	\quad
 	\begin{subfigure}{0.48\textwidth}
 	\centering
 		\includegraphics[width=\linewidth,trim={0cm 0cm 0cm 0cm},clip]{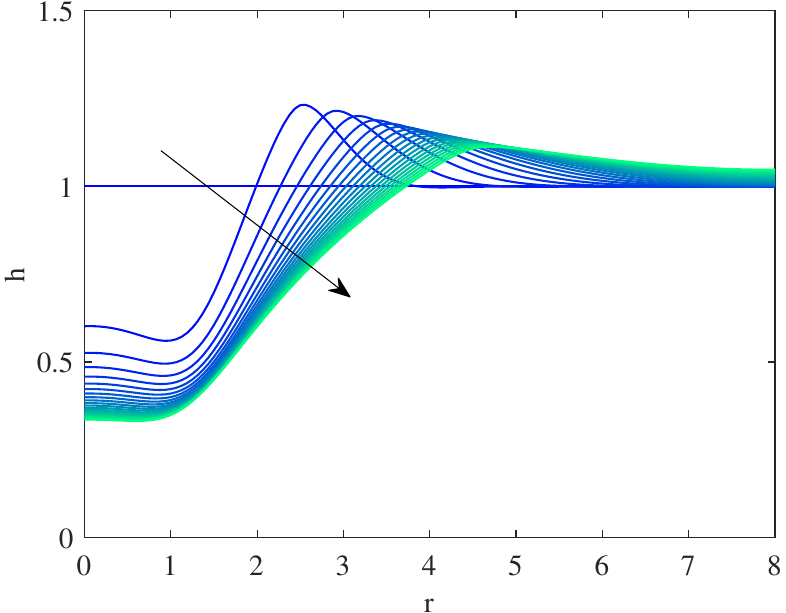}
 		\caption{}
    \end{subfigure}	
 \caption{\bl{Isogeometric solution for spreading of a surfactant drop on a thin liquid film rested on a smooth solid substrate: The evolution of (a) surfactant concentration $c$ and (b) film height $h$ for $t=[0,50]$ in time steps of $2.5$. The arrows indicate the direction of increasing time. A color grading is also used to distinguish different curves in time.}}	
 \label{fig:surfactantSpreading-drop-ex2}	 		
 \end{figure}

Next, we study the spreading of a surfactant drop on a thin liquid film by considering the full PDEs, adding the contributions corresponding to capillary, gravity and surface diffusion terms by setting $\mathscr{C}=0.013$, $\mathscr{G}=20.846$, $\mathrm{Pe}=10^5/3$, while assuming a smooth solid substrate $f(x,y)=0$. The computational domain is $[-8,8]\times[-8,8]$ which is discretized using $256\times256$ cubic B-spline elements. Fixed time step sizes of $\Delta t = 0.01$ is used. The  initial conditions are:
\begin{align*}
& h(x,y,0)=1, \qquad
 c(x,y,0)= \frac{1}{2} (1-\tanh(4(r-1))) , \qquad  \text{where } r=\sqrt{x^{2}+y^{2}}.
\end{align*}
The results which show the evolution of the surfactant concentration and the film height for $t=[0,50]$ in time steps of $2.5$ are presented in \fref{fig:surfactantSpreading-drop-ex2}. The initial steep gradient in surfactant concentration induces surface tension gradients which drives a Marangoni flow. According to this flow, the surfactants are spread in the direction of surfactant-free regions of the thin film surface. This surfactant spreading relaxes the concentration gradient and accordingly the driving force of the Marangoni flow becomes attenuated. During this spreading process, the liquid film undergoes large deformations; the film height reaches its maximum almost at the leading edge of surfactant monolayer. Due to the mass conservation, this elevation in film height is compensated by a substantial decrease in film height at the upstream. 
We note that this problem was already studied in \cite{Wong2011,LIU2019429} using finite volume/finite difference and mixed finite element methods, respectively. The presented isogeometric results are in perfect agreement with those from \cite{Wong2011,LIU2019429}.  
\subsection{Verification example: The effect of finite--amplitude perturbations on the surfactant spreading}
\begin{figure}
    	\centering
    	\begin{subfigure}{0.48\textwidth}
    	\centering
    	\includegraphics[width=\linewidth,trim={0.5cm 0cm 1cm 3cm},clip]{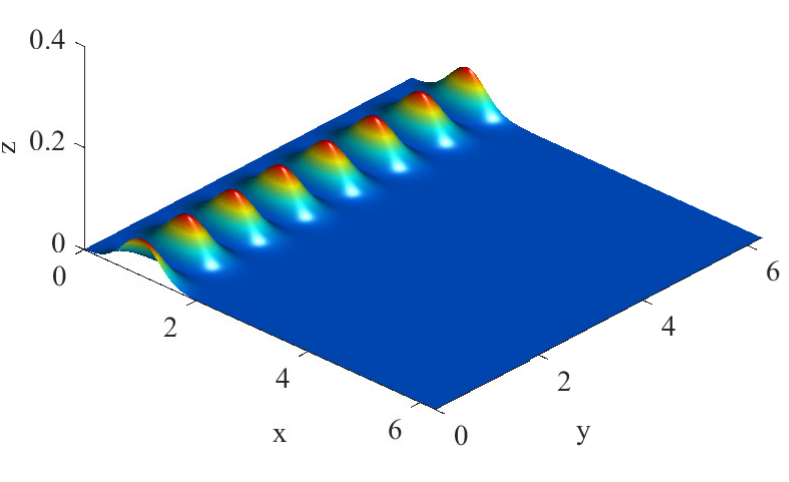}
    	\caption{}
    	\end{subfigure}
    	\quad
    	\begin{subfigure}{0.48\textwidth}
    	\centering
    		\includegraphics[width=\linewidth,trim={0.5cm 0cm 1cm 3cm},clip]{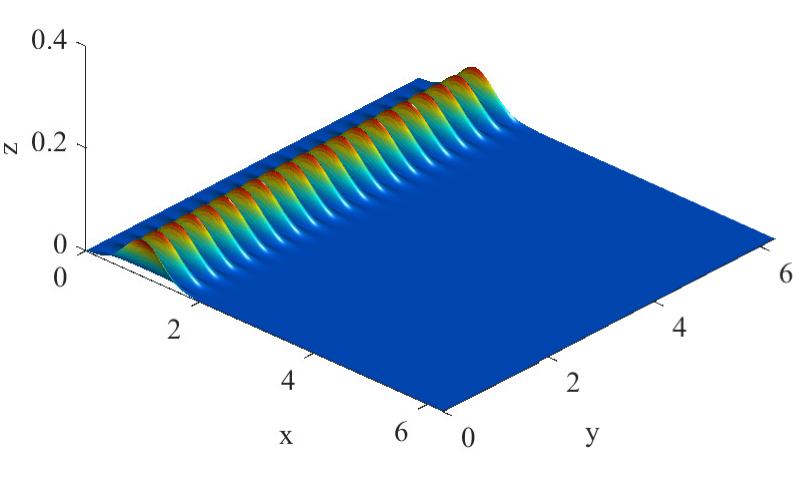}
    		\caption{}
       \end{subfigure}
       \begin{subfigure}{0.45\textwidth}
           	\centering
           		\includegraphics[width=\linewidth,trim={0.5cm 0cm 1cm 2cm},clip]{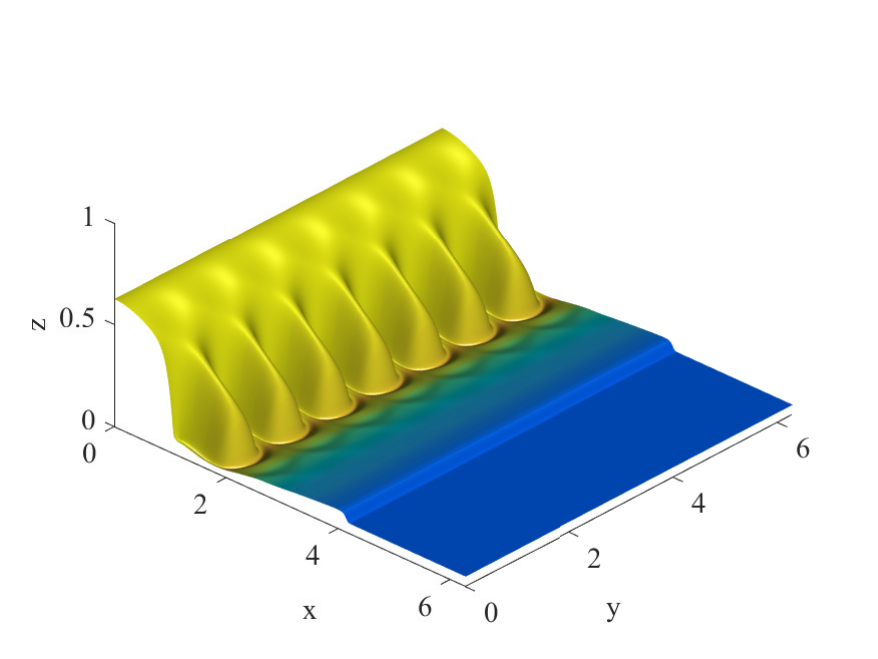}
           		\caption{}
              \end{subfigure}
              \quad
       \begin{subfigure}{0.45\textwidth}
           	\centering
           		\includegraphics[width=\linewidth,trim={0.5cm 0cm 1cm 2cm},clip]{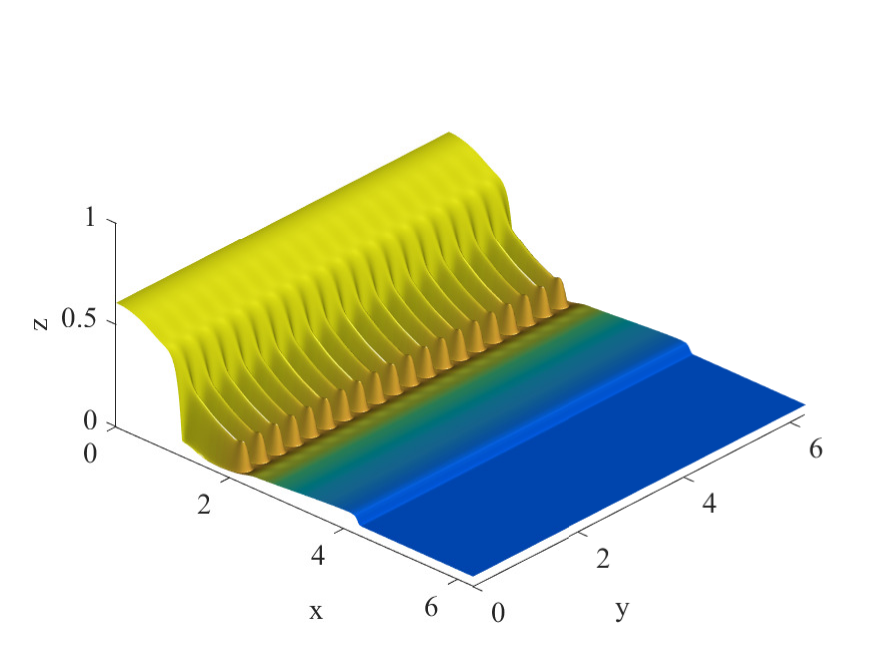}
           		\caption{}
              \end{subfigure}       
       \begin{subfigure}{0.05\textwidth}
                  	\centering
                  		\includegraphics[width=\linewidth,trim={0cm 0cm 0cm 0cm},clip]{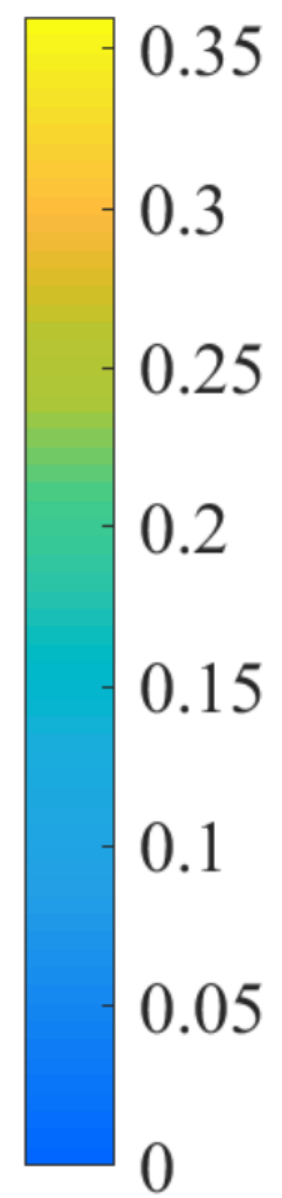}
                  		\caption*{$c$}
                     \end{subfigure}         	
    \caption{\bl{(a) and (b): Surface plots of the roughness function for $\lambda=7$ and $\lambda=20$, respectively. (c) and (d): Surface plots of the film height at $t = 100$ for $\lambda=7$ and $\lambda=20$, respectively, colored by the surfactant concentration. Uniform meshes with $512\times512$ elements are employed.}}
    \label{fig:fingeringInstability-ex00}	 		
    \end{figure}

Introducing perturbations to the surfactant spreading problem is essential for triggering fingering instabilities. It has been shown that  finite--amplitude perturbations, when localized at the edge of the surfactant deposition and applied either to the initial film height or to the roughness function of the substrate, can significantly affect the surfactant spreading and the deformation of liquid film \cite{warner_craster_matar_2004,LIU2019429}. In fact, defining perturbations by trigonometric functions with a single wavenumber would lead to the formation of several uniform fingers, the number of which is equal to the wavenumber (cf. \cite{LIU2019429}). We use this feature to further verify our code. Following \cite{LIU2019429}, we consider two cases and in each case we perturb the roughness function of the substrate by considering a trigonometric function with a single mode. The computational domain is $[0,2\pi]\times[0,2\pi]$ and is discretized by $512\times512$ cubic B-spline elements. The roughness function and the initial conditions are chosen as the following
\begin{equation*}
\begin{split}
& f(x,y)=\bar{A}\exp\left(-B(x-1)^{2}\right) \left( \cos(\lambda y)+1\right), \\
& h(x,y,0)=(1-x^{2}+b) H(1-x)+b H(x-1) +\bar{A}\exp \left(-B(x-1)^{2}\right), \\
& c(x,y,0)=H(1-x),\\
& H(x)=\frac{1}{2}\left(1 + \tanh \left( K x\right)\right).
\end{split}
\end{equation*}
where $H(x)$ is a generalized Heaviside function and the parameters are taken as $\bar{A}=0.035$, $b=0.05$, $K=20$, $B=5$, $\mathscr{C}=10^{-4}$, $\mathscr{G}=0$, $\mathrm{Pe}=10^{4}$. \bl{Here, $\bar{A}$ controls the magnitude of the substrate roughness which in this case equals $2\bar{A}$. $b=\frac{\mathcal{H}_b}{\mathcal{H}}$ is a geometrical parameter which denotes the ratio between the thickness of the precursor film to the characteristic height of the film; when $\mathcal{H}\approx1$, $b$ represents the thickness of the precursor film. $K$ and $B$ are parameters that control the smoothing lengths of the smoothed-out Heaviside and delta functions, respectively.\\}
For the first case, we assume the wavenumber to be $\lambda = 7$ while $\lambda = 20$ for the second case. We solve the problem until $t=100$. \fref{fig:fingeringInstability-ex00} presents the surface plots of the roughness function for $\lambda = 7$ and $\lambda = 20$ together with the corresponding surface plots of the film height, colored by the \bl{surfactant} concentration, at the final time. As expected, the number of formed fingers equals the chosen wavenumber in each case. 
    
\subsection{Fingering instability}
In this section, we consider the fingering instability that may develop behind the spreading front when the surfactant is spreading on prewetted substrates \cite{Afsar-Siddiqui2003a,Afsar-Siddiqui2003b}. The physical mechanism of this Marangoni-driven fingering instability in the absence of long-range
intermolecular forces (van der Waals forces) was first conceptualized in \cite{warner_craster_matar_2004} (see also Sec. VI.B in \cite{Craster2009}) and is explained in the following. Consider there are some local height increases at the edge of a drop of surfactant deposited on a thin liquid film. These local height increases in the thinning region yields an increase in the local surface velocity. As a result, the amount of surfactant transported from this region through advection will be increased which results in a local decrease of surfactant concentration. The local decrease in surfactant concentration creates a Marangoni flow from the neighboring depressed surfactant-rich regions to the locally elevated regions, which brings more surfactant and flow to these regions. Therefore, at the drop edge, relatively thick and thin liquid regions are juxtaposed. The faster spreading of thick regions with respect to the neighboring thinned ones forms finger-like patterns.\\
As discussed earlier, to model fingering instability numerically, it is customary to introduce perturbations in the initial height of the film \cite{warner_craster_matar_2004} which then trigger the instability. Recently \cite{LIU2019429} showed that the roughness of the substrate can also initiate fingering instability. Here, we model both height--induced and roughness--induced fingering instabilities by considering several examples from \cite{warner_craster_matar_2004,LIU2019429}. In particular, we intentionally focus on the numerical examples provided in \cite{LIU2019429} in order to demonstrate that IGA is a promising alternative to the mixed finite element methods already available in the literature. Henceforth, we assume $\mathscr{G}=0$ neglecting gravitational effects that is a valid assumption for very thin liquid films. Consistent with this assumption and similar to \cite{LIU2019429}, we assume a small capillarity number $\mathscr{C}=10^{-4}$ which can be the case for very thin liquid films. The surface Peclet number is set to $\mathrm{Pe}=10^{4}$.
\subsubsection{The problem of a planar strip of surfactant} \label{sec:planarStripOfSurfactant}
\begin{figure}
 	\centering
 	\begin{subfigure}{0.48\textwidth}
 	\centering
 	\includegraphics[width=\linewidth,trim={10cm 4cm 4cm 4cm},clip]{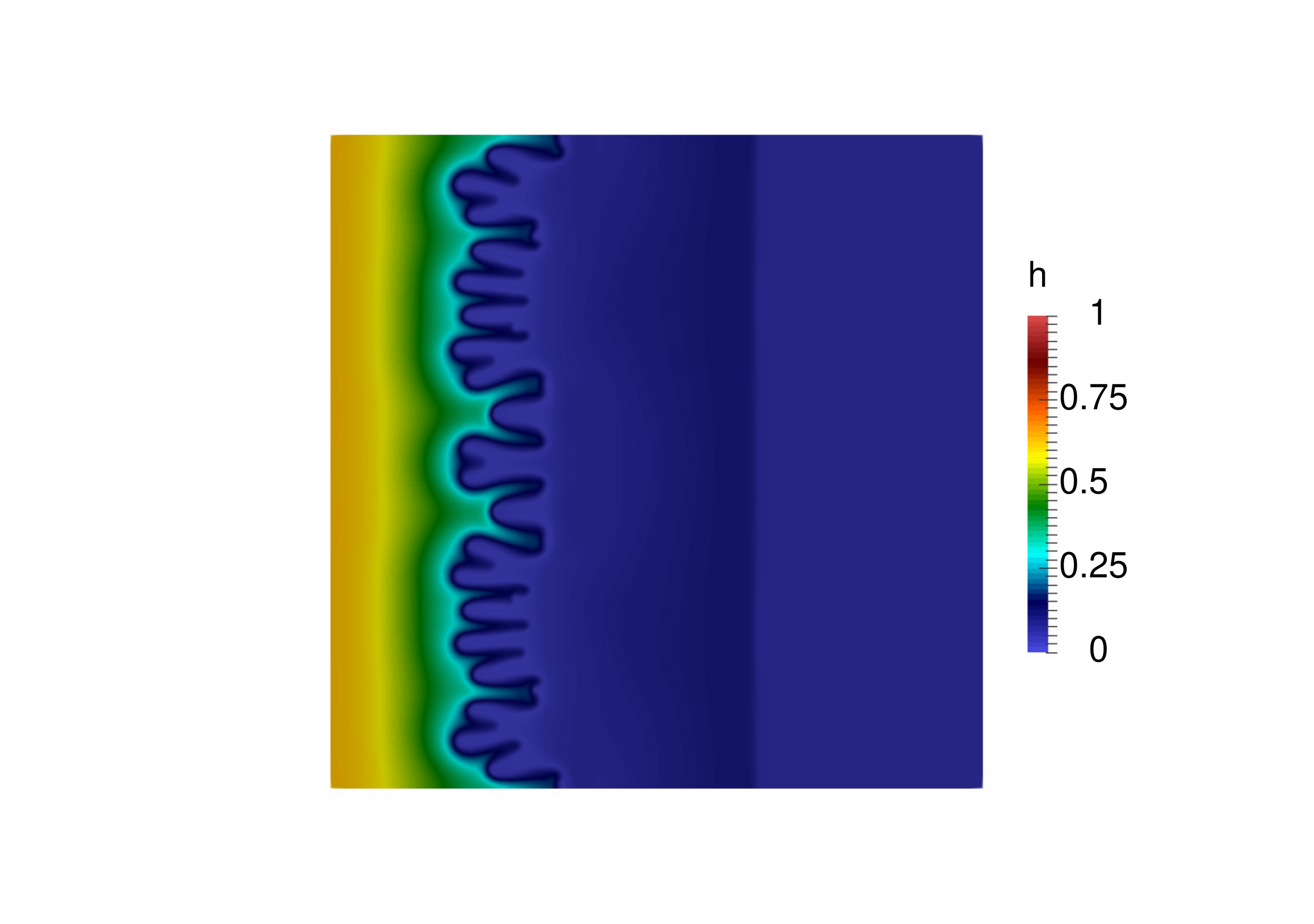}
 	\caption{}
 	\end{subfigure}
 	\quad
 	\begin{subfigure}{0.48\textwidth}
 	\centering
 		\includegraphics[width=\linewidth,trim={10cm 4cm 4cm 4cm},clip]{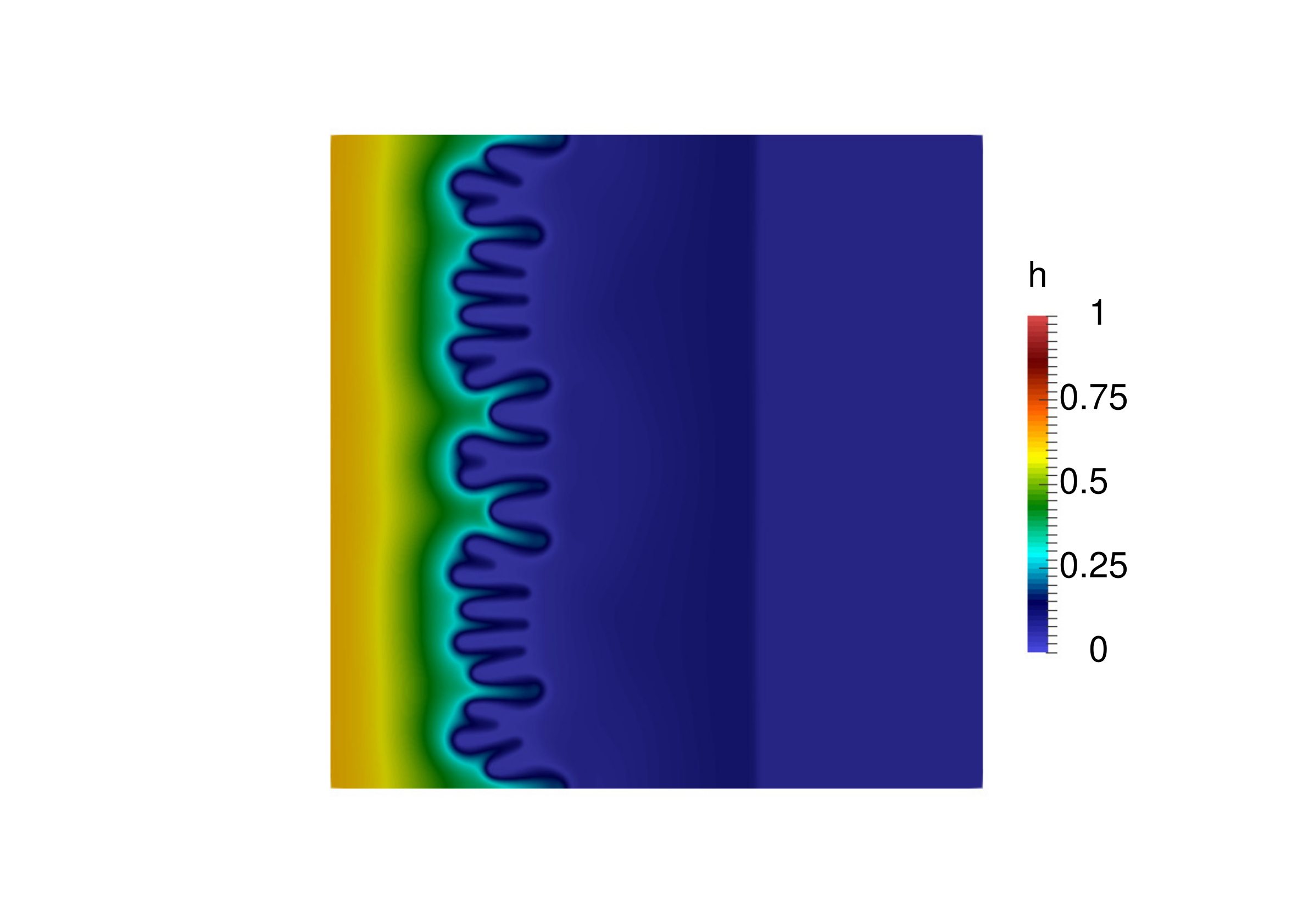}
 		\caption{}
    \end{subfigure}	    
 \caption{\bl{Isogeometric solution for fingering instability initiated by film height perturbations. The plots show the height of a thin liquid film rested on a smooth solid substrate at $t=100$. Uniform meshes with (a) $256\times256$ and (b) $512\times512$ elements are used.}}	
 \label{fig:fingeringInstability-ex01}	 		
 \end{figure}

A rectangular strip of surfactant deposited on the left edge of a square computational domain $[0,2\pi]\times[0,2\pi]$ is considered. \bl{The computational domain is assumed to have periodic boundary conditions on its top and bottom surfaces}. We first model this problem on a smooth substrate by setting $f(x,y)=0$ and initiate the instability solely by perturbing the initial film height. Therefore, the initial conditions are set to    
\begin{equation*}
\begin{split}
& h(x,y,0)=(1-x^{2}+b) H(1-x)+b H(x-1) +\exp \left(-B(x-1)^{2}\right) \sum_{k=1}^{N} A_{k} \cos(\lambda_{k}y), \\
& c(x,y,0)=H(1-x).
\end{split}
\end{equation*}
where the parameters are taken as $b = 0.05, K=20, B = 5, A_1=0.01, A_2= 0.01, A_3=0.01, A_4=0.005, \lambda_1 = 2, \lambda_2 = 5, \lambda_3 = 7, \lambda_4 = 20$. Note that the last term in the definition of initial height is a finite--amplitude perturbation localized at the edge of surfactant deposition. We model the problem using $256\times256$ and $512\times512$ cubic B-spline elements. The film height computed at $t=100$ for these two different discretizations are presented in \fref{fig:fingeringInstability-ex01} which clearly shows the fingering instability and the ability of our formulation to accurately capture all the main processes including spreading of fingers, finger tip splitting, and shielding. Despite some fine details at the finger tips, the results for these two meshes are almost visually indistinguishable; this confirms that it is indeed the initial height perturbations that caused the instability and not some numerical errors. The observed fingering patterns are nearly the same as those presented in \cite{LIU2019429} using a mixed finite element method, and in \cite{Wong2011} using a finite volume method.  

\begin{figure}
  	\centering
  	\begin{subfigure}{0.48\textwidth}
  	\centering
  	\includegraphics[width=\linewidth,trim={10cm 3cm 4cm 3cm},clip]{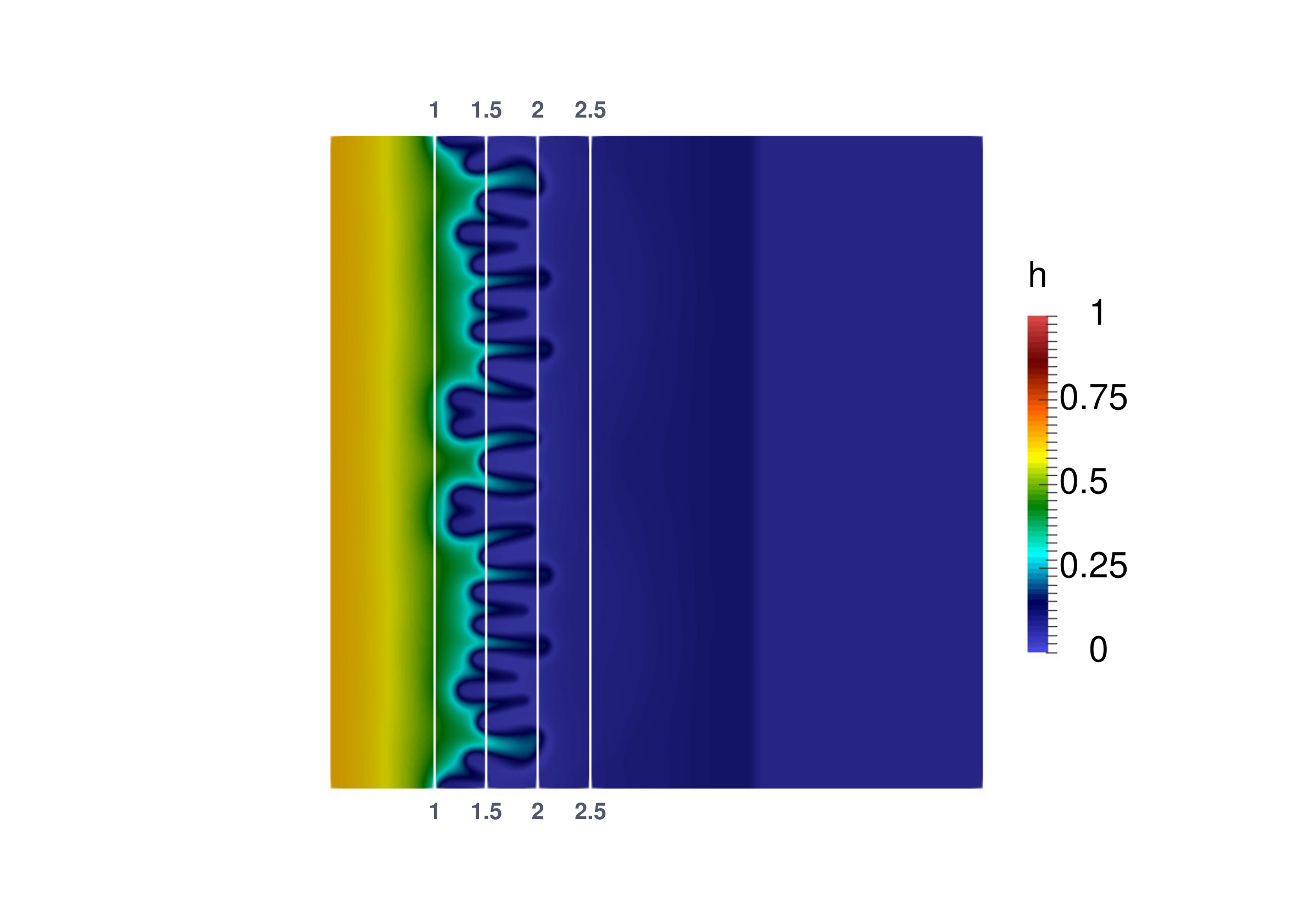}
  	\caption{}
  	\end{subfigure}
  	\quad
  	\begin{subfigure}{0.48\textwidth}
  	\centering
  		\includegraphics[width=\linewidth,trim={10cm 3cm 4cm 3cm},clip]{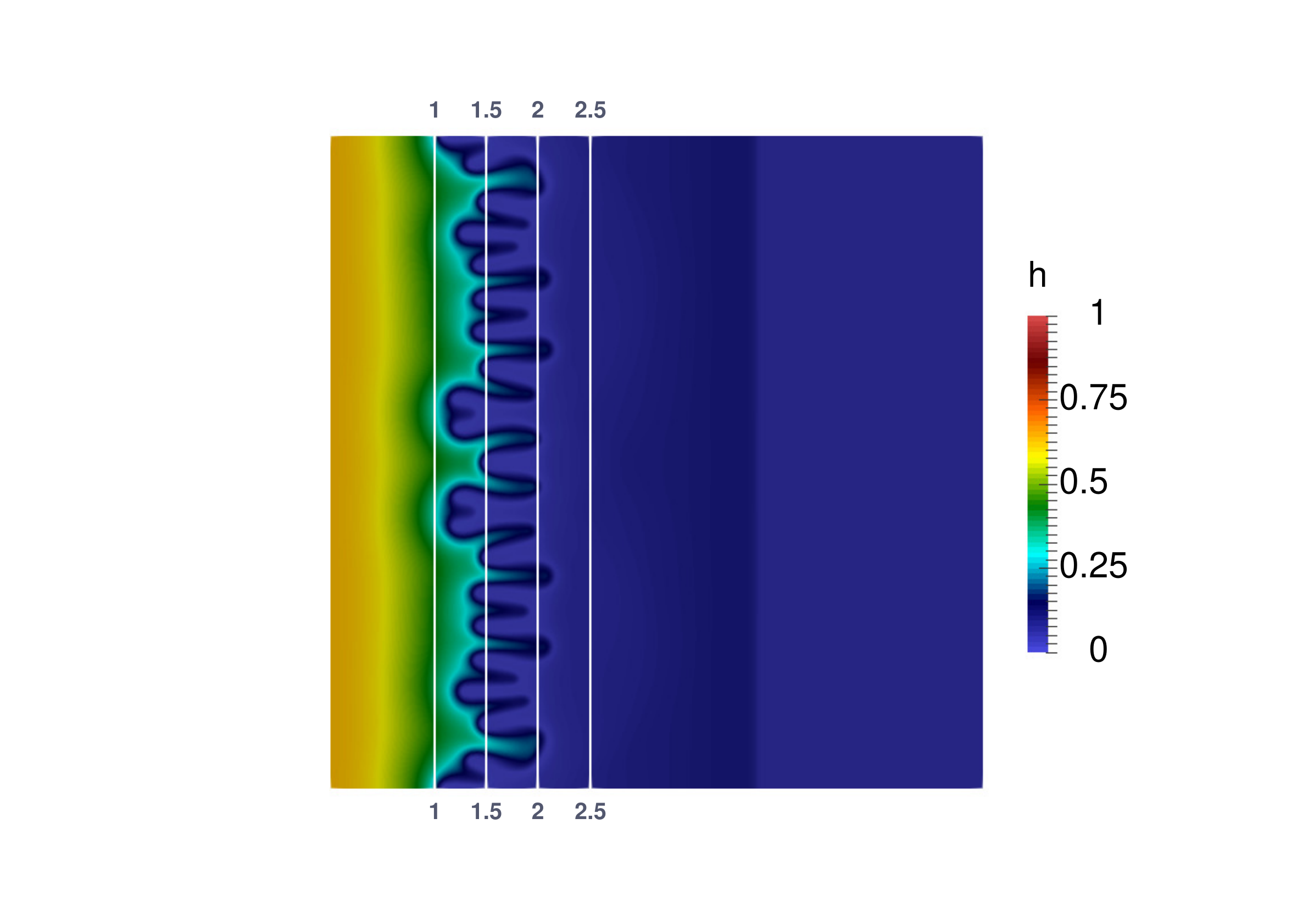}
  		\caption{}
     \end{subfigure}	
  \caption{\bl{Isogeometric solution for fingering instability initiated by: (a) roughness of the substrate (b) initial height perturbations. The initial physical height $h_p$ and concentration are the same in both cases. The plots show the height of the thin film at $t=100$. Uniform meshes with $512\times512$ elements are employed.}}	
  \label{fig:fingeringInstability-ex02}	 		
  \end{figure}

\begin{figure}
    	\centering
    	\begin{subfigure}{0.48\textwidth}
    	\centering
    	\includegraphics[width=\linewidth,trim={0cm 0cm 0cm 0cm},clip]{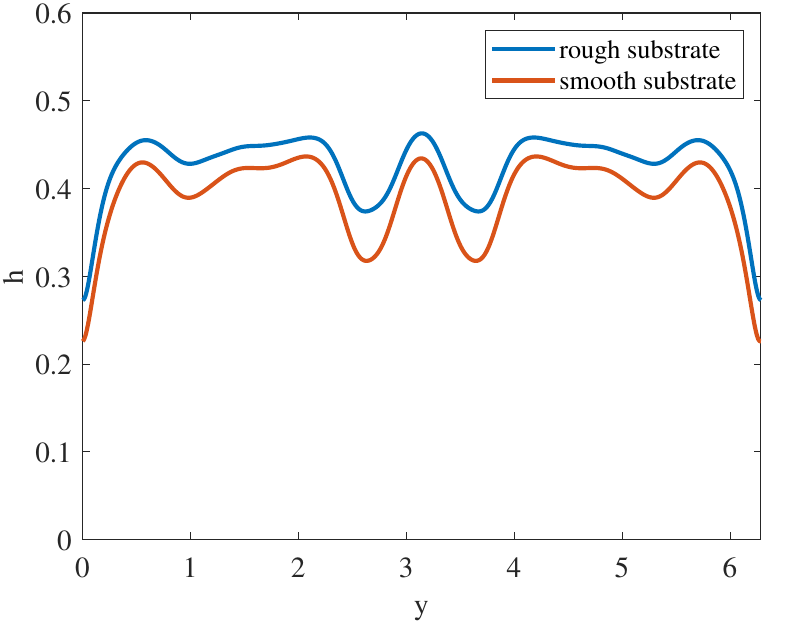}
    	\caption{x=1.0}
    	\end{subfigure}
    	\quad
    	\begin{subfigure}{0.48\textwidth}
    	\centering
    		\includegraphics[width=\linewidth,trim={0cm 0cm 0cm 0cm},clip]{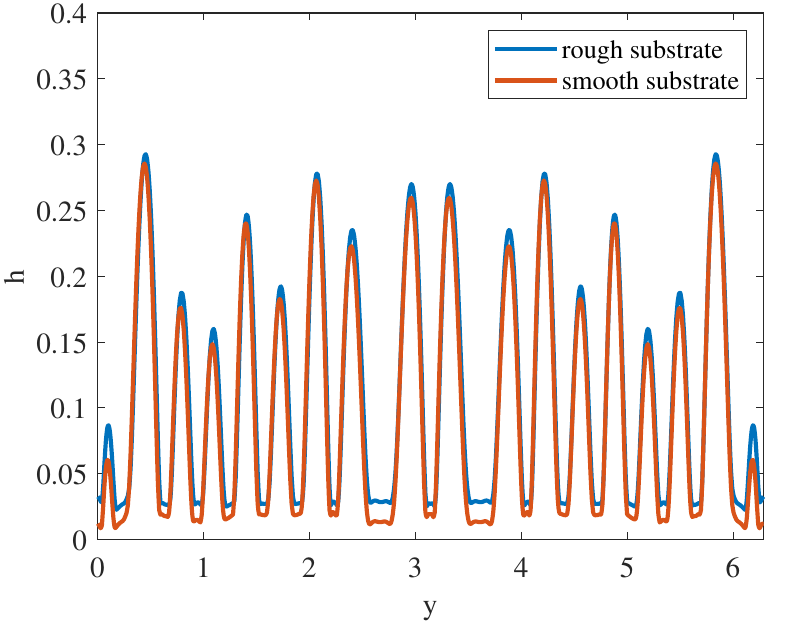}
    		\caption{x=1.5}
       \end{subfigure}
       \begin{subfigure}{0.48\textwidth}
           	\centering
           		\includegraphics[width=\linewidth,trim={0cm 0cm 0cm 0cm},clip]{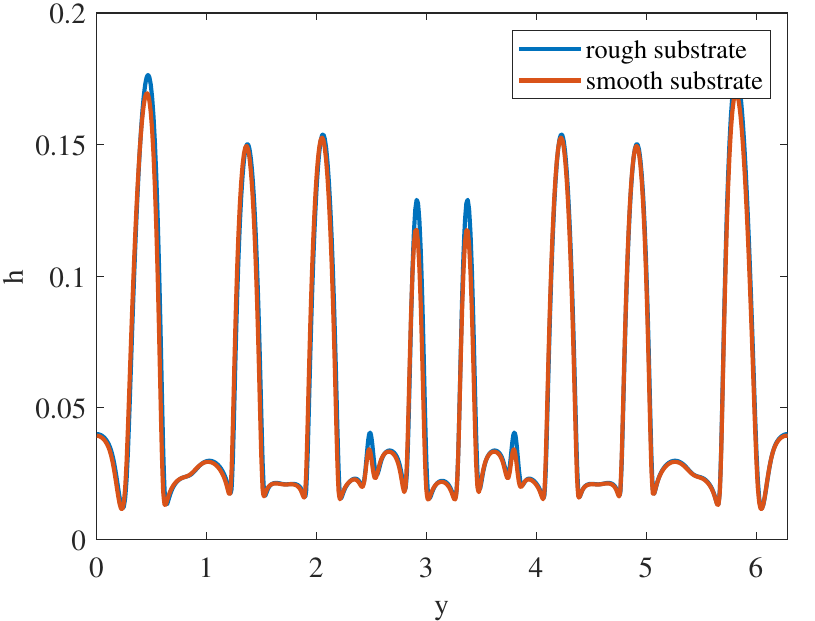}
           		\caption{x=2.0}
              \end{subfigure}
              \quad
       \begin{subfigure}{0.48\textwidth}
           	\centering
           		\includegraphics[width=\linewidth,trim={0cm 0cm 0cm 0cm},clip]{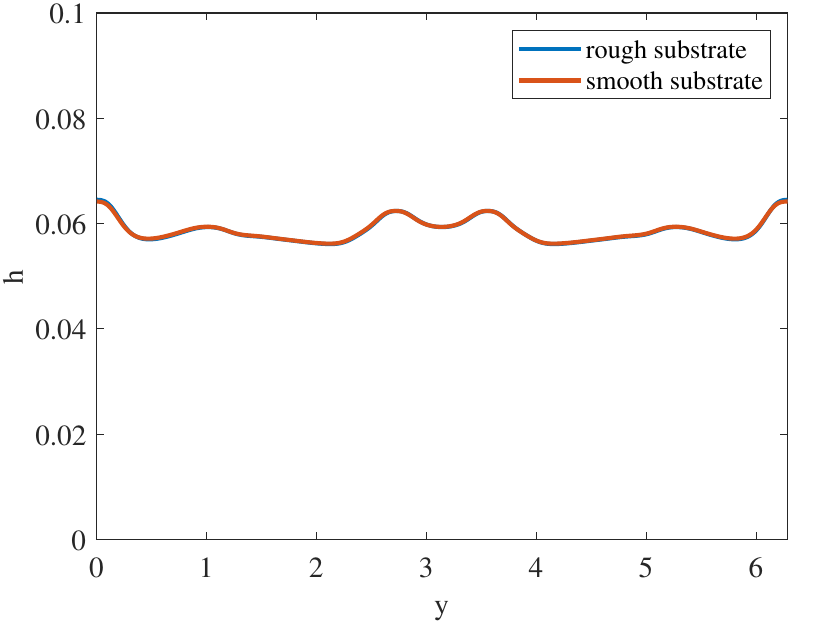}
           		\caption{x=2.5}
              \end{subfigure}       	
    \caption{\bl{Comparison of the height profiles at different cross-sections for the fingering instabilities initiated either by roughness of the substrate (rough substrate) \bl{or} by initial height perturbations (smooth substrate). The plots show the height profiles at $t=100$. Uniform meshes with $512\times512$ elements are employed.}}
    \label{fig:fingeringInstability-ex02-2}	 		
    \end{figure}

\begin{figure}
      	\centering
      	\begin{subfigure}{0.48\textwidth}
      	\centering
      	\includegraphics[width=\linewidth,trim={10cm 4cm 4cm 4cm},clip]{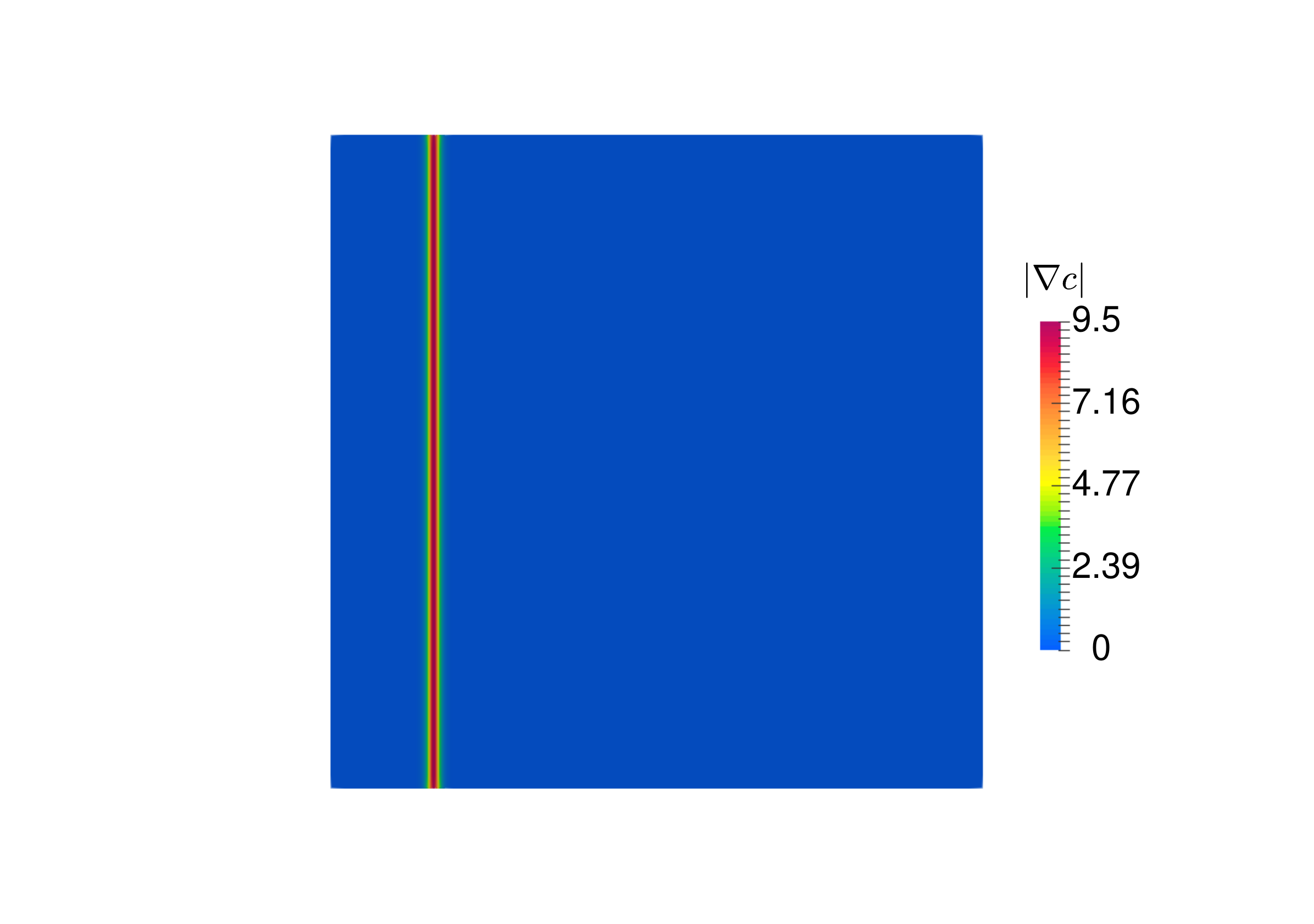}
      	\caption{$t=0$}
      	\end{subfigure}
      	\quad
      	\begin{subfigure}{0.48\textwidth}
      	\centering
      		\includegraphics[width=\linewidth,trim={10cm 4cm 4cm 4cm},clip]{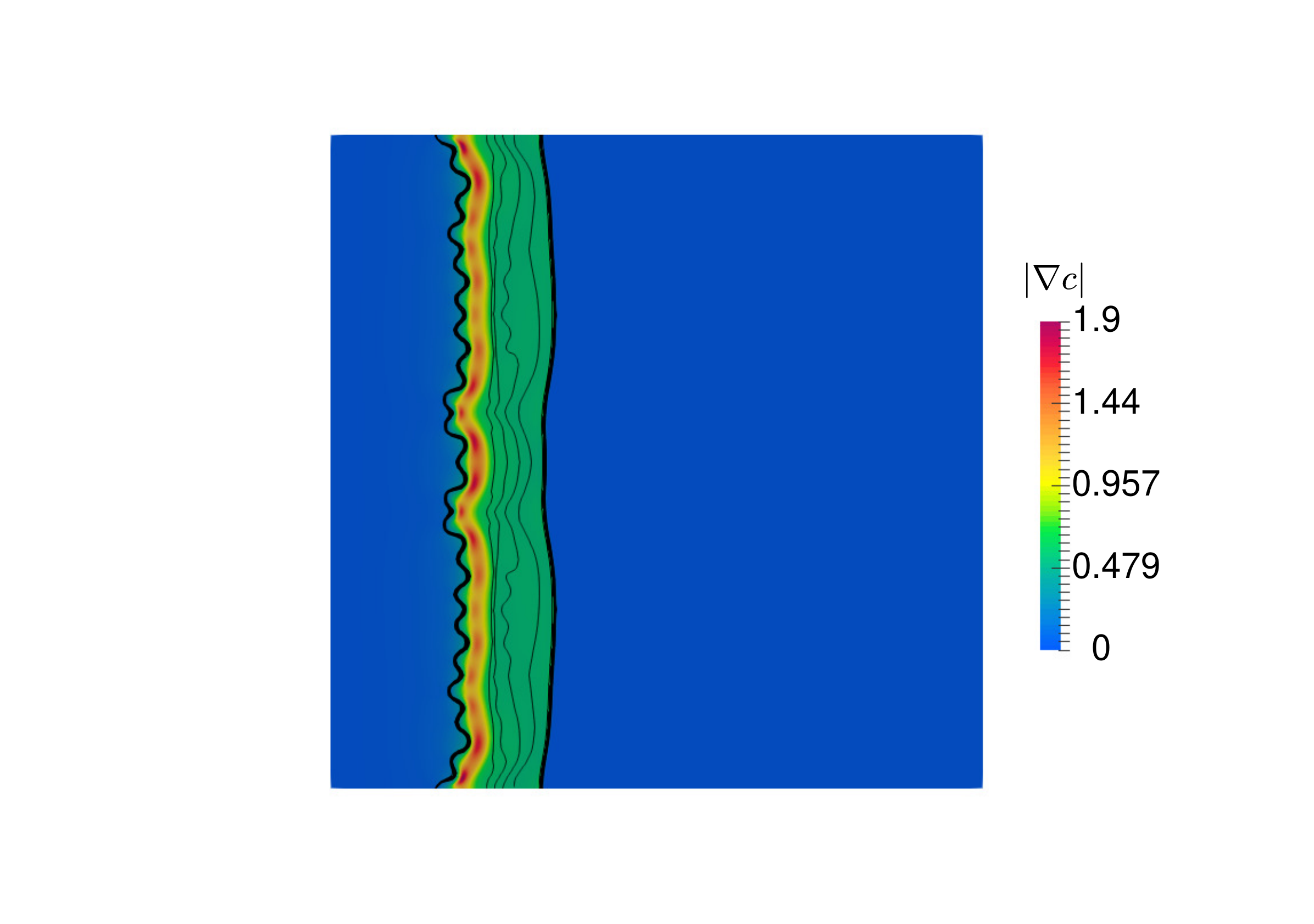}
      		\caption{$t=10$}
         \end{subfigure}	
        \begin{subfigure}{0.48\textwidth}
              	\centering
              	\includegraphics[width=\linewidth,trim={10cm 4cm 4cm 4cm},clip]{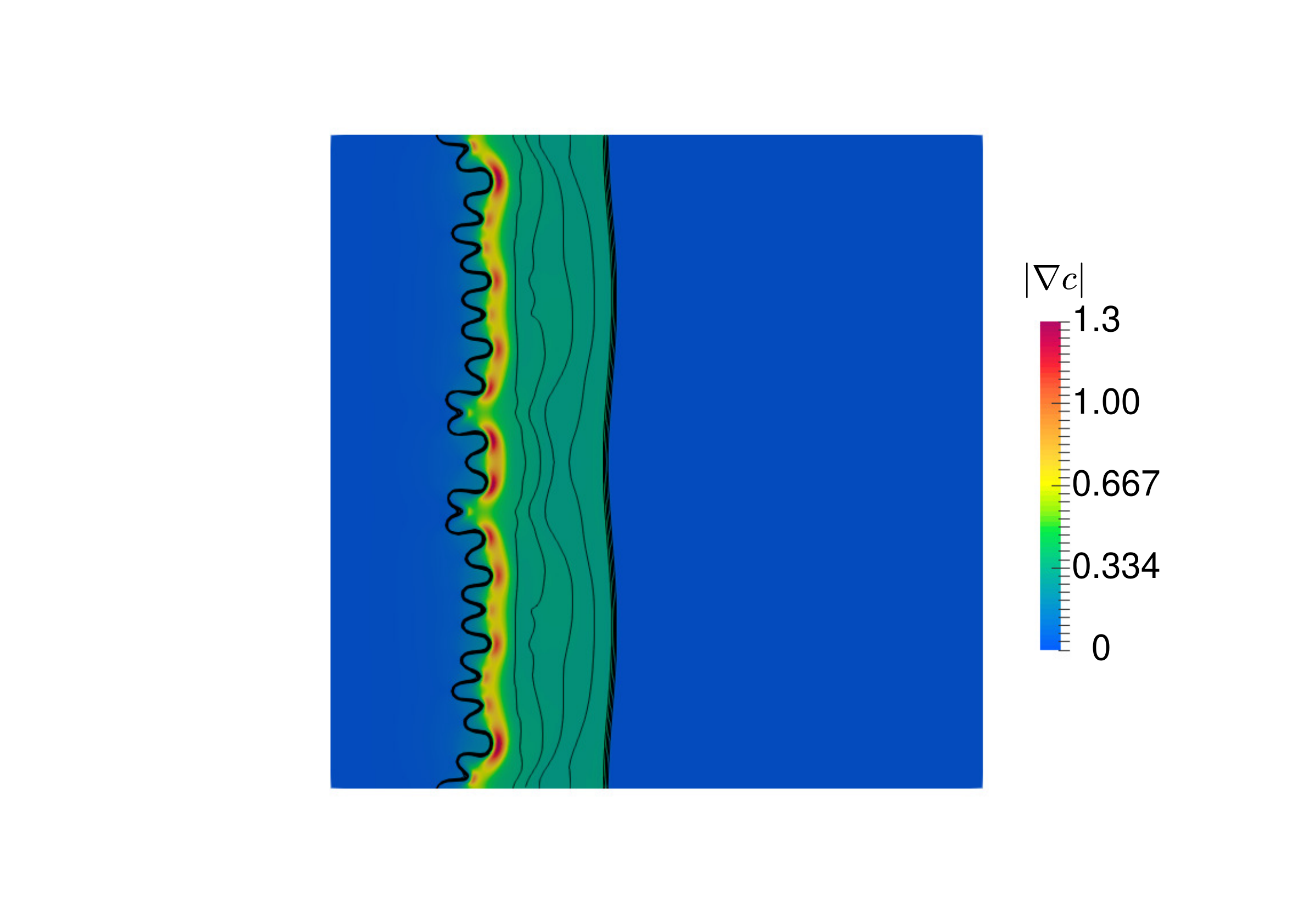}
              	\caption{$t=25$}
              	\end{subfigure}
              	\quad
              	\begin{subfigure}{0.48\textwidth}
              	\centering
              		\includegraphics[width=\linewidth,trim={10cm 4cm 4cm 4cm},clip]{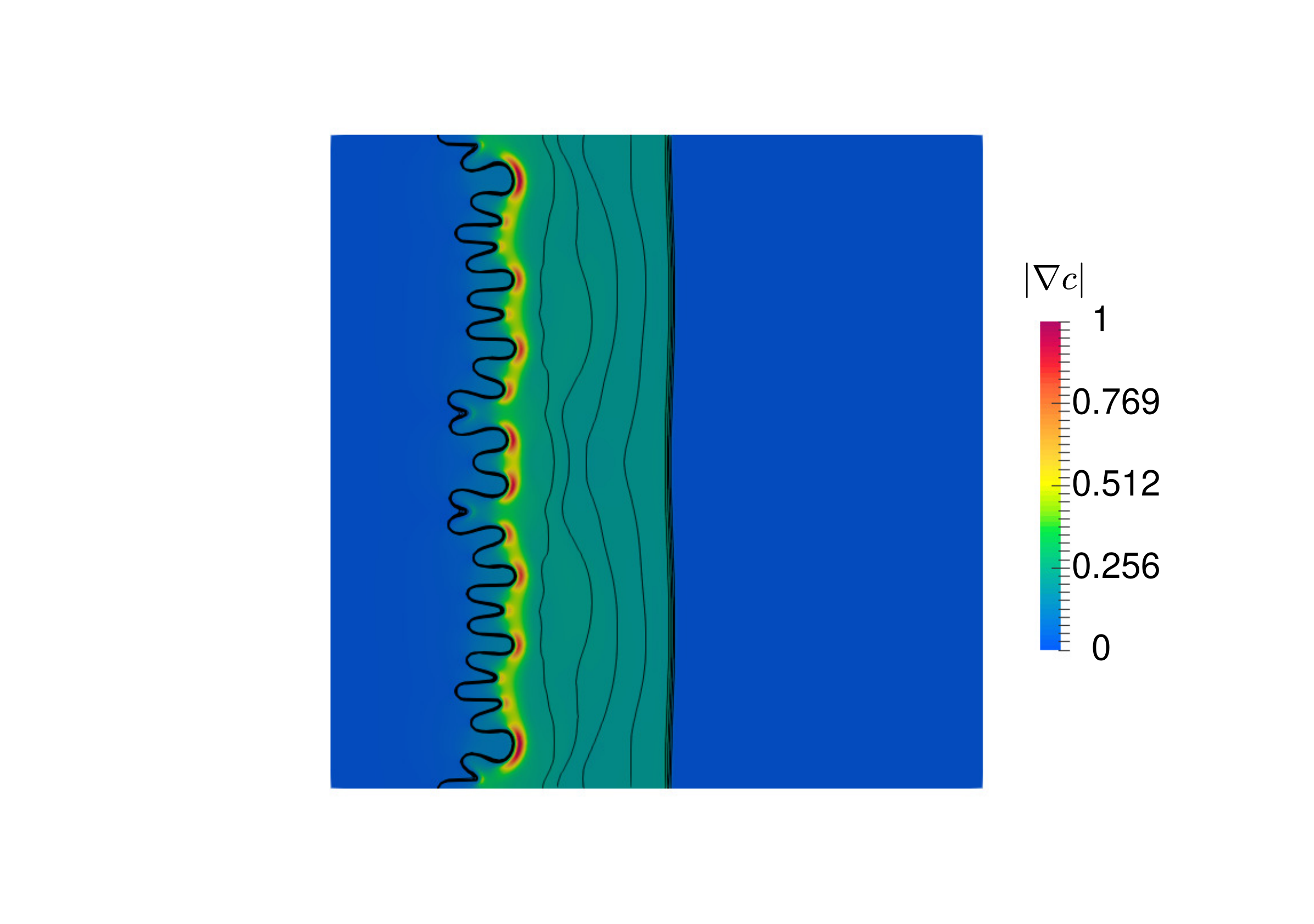}
              		\caption{$t=50$}
                 \end{subfigure}
       \begin{subfigure}{0.48\textwidth}
                     	\centering
                     	\includegraphics[width=\linewidth,trim={10cm 4cm 4cm 4cm},clip]{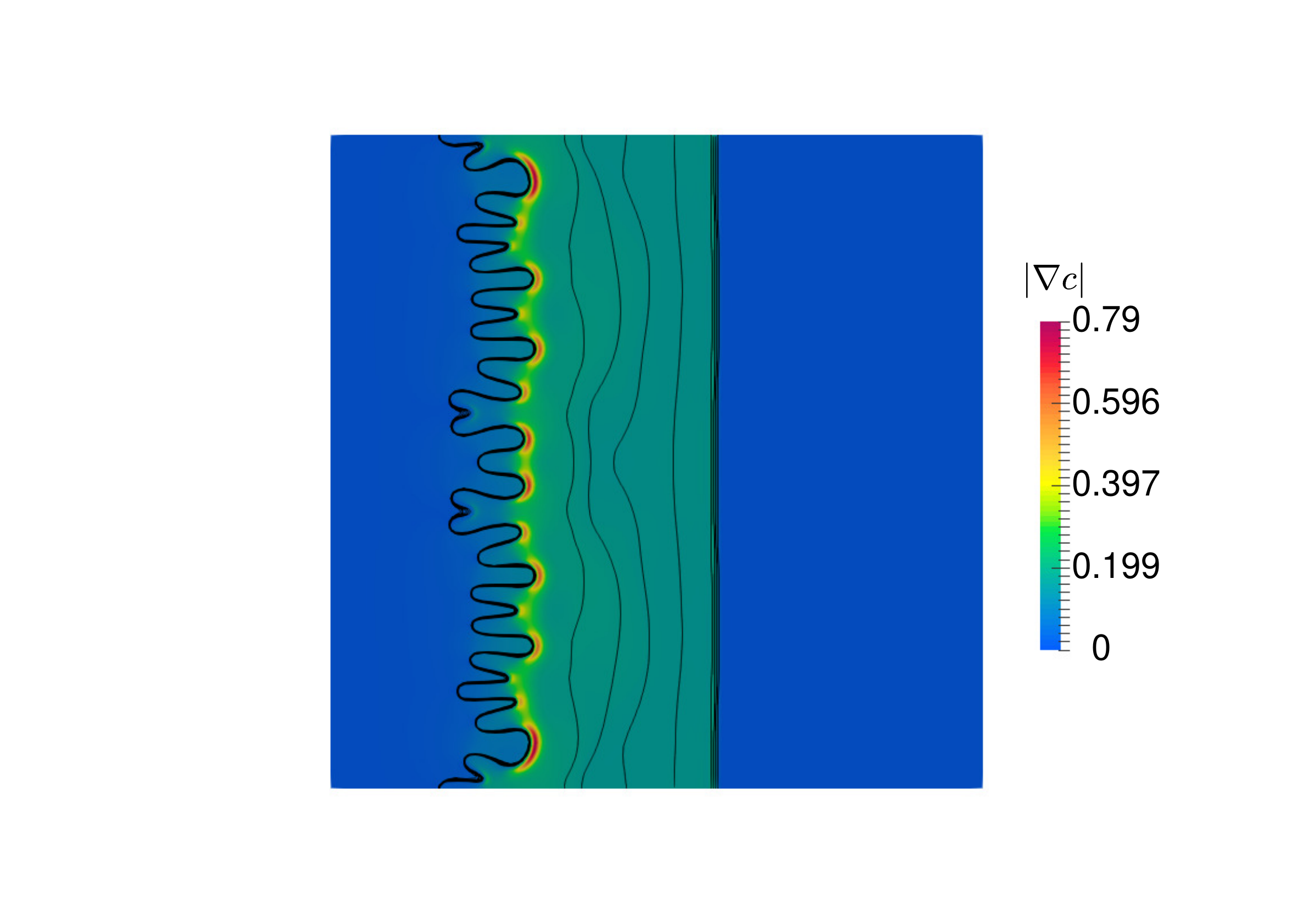}
                     	\caption{$t=74$}
                     	\end{subfigure}
                     	\quad
                     	\begin{subfigure}{0.48\textwidth}
                     	\centering
                     		\includegraphics[width=\linewidth,trim={10cm 4cm 4cm 4cm},clip]{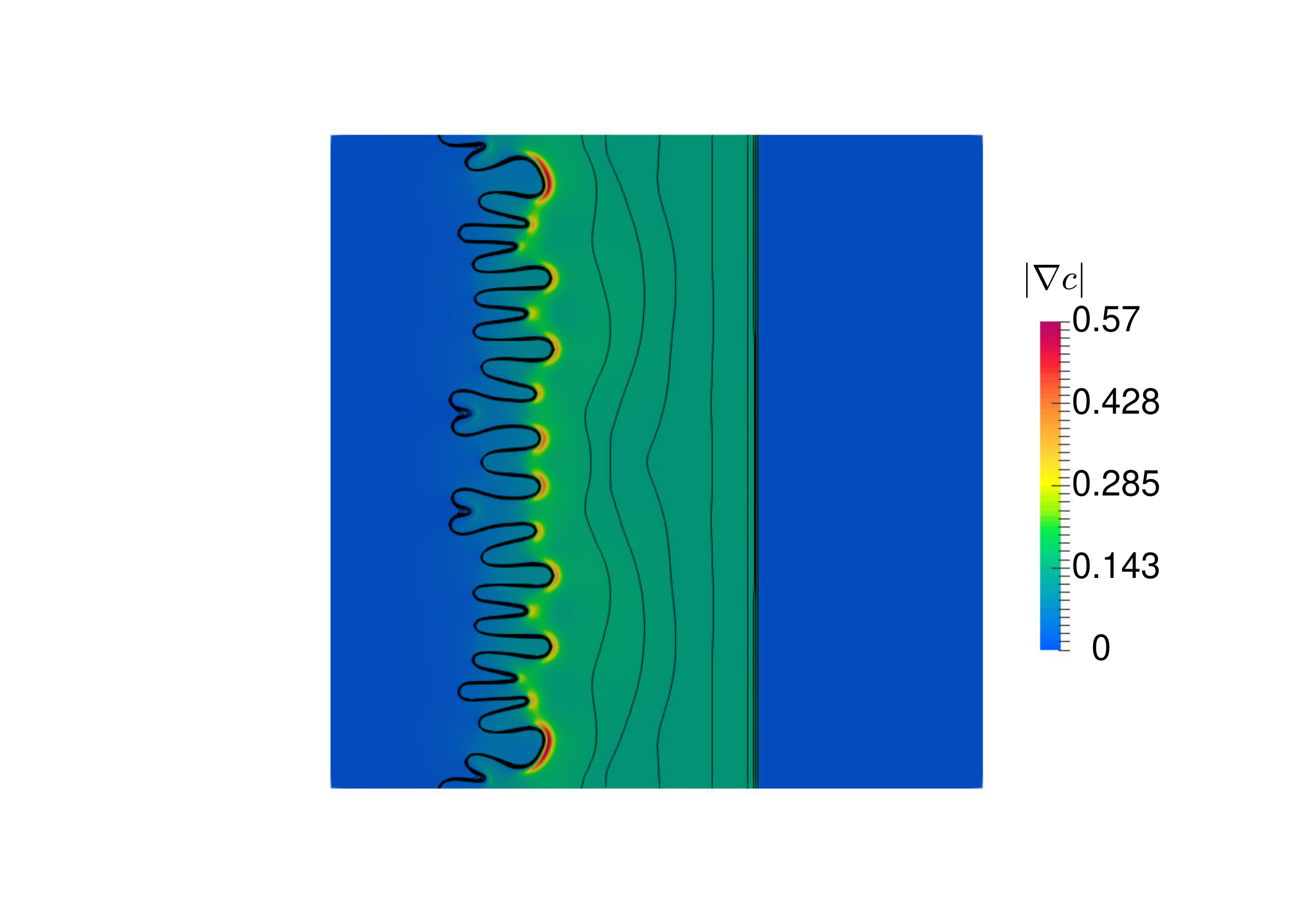}
                     		\caption{$t=100$}
                        \end{subfigure}            
      \caption{\bl{Isogeometric solution for fingering instability initiated by initial height perturbations at different time instants. The height field in the vicinity of the fingers is plotted by black contour lines. The magnitude of gradient of surfactant concentration are plotted at different time instants which show that the higher magnitude of concentration gradient at the tips of fingers is a driving mechanism for fingers to grow. Uniform meshes with $512\times512$ elements are employed.}}	
      \label{fig:fingeringInstability-ex02-3}	 		
      \end{figure}

Next, we show that for the previous example the instability can be initiated by considering a rough substrate. For this purpose, we consider the roughness function and the initial height and concentration to be defined as
\begin{equation*}
\begin{split}
& f(x,y)=\exp\left(-B(x-1)^{2}\right) \sum_{k=1}^{N} A_{k}( \cos(\lambda_{k}y)+1), \\
& h(x,y,0)=(1-x^{2}+b) H(1-x)+b H(x-1) +\bar{A}\exp \left(-B(x-1)^{2}\right), \\
& c(x,y,0)=H(1-x).
\end{split}
\end{equation*}
where $\bar{A} =\sum_{k=1}^{N} A_{k}$. All the parameters are the same as the previous example. For comparison purposes and similar to \cite{LIU2019429}, we consider a reference case with a perturbed initial height and a smooth substrate such that the initial physical height ($h_p = h-f$) and concentration in both cases are the same. The roughness function, and the initial height and concentration for this reference case are \bl{therefore} defined as follows: 
\begin{equation*}
\begin{split}
& f(x,y)=0, \\
& h(x,y,0)=(1-x^{2}+b) H(1-x)+b H(x-1) -\exp\left(-B(x-1)^{2}\right) \sum_{k=1}^{N} A_{k} \cos(\lambda_{k}y), \\
& c(x,y,0)=H(1-x).
\end{split}
\end{equation*}
The film height for these two cases, i.e. roughness--induced and perturbed height--induced fingering instabilities, are reported in \fref{fig:fingeringInstability-ex02} together with the profiles of height $h$ at different cross-sections ($x=1.0, 1.5, 2.0, 2.5$) in \fref{fig:fingeringInstability-ex02-2}. Although the results in \fref{fig:fingeringInstability-ex02} seem visually imperceptible, some differences between the two cases can be observed in \fref{fig:fingeringInstability-ex02-2}, particularly in regions well behind the thinned zone (at $x=1$). However, while not plotted here, the differences between the two cases for the profiles of the physical height $h_p$ at the same cross-sections would be minimal; this seems reasonable as in both cases we assumed the same initial physical height and concentration. Nevertheless, this is in contrast to \cite{LIU2019429} where the authors reported thicker and longer fingers over the rough substrate. \bl{In \fref{fig:fingeringInstability-ex02-3}, we plot the magnitude of the gradient of surfactant concentration, $|\nabla c|$, for the case of height--induced fingering instability at different time instants, while the height field in the vicinity of the fingers, shown by black contour lines, is superimposed on the top of the plot.} It can be observed that as time passes little fingers appear at the edge of the drop, in the thinning region. There are some fingers which at their tip $|\nabla c|$ is higher than some neighboring fingers. The Marangoni stresses are more pronounced at these fingers and they are a preferential candidate for Marangoni flow and as a result these fingers are those which keep spreading. In addition, the rate of decrease of $|\nabla c|$ is much faster at early times and becomes considerably slower at later times.    
 
    
\subsubsection{The problem of a drop of surfactant}
\begin{figure}
    	\centering
    	\begin{subfigure}{0.3\textwidth}
    	\centering
    	\includegraphics[width=\linewidth,trim={10cm 4cm 10cm 4cm},clip]{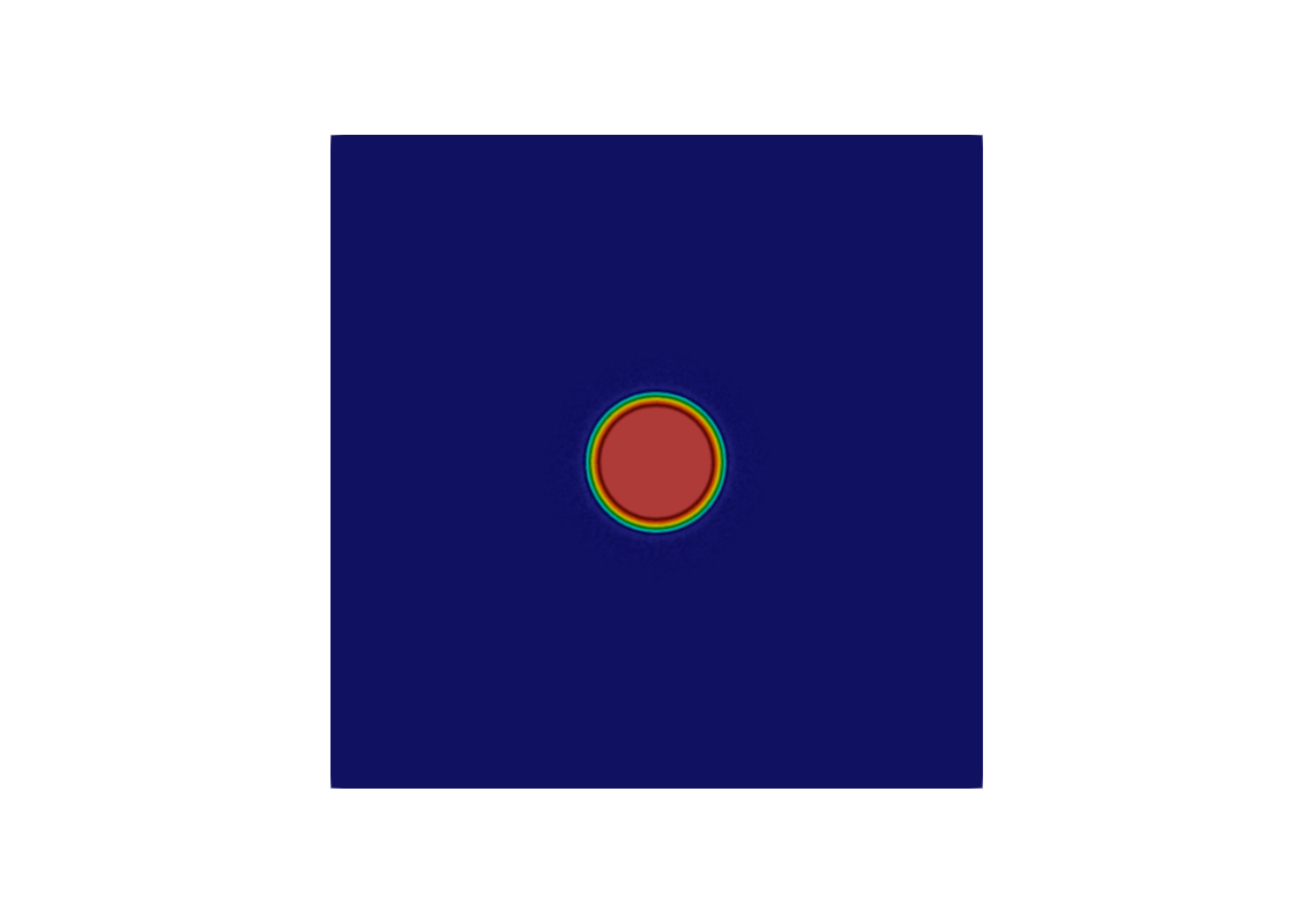}
    	\caption{t=0}
    	\end{subfigure}\quad
    	\begin{subfigure}{0.3\textwidth}
    	\centering
    		\includegraphics[width=\linewidth,trim={10cm 4cm 10cm 4cm},clip]{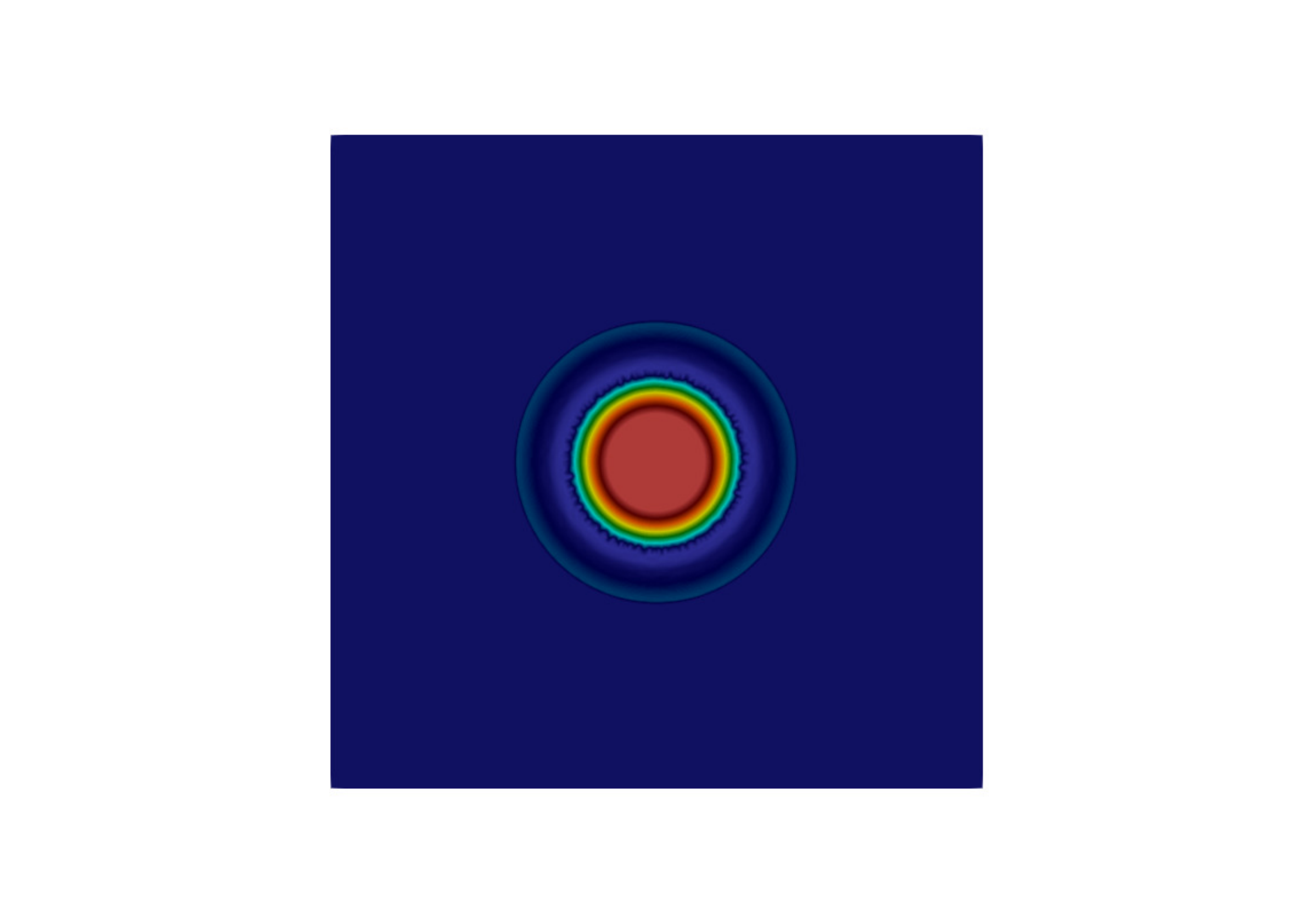}
    		\caption{t=10}
       \end{subfigure} \quad
       \begin{subfigure}{0.3\textwidth}
           	\centering
           		\includegraphics[width=\linewidth,trim={10cm 4cm 10cm 4cm},clip]{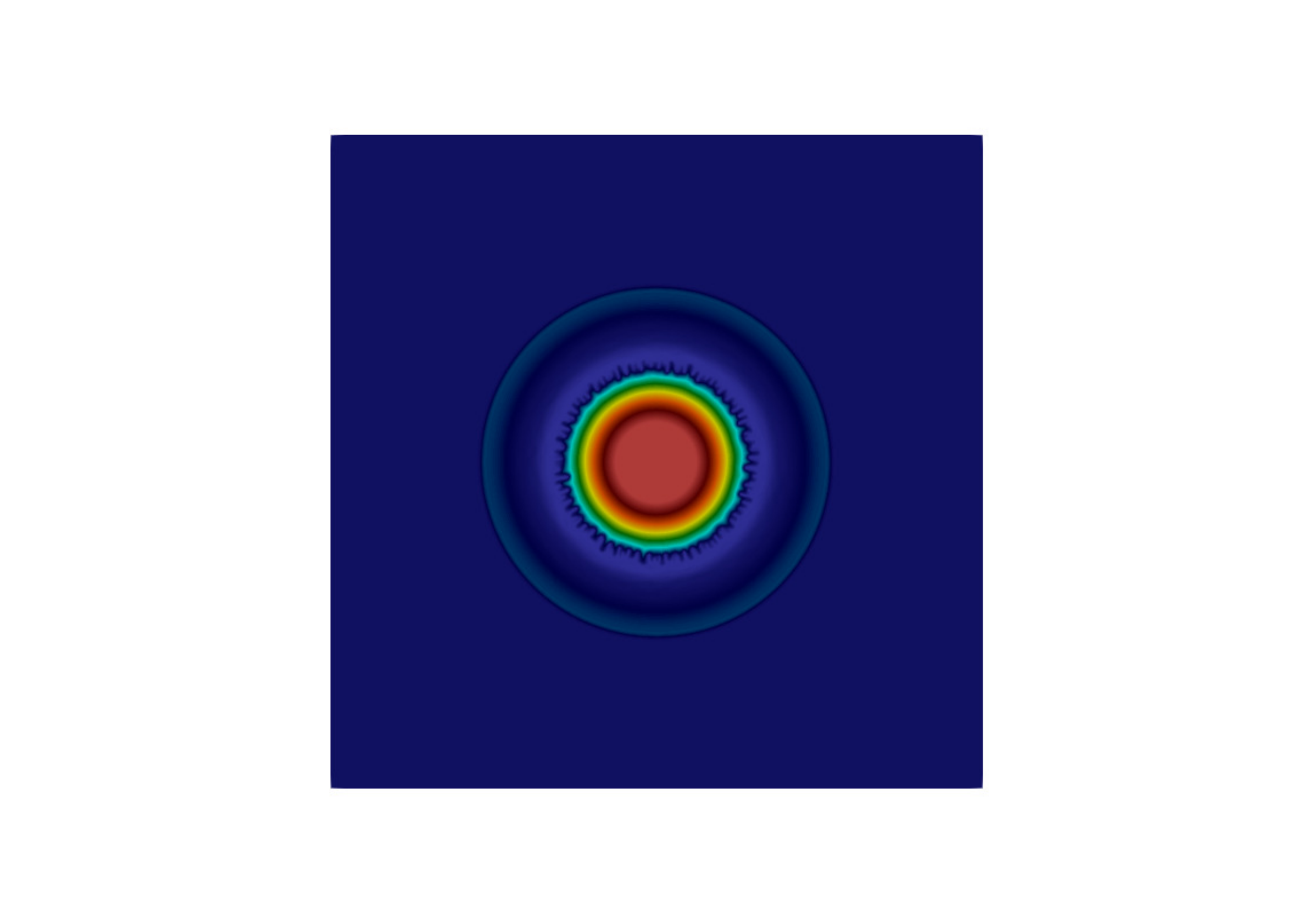}
           		\caption{t=30}
              \end{subfigure} \quad
        \begin{subfigure}{0.3\textwidth}
                   	\centering
                   		\includegraphics[width=\linewidth,trim={10cm 4cm 10cm 4cm},clip]{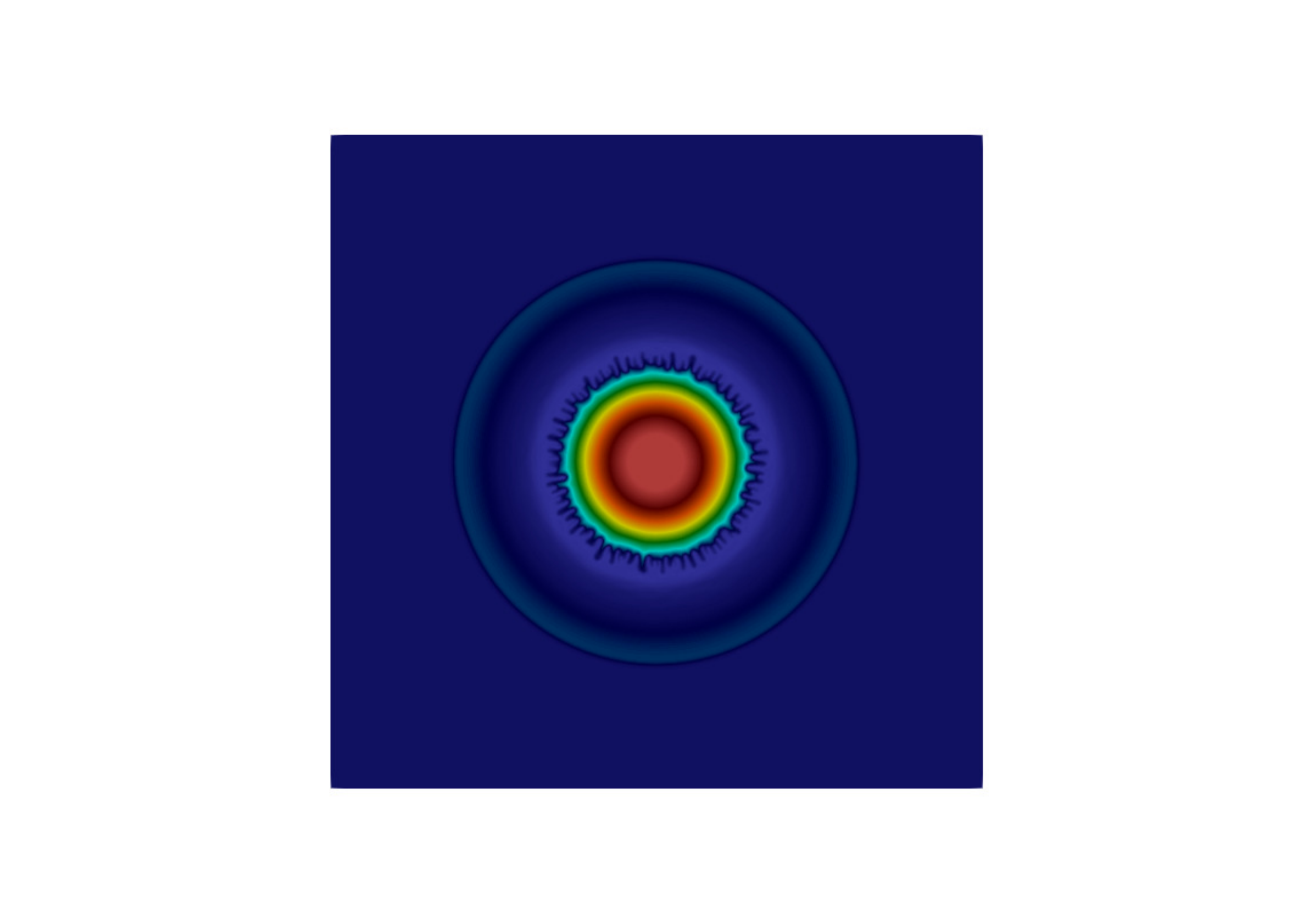}
                   		\caption{t=59}
                      \end{subfigure} \quad      
       \begin{subfigure}{0.3\textwidth}
           	\centering
           		\includegraphics[width=\linewidth,trim={10cm 4cm 10cm 4cm},clip]{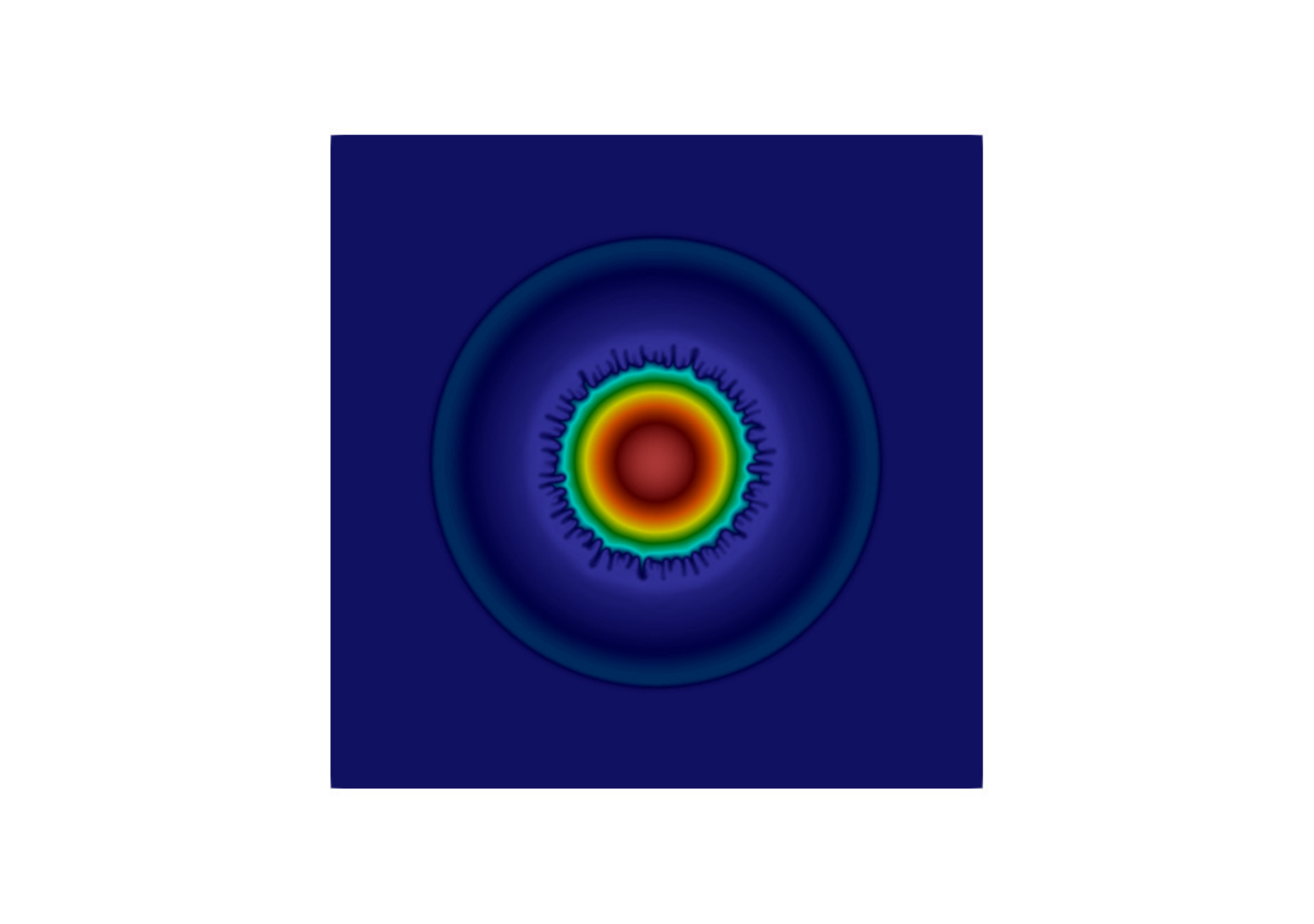}
           		\caption{t=93}
              \end{subfigure} \quad
     \begin{subfigure}{0.3\textwidth}
                	\centering
                		\includegraphics[width=\linewidth,trim={10cm 4cm 10cm 4cm},clip]{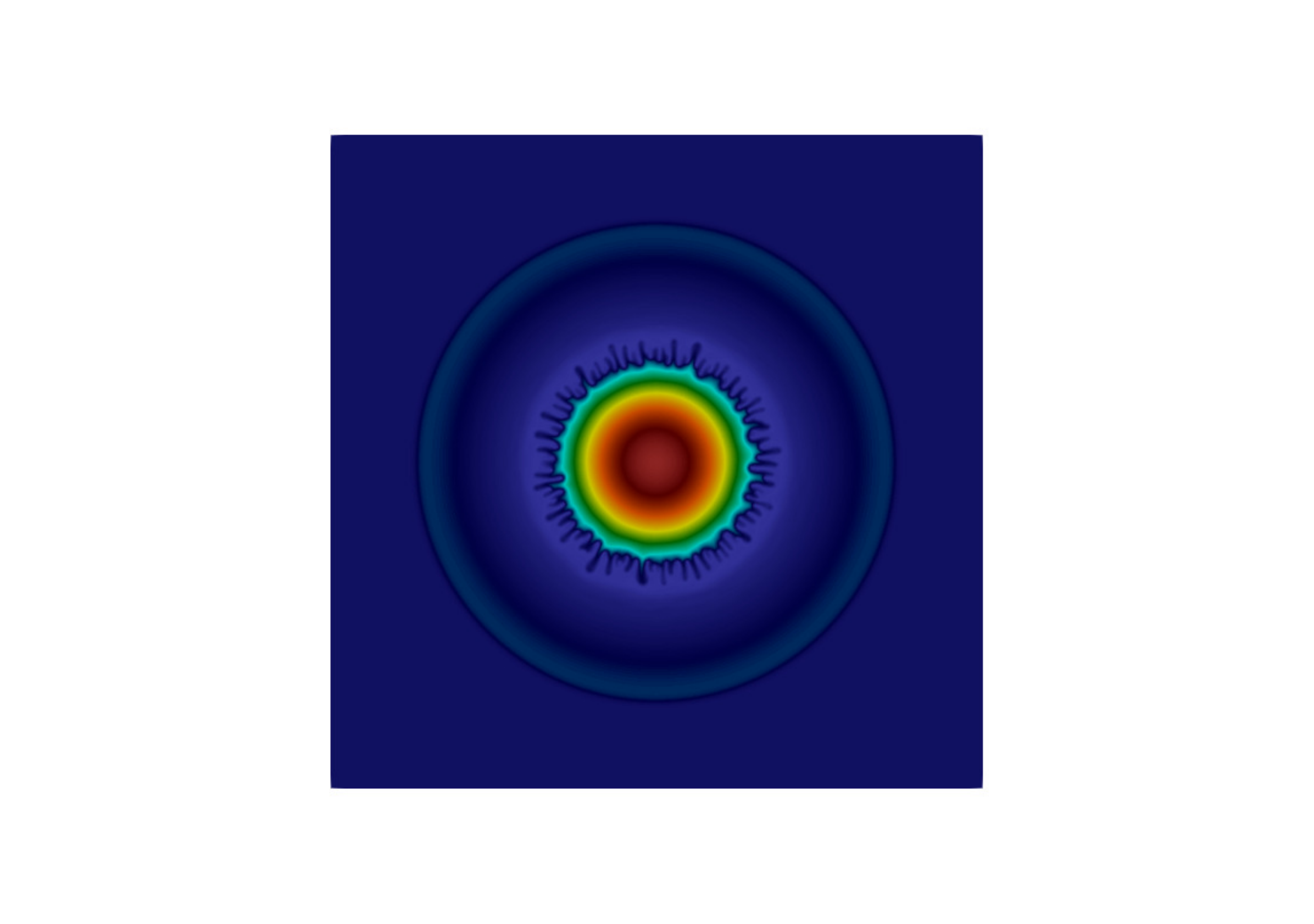}
                		\caption{t=120}
                   \end{subfigure}
                        \begin{subfigure}{0.5\textwidth}
                           \centering        		\includegraphics[width=\linewidth,trim={9cm 3cm 4cm 3cm},clip]{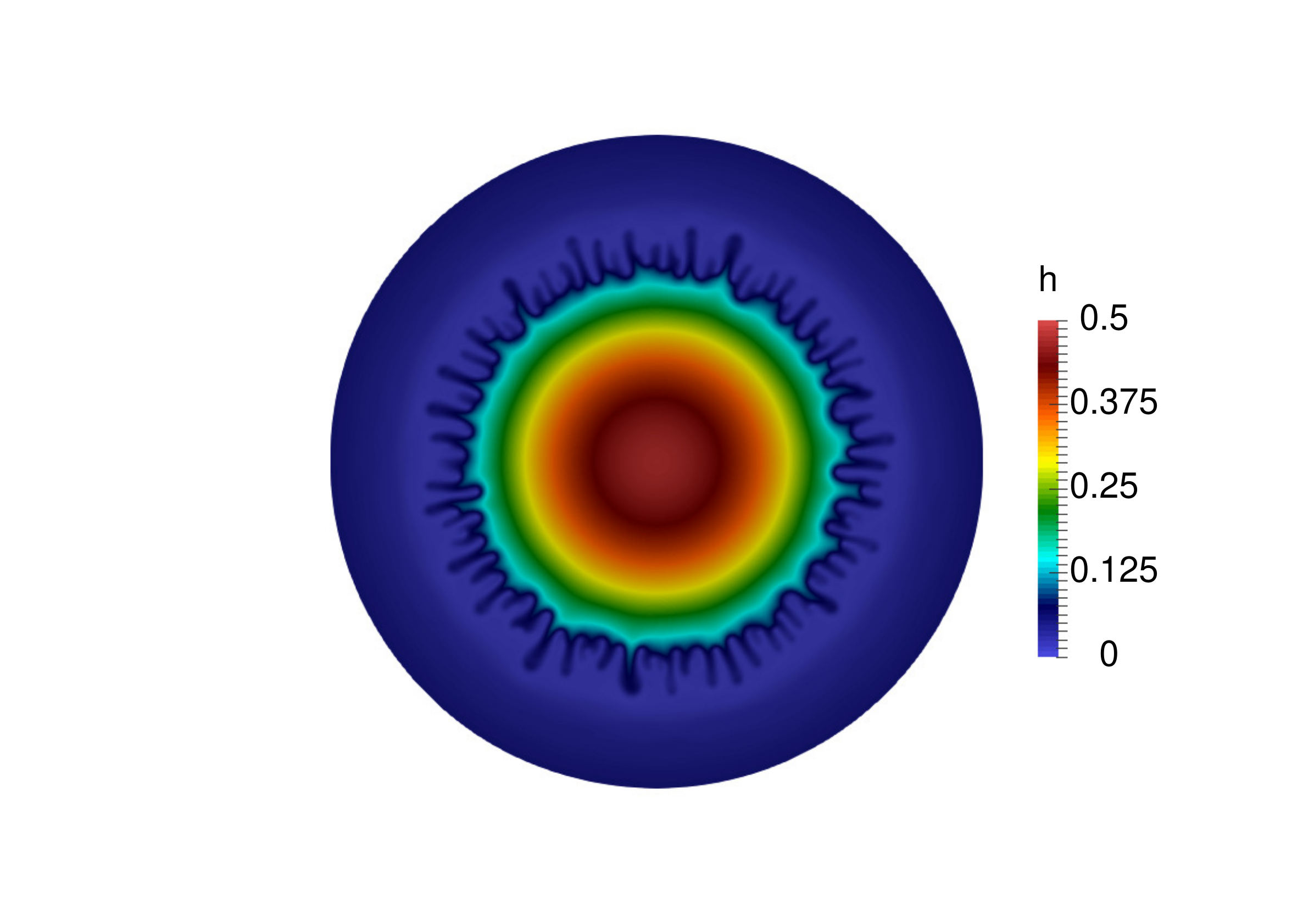}
                           \caption{Zoomed-in view at $t=120$}
                       \end{subfigure} 
    \caption{\bl{Isogeometric solution for fingering instability induced by perturbing the initial height of a thin liquid film rested on a smooth solid substrate at different time instants. Uniform meshes with $512\times512$ elements are employed.}}	
    \label{fig:fingeringInstability-ex03}	 		
    \end{figure}

\begin{figure}
        	\centering
        	\begin{subfigure}{0.3\textwidth}
        	\centering
        	\includegraphics[width=\linewidth,trim={10cm 4cm 10cm 4cm},clip]{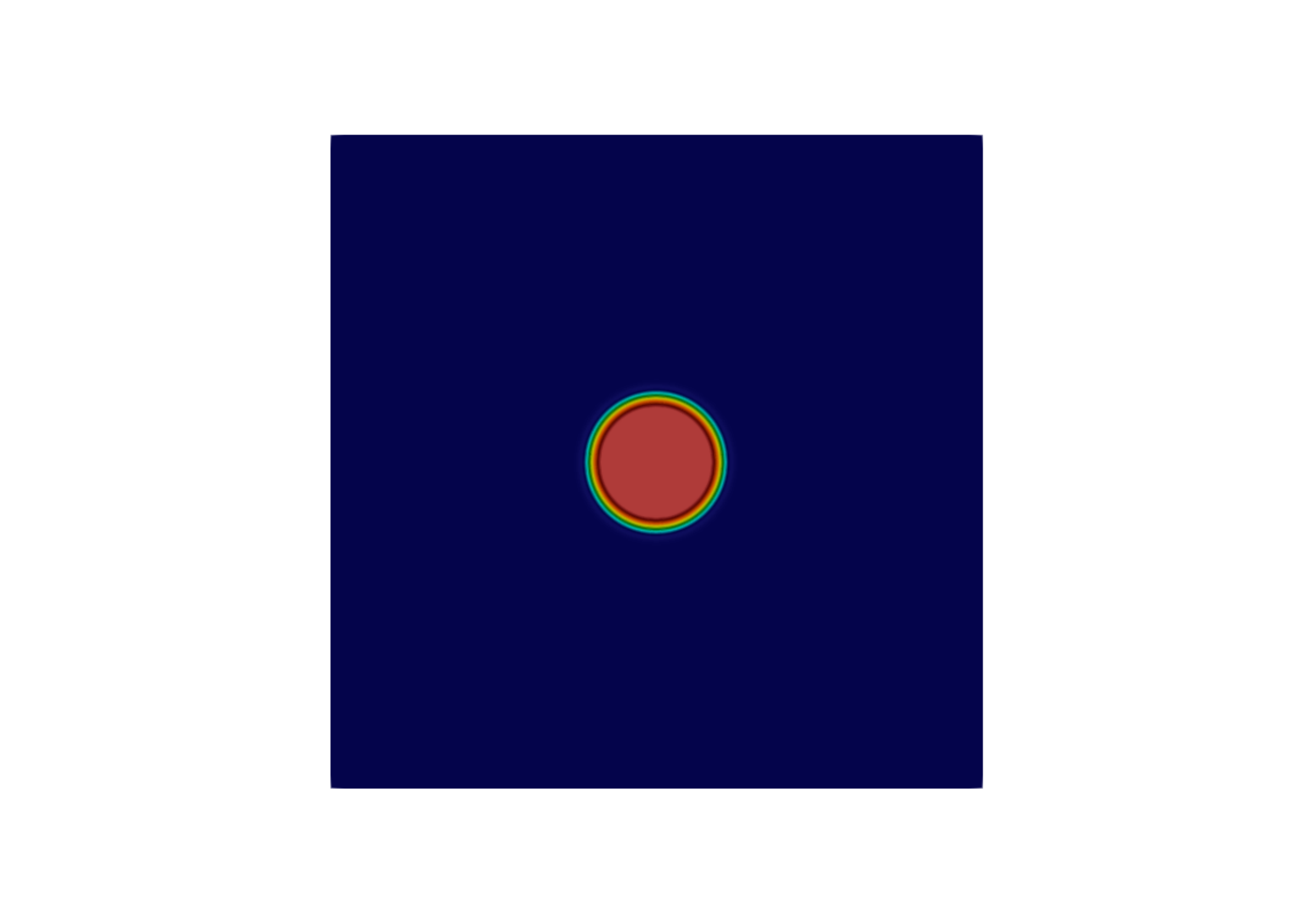}
        	\caption{t=0}
        	\end{subfigure}\quad
        	\begin{subfigure}{0.3\textwidth}
        	\centering
        		\includegraphics[width=\linewidth,trim={10cm 4cm 10cm 4cm},clip]{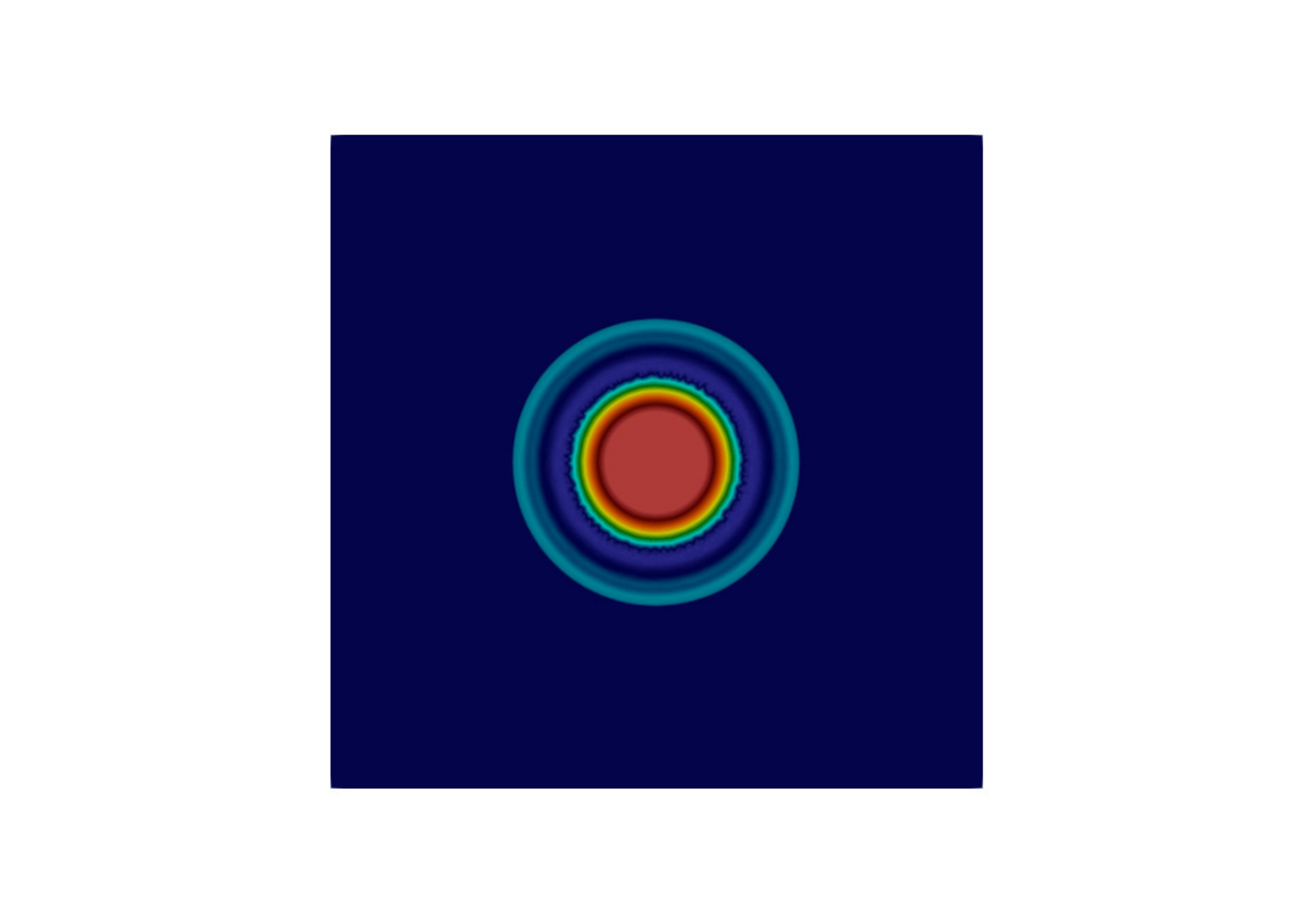}
        		\caption{t=10}
           \end{subfigure} \quad
           \begin{subfigure}{0.3\textwidth}
                   	\centering
                   		\includegraphics[width=\linewidth,trim={10cm 4cm 10cm 4cm},clip]{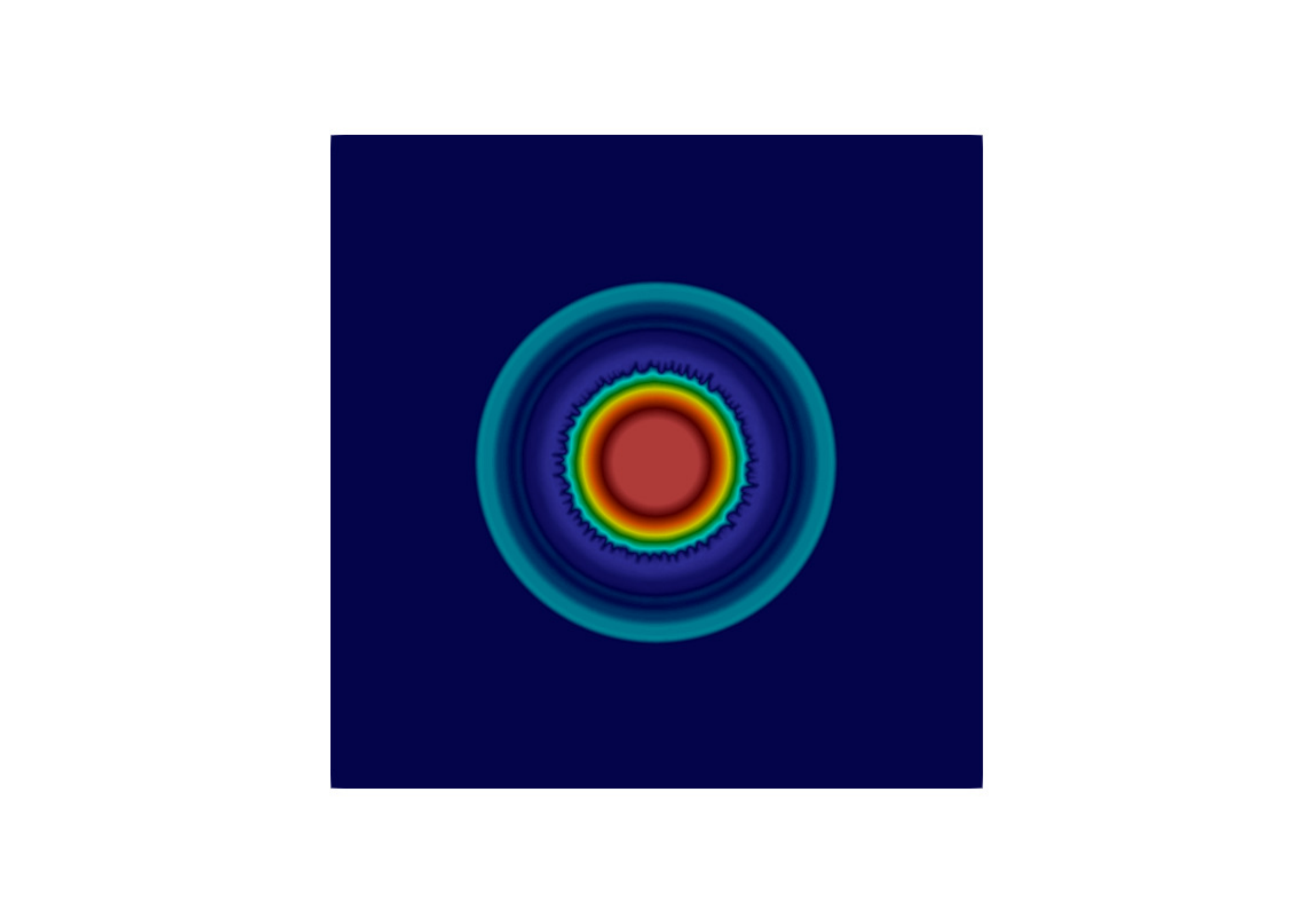}
                   		\caption{t=31}
                      \end{subfigure} \quad
           \begin{subfigure}{0.3\textwidth}
               	\centering
               		\includegraphics[width=\linewidth,trim={10cm 4cm 10cm 4cm},clip]{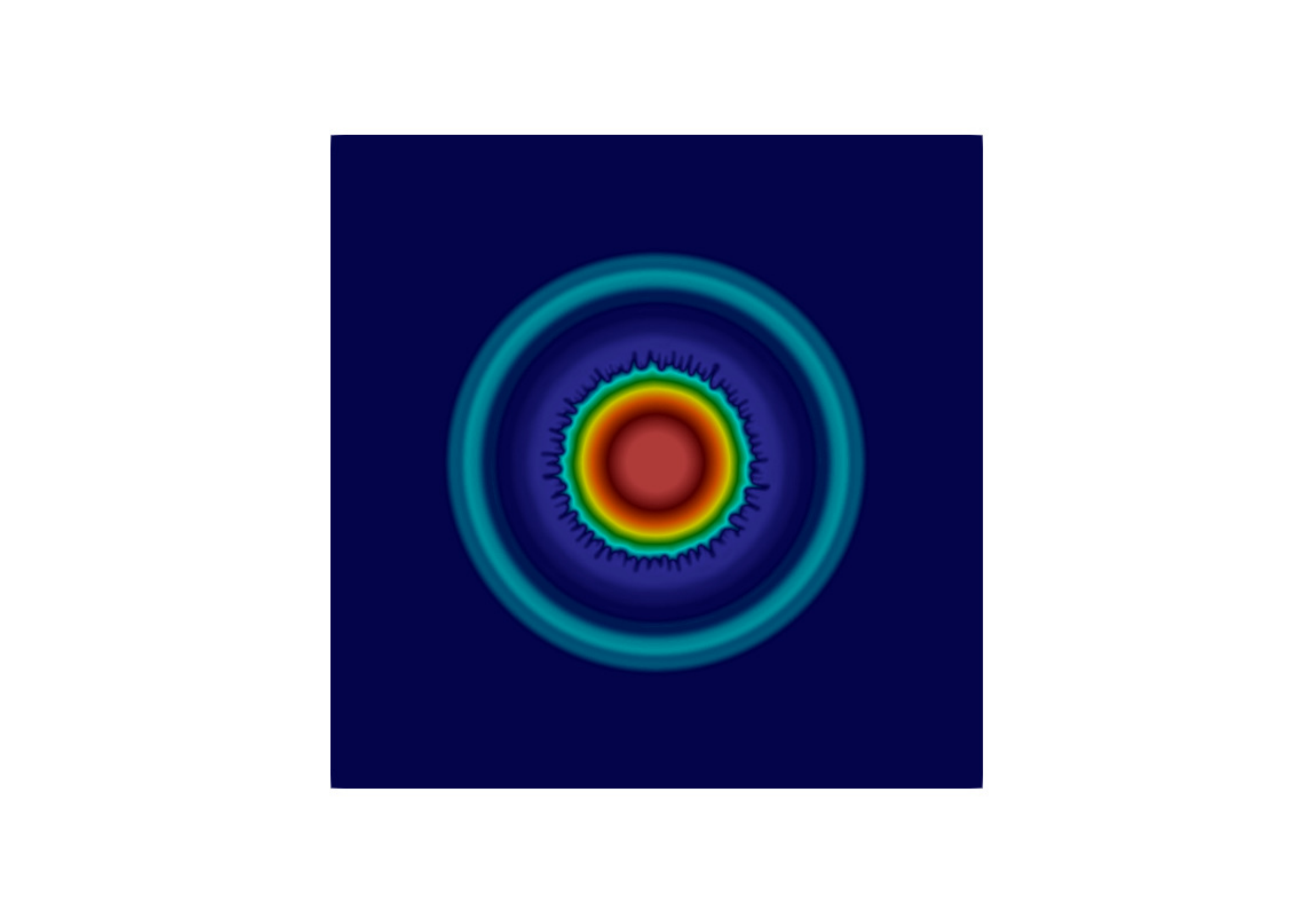}
               		\caption{t=62}
                  \end{subfigure} \quad
           \begin{subfigure}{0.3\textwidth}
               	\centering
               		\includegraphics[width=\linewidth,trim={10cm 4cm 10cm 4cm},clip]{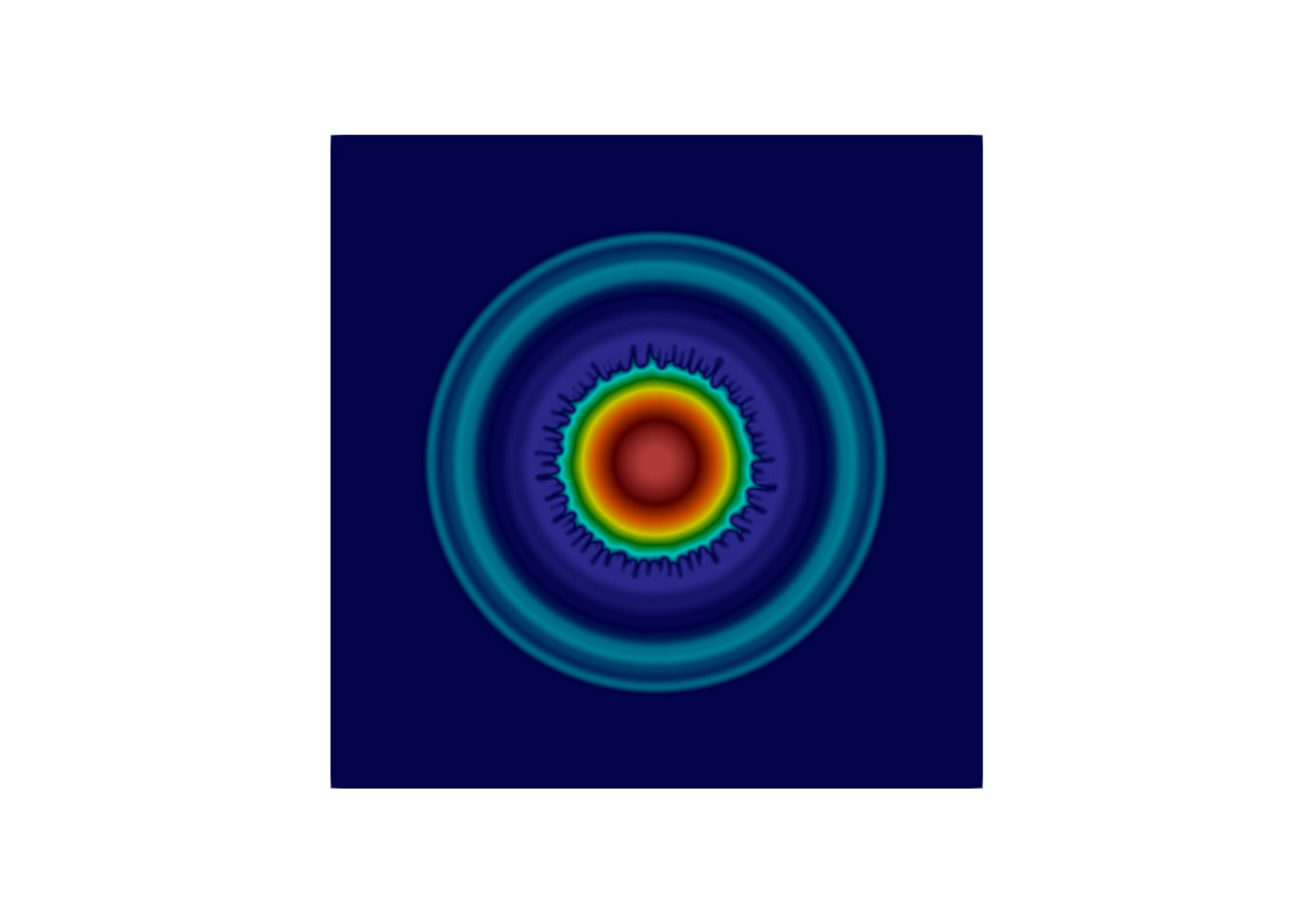}
               		\caption{t=92}
                  \end{subfigure} \quad
         \begin{subfigure}{0.3\textwidth}
                    	\centering
                    		\includegraphics[width=\linewidth,trim={10cm 4cm 10cm 4cm},clip]{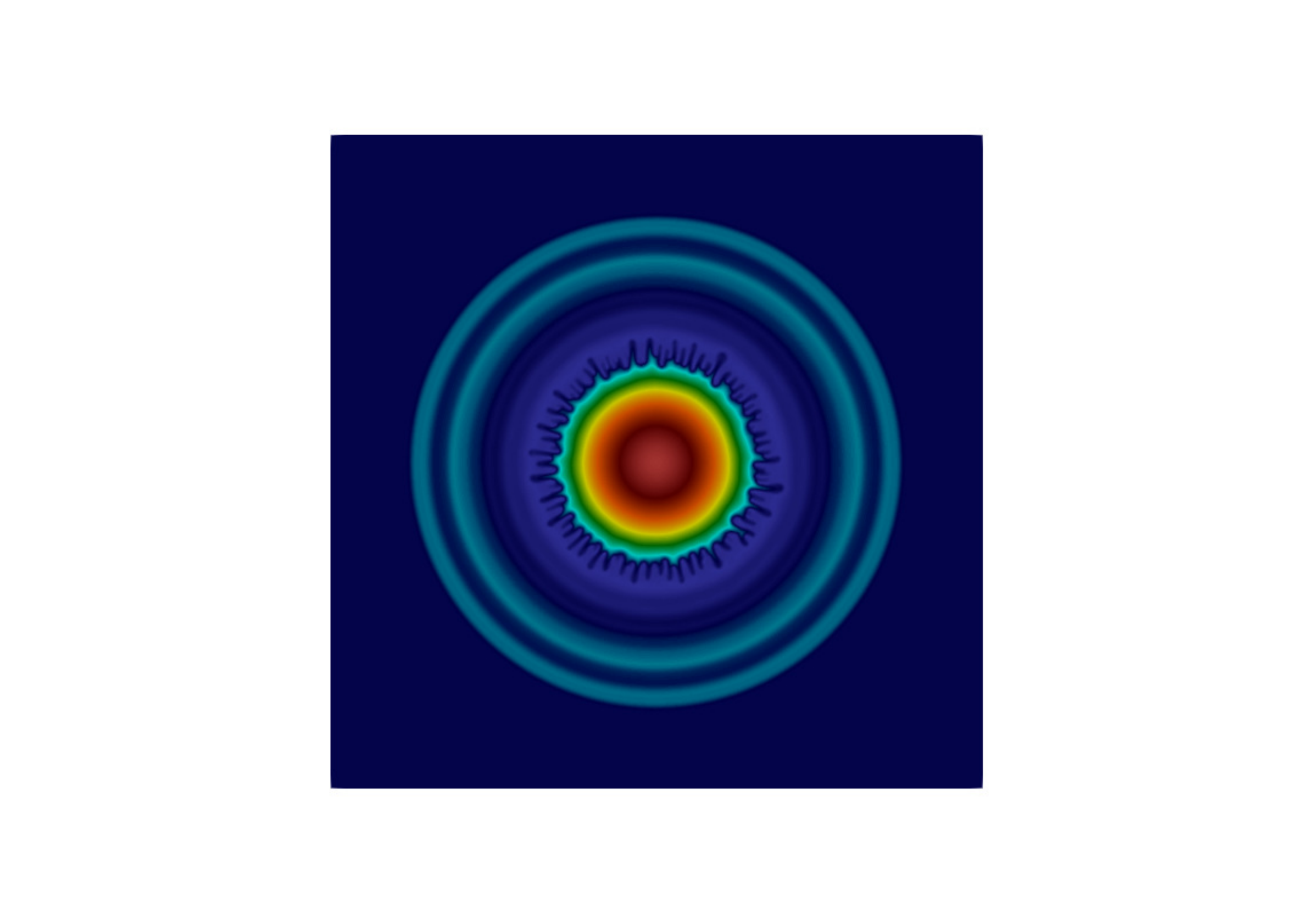}
                    		\caption{t=120}
                       \end{subfigure}
                            \begin{subfigure}{0.5\textwidth}
                               \centering        		\includegraphics[width=\linewidth,trim={9cm 3cm 4cm 3cm},clip]{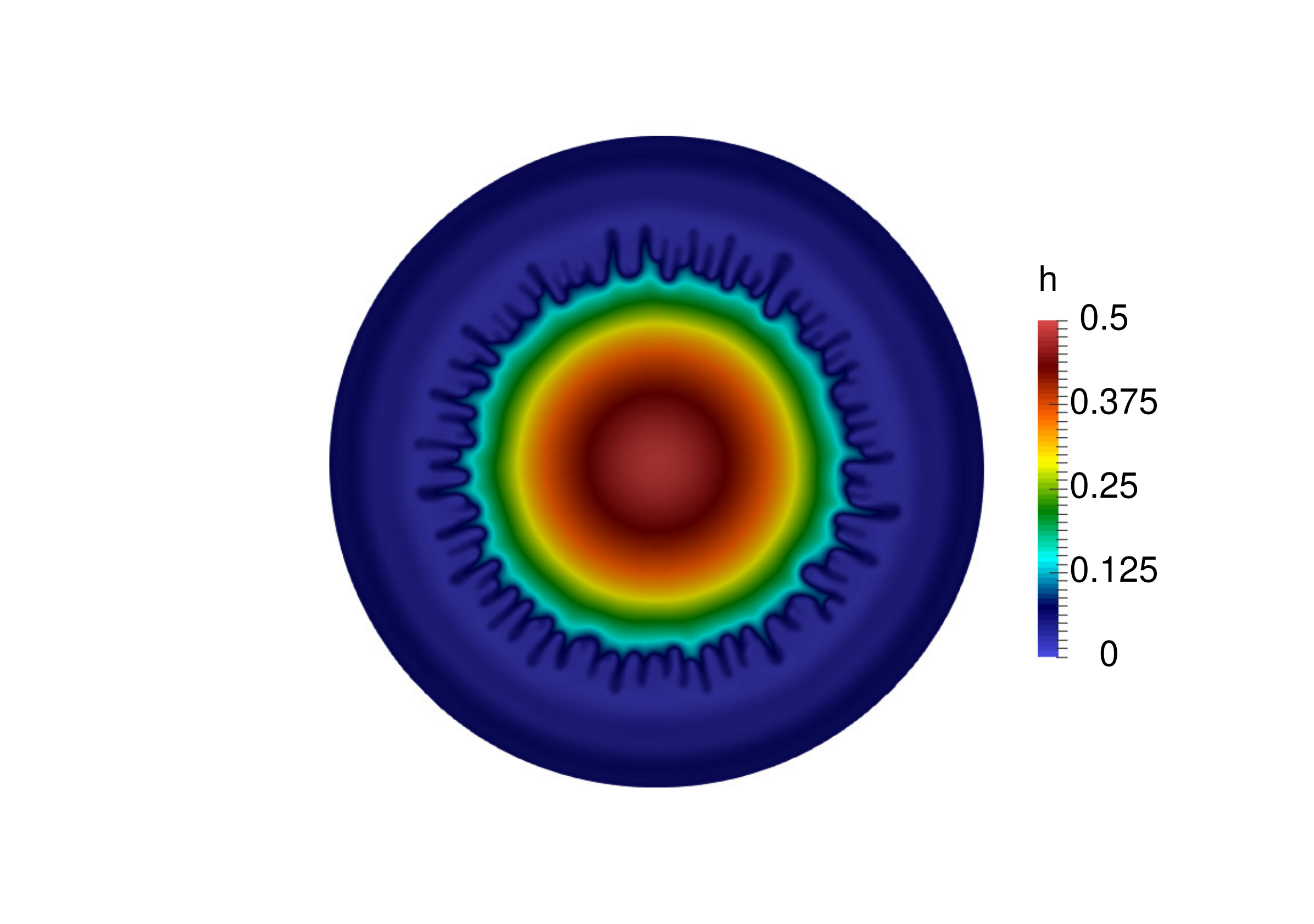}
                               \caption{Zoomed-in view at $t=120$}
                           \end{subfigure} 
        \caption{\bl{Isogeometric solution for fingering instability induced by the roughness of solid substrate at different time instants. Uniform meshes with $512\times512$ elements are employed.}}	
        \label{fig:fingeringInstability-ex04}	 		
        \end{figure}

\begin{figure}
      	\centering
      	\begin{subfigure}{0.48\textwidth}
      	\centering
      	\includegraphics[width=\linewidth,trim={0cm 0cm 0cm 0cm},clip]{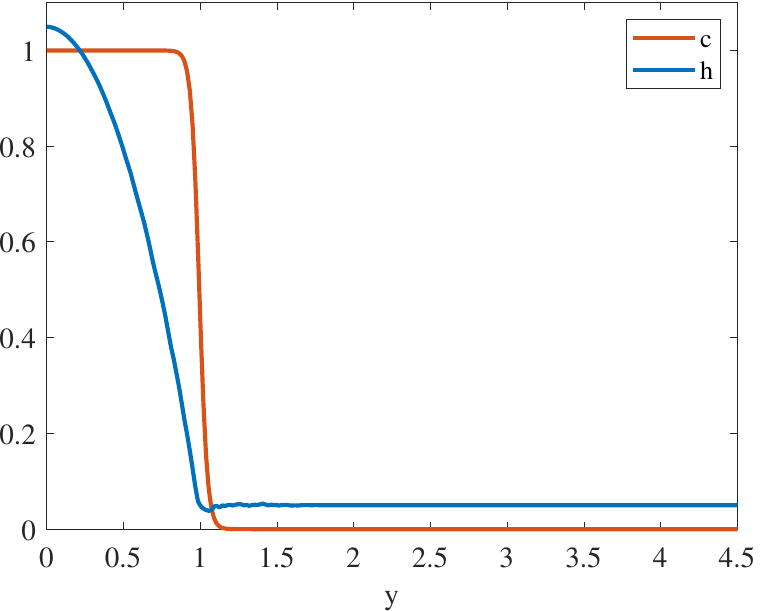}
      	\caption{$t=0$}
      	\end{subfigure}
      	\quad
      	\begin{subfigure}{0.48\textwidth}
      	\centering
      		\includegraphics[width=\linewidth,trim={0cm 0cm 0cm 0cm},clip]{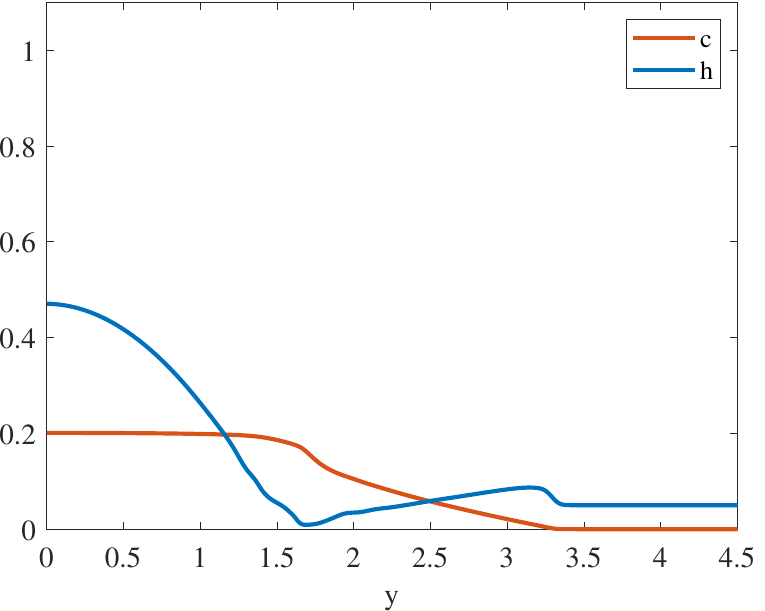}
      		\caption{$t=120$}
         \end{subfigure}\\ 
         \begin{subfigure}{0.48\textwidth}
               	\centering
               	\includegraphics[width=\linewidth,trim={0cm 0cm 0cm 0cm},clip]{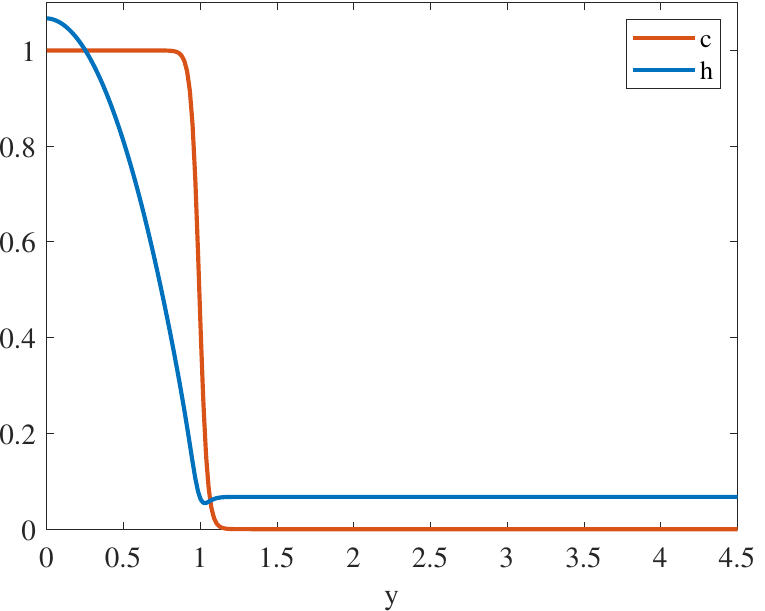}
               	\caption{$t=0$}
               	\end{subfigure}
               	\quad
               	\begin{subfigure}{0.48\textwidth}
               	\centering
               		\includegraphics[width=\linewidth,trim={0cm 0cm 0cm 0cm},clip]{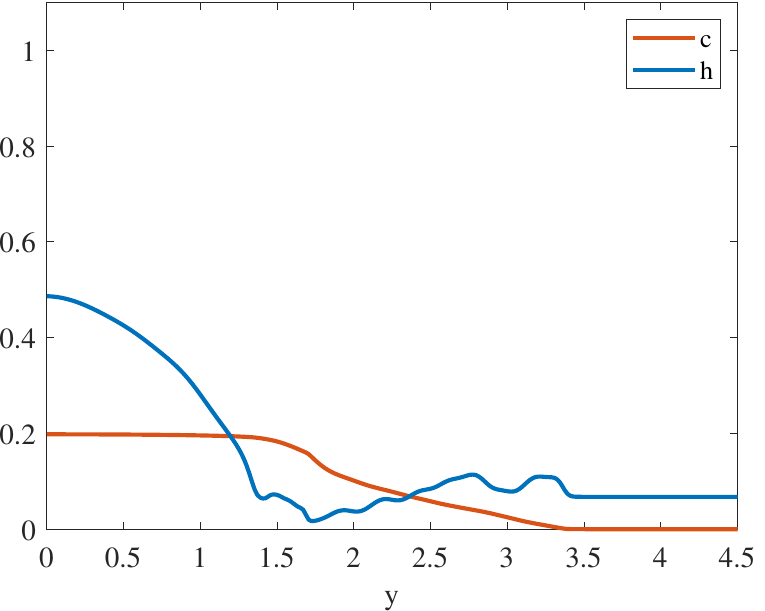}
               		\caption{$t=120$}
                  \end{subfigure}	
      \caption{\bl{The profiles of film height and concentration along the line $x=0$ at time instants $t=0$ and $t=120$ for the smooth (sub-figures (a) and (b)) and rough substrates (sub-figures (c) and (d)).}}	
      \label{fig:fingeringInstability-ex03-ex04-crossSection}	 		
      \end{figure}

\begin{figure}
            	\centering
            	\begin{subfigure}{0.48\textwidth}
            	\centering
            	\includegraphics[width=\linewidth,trim={10cm 4cm 4cm 4cm},clip]{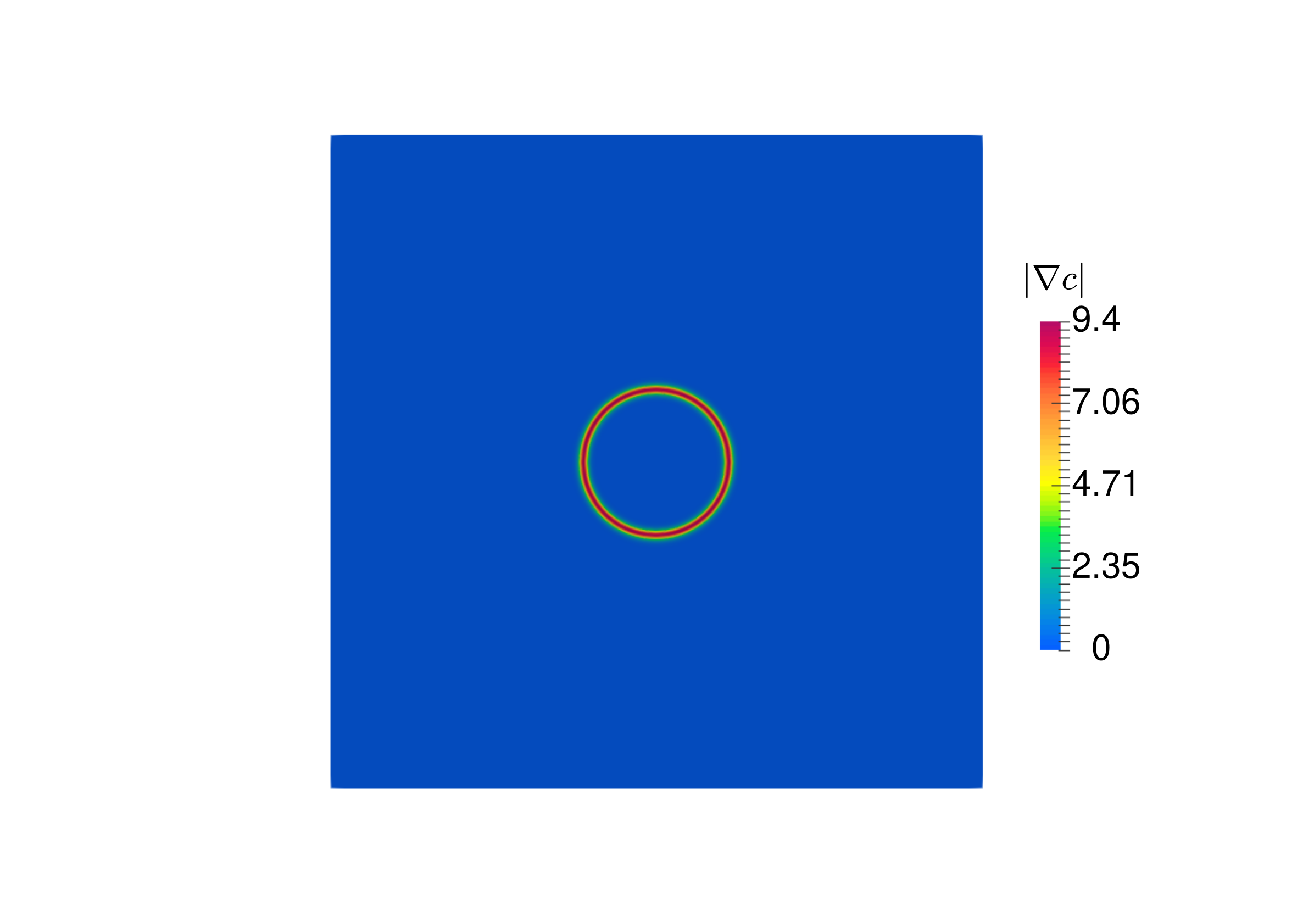}
            	\caption{$t=0$}
            	\end{subfigure}
            	\quad
            	\begin{subfigure}{0.48\textwidth}
            	\centering
            		\includegraphics[width=\linewidth,trim={10cm 4cm 4cm 4cm},clip]{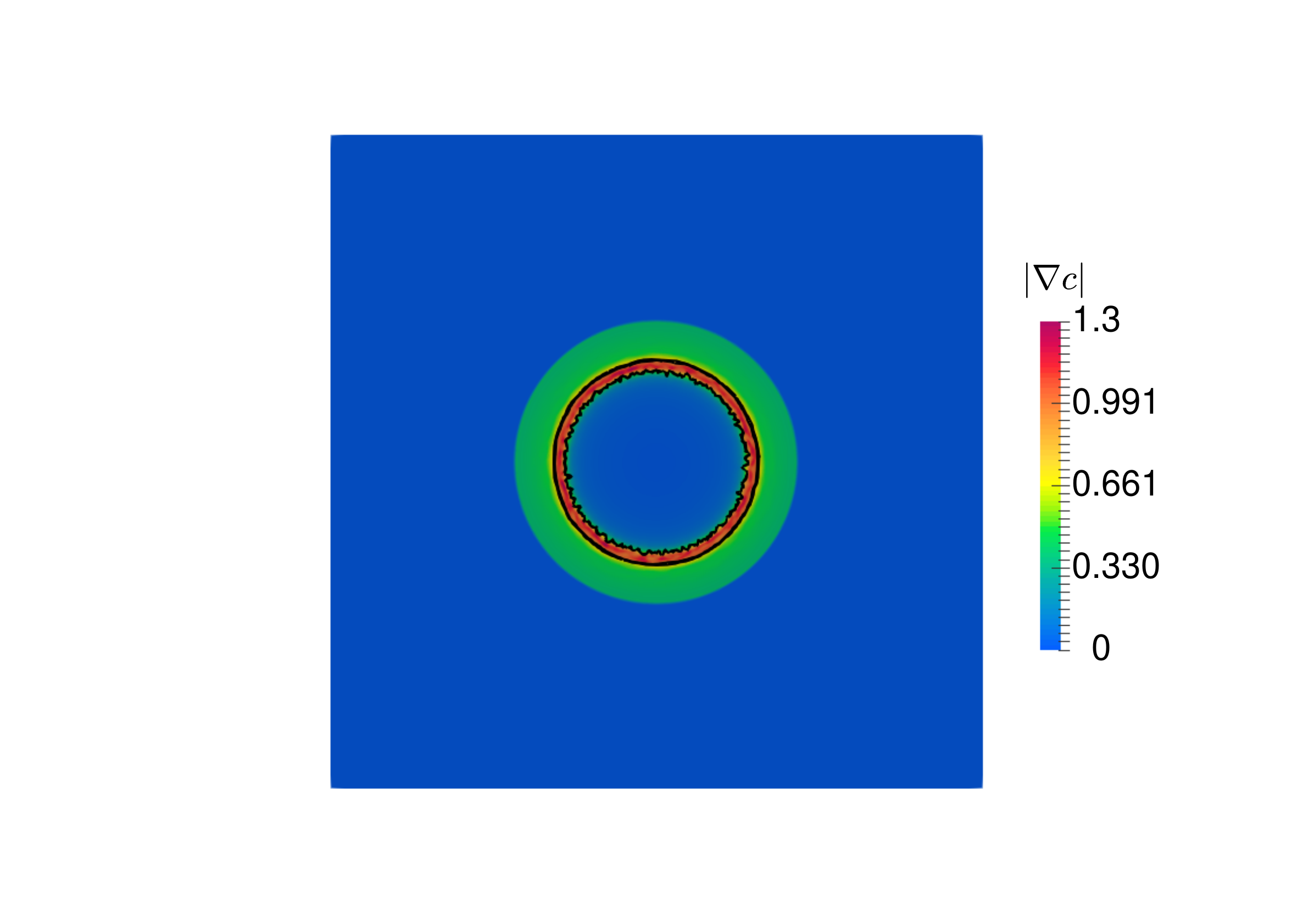}
            		\caption{$t=10$}
               \end{subfigure}	
              \begin{subfigure}{0.48\textwidth}
                    	\centering
                    	\includegraphics[width=\linewidth,trim={10cm 4cm 4cm 4cm},clip]{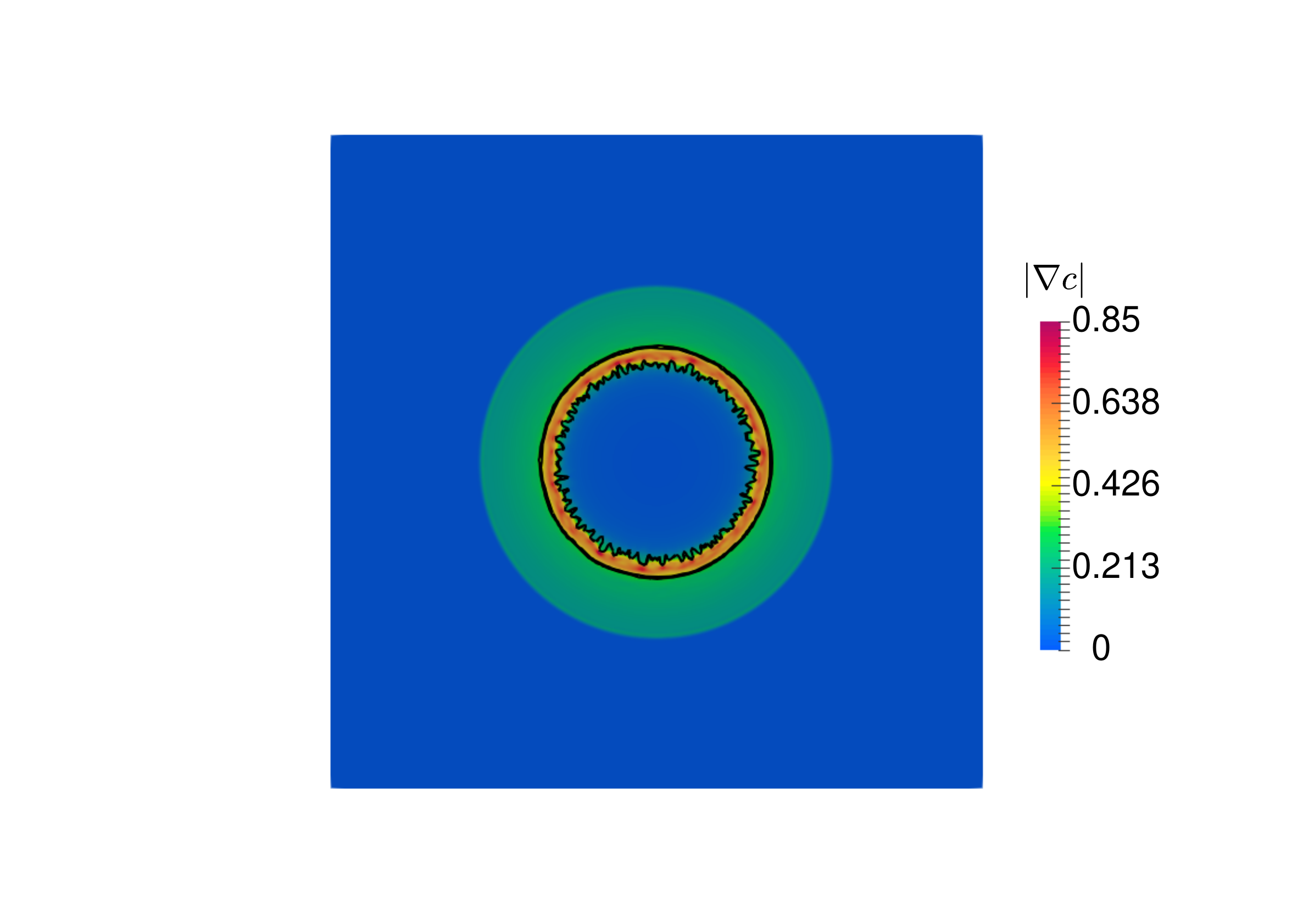}
                    	\caption{$t=30$}
                    	\end{subfigure}
                    	\quad
                    	\begin{subfigure}{0.48\textwidth}
                    	\centering
                    		\includegraphics[width=\linewidth,trim={10cm 4cm 4cm 4cm},clip]{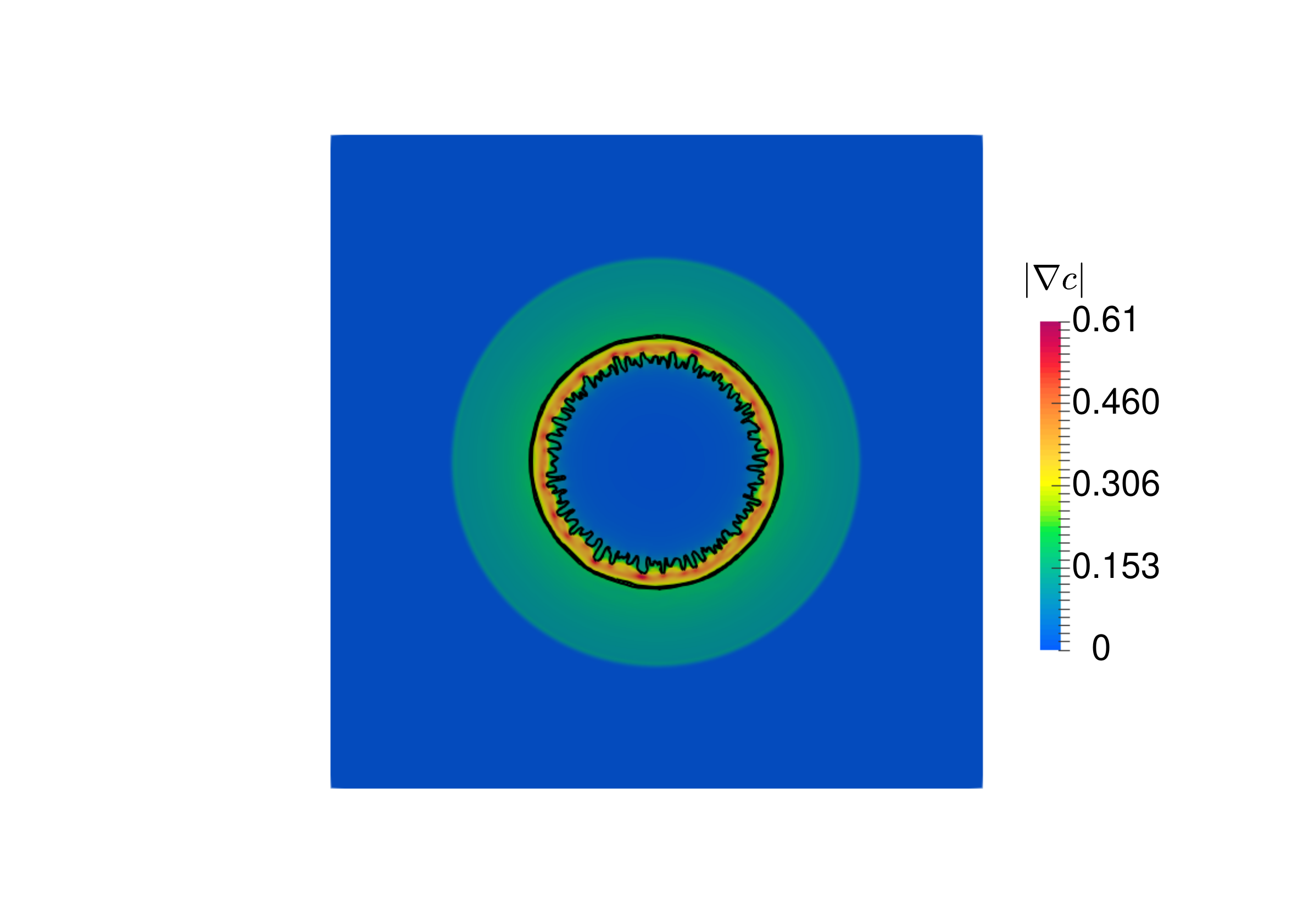}
                    		\caption{$t=59$}
                       \end{subfigure}
             \begin{subfigure}{0.48\textwidth}
                           	\centering
                           	\includegraphics[width=\linewidth,trim={10cm 4cm 4cm 4cm},clip]{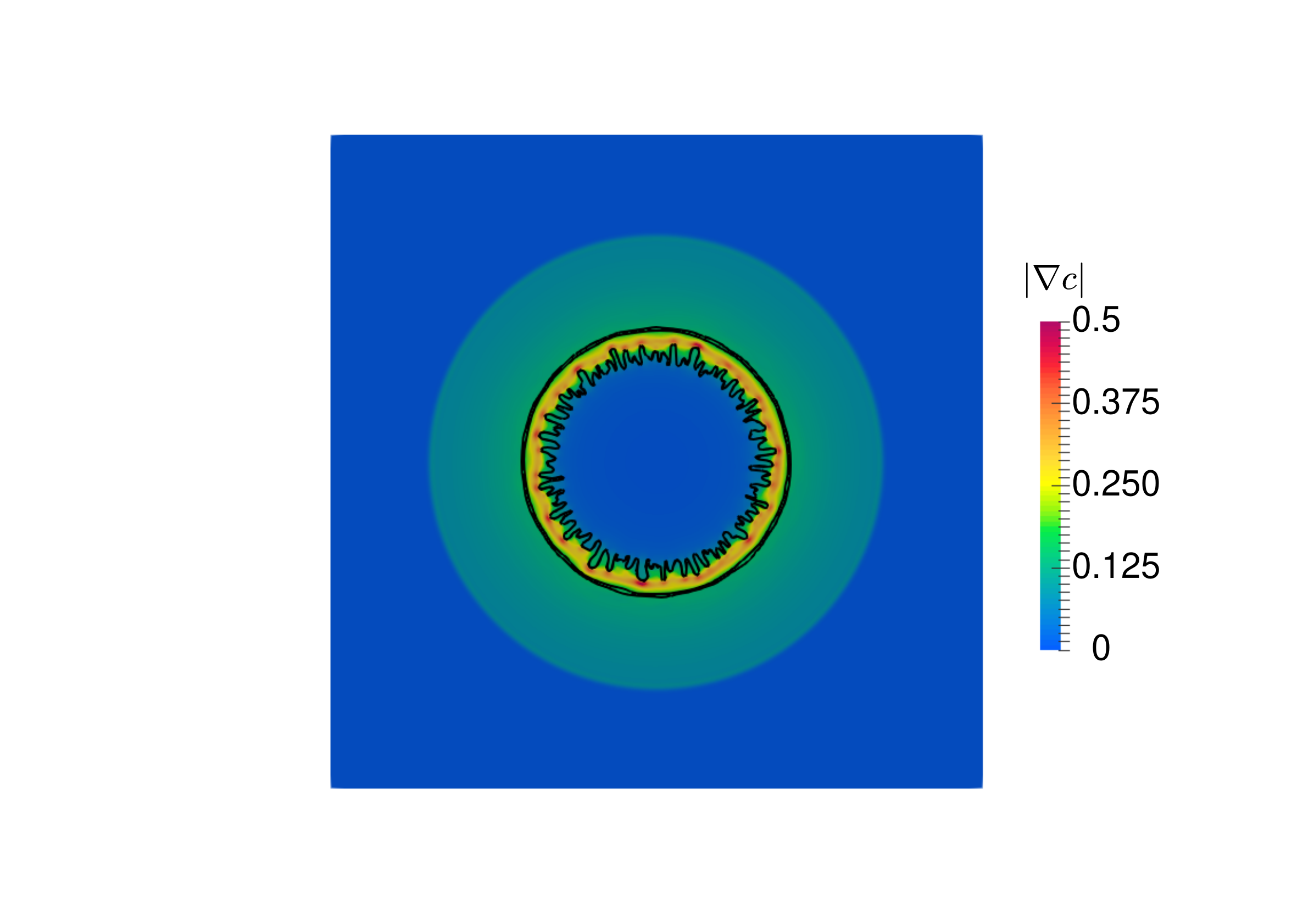}
                           	\caption{$t=93$}
                           	\end{subfigure}
                           	\quad
                           	\begin{subfigure}{0.48\textwidth}
                           	\centering
                           		\includegraphics[width=\linewidth,trim={10cm 4cm 4cm 4cm},clip]{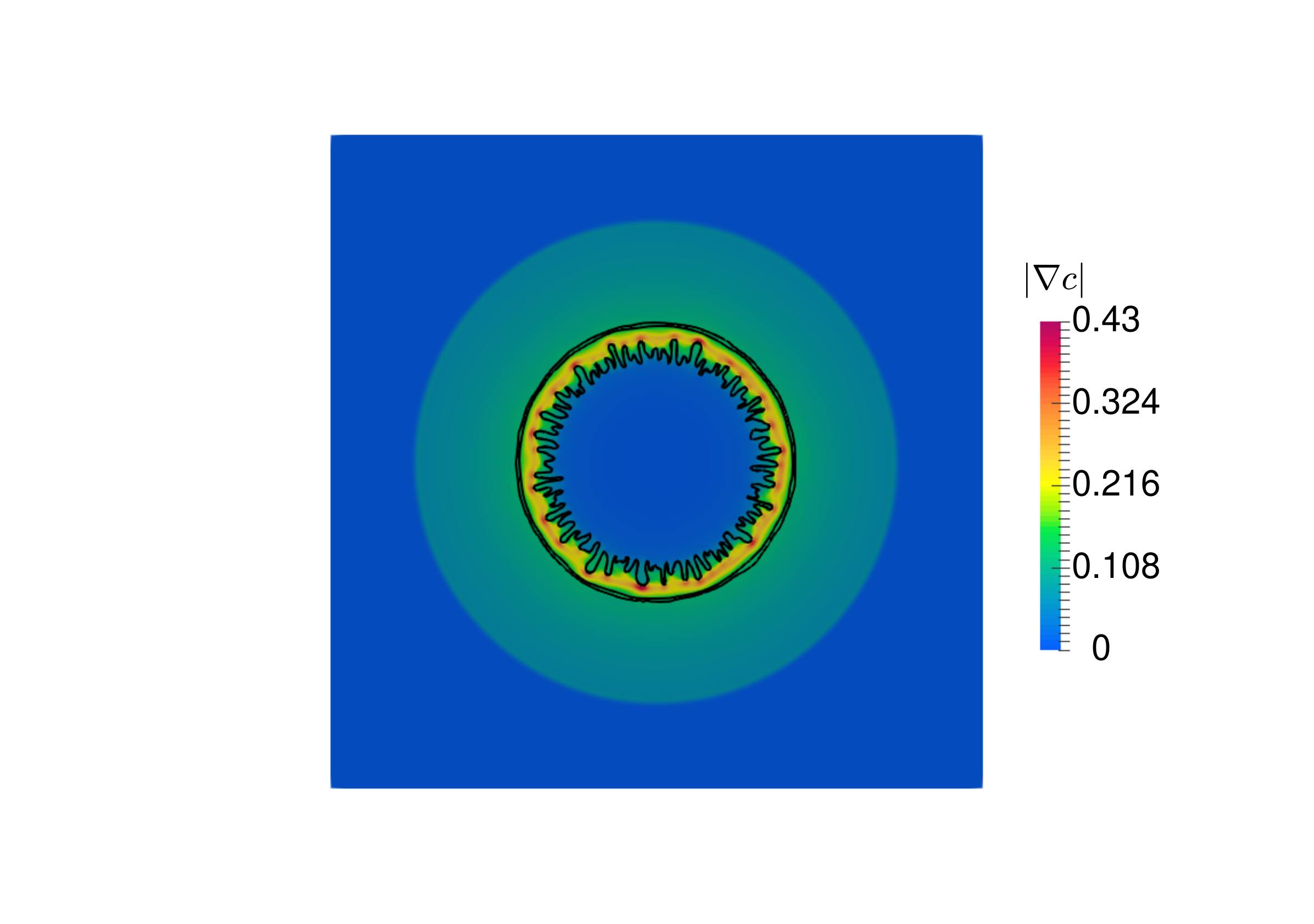}
                           		\caption{$t=120$}
                              \end{subfigure}            
            \caption{\bl{Isogeometric solution for fingering instability initiated by initial height perturbations at different time instants. The height field in the vicinity of the fingers is plotted by black contour lines. The magnitude of gradient of surfactant concentration are plotted at different time instants which show that the higher magnitude of concentration gradient at the tips of fingers is a driving mechanism for fingers to grow. Uniform meshes with $512\times512$ elements are employed.}}	
            \label{fig:fingeringInstability-ex03-3}	 		
            \end{figure}

We now consider the problem of deposition of a drop of surfactant on prewetted smooth and rough substrates. Similar to the previous examples, the instabilities are initiated by introducing perturbations to the initial film height in the case of smooth substrate, and by the roughness of the substrate in the other case. The drop has a radius of unity and the computational domain is a square defined as $[-4.5,4.5]\times[-4.5,4.5]$. 
For the case of a smooth substrate, the roughness function, the initial film height and concentration are defined as
\begin{equation*}
\begin{split}
& f(x,y)=0, \\
& h(x,y,0)=(1-r^{2}+b) H(1-r)+b H(r-1) +A(\varphi)\exp \left(-B(r-1)^{2}\right), \\
& c(x,y,0)=H(1-r).
\end{split}
\end{equation*}
where $r=\sqrt{x^2+y^2}$, $\varphi=\arctan(y/x)$, and $A(\varphi)$ are random
values selected from a uniform distribution on $\left[-0.01,0.01\right]$. Other parameters are the same as those in Sec.~\ref{sec:planarStripOfSurfactant}. For the case of a rough substrate, the roughness function, the initial film height and concentration are defined as
\begin{equation*}
\begin{split}
& f(x,y)=\frac{1}{2}\sum_{k=1}^{N} A_{k}( \cos(\lambda_{k}r)+1) + A(\varphi)\exp\left(-B(r-1)^2\right)-\delta, \\
& h(x,y,0)=(1-r^{2}+b) H(1-r)+b H(r-1) +\frac{1}{2}\bar{A}, \\
& c(x,y,0)=H(1-r).
\end{split}
\end{equation*}
where $\delta=0.005$. In both cases, we use $512\times512$ cubic B-spline elements for discretizing the computational domain \bl{and apply periodic boundary conditions in the $x$ and $y$ directions on the domain boundaries.}. The results of our isogeometric formulation for the fingering patterns on smooth and rough substrates at different time instants are reported in \fref{fig:fingeringInstability-ex03} and \fref{fig:fingeringInstability-ex04}, respectively. The results are very similar to those reported in \cite{LIU2019429} and clearly show the formation of single and branching fingers in a growing-in-time thinned zone. The fingering patterns show a good agreement with the experimentally observed patterns in \cite{Afsar-Siddiqui2003a}.
To better examine the characteristics of spreading processes in these two cases, we show the profiles of film height and concentration in \fref{fig:fingeringInstability-ex03-ex04-crossSection}. The main difference is in the structure of the ramped region which is a smooth line for the smooth substrate and a bumpy line for the rough substrate. Interestingly, in both cases, the peak value of film height in the ramped region is about 1.7 times the height of undisturbed film. \bl{\fref{fig:fingeringInstability-ex03-3} shows the magnitude of the gradient of surfactant concentration, $|\nabla c|$, for the case of height--induced fingering instability at different time instants, while the height field in the vicinity of the fingers, shown by black contour lines, is superimposed on the top of the plot.} It is seen that the higher magnitude of concentration gradient at the tips of some fingers is a driving mechanism for them to spread. 
    
      
    \section{Conclusions} \label{sec:conclusions}
    In this paper, we have described the use of isogeometric analysis for solving a system of time--dependent nonlinear fourth--order PDEs governing the spreading of insoluble surfactants on the free surface of a thin liquid film. Applying a standard Galerkin finite element method to the problem, the weak form involves several terms with second--order spatial partial differential operators which are well--defined and integrable only if the basis functions are piecewise smooth and globally $C^1$ continuous. Exploiting the high continuity of high-order spline basis functions in IGA, the variational formulation of the problem is treated in a straightforward fashion, without the need for mixed finite element methods. \\
    We have first reviewed the theory underlying surfactant spreading on thin films over smooth/rough substrates. Starting from a free boundary problem of
    Navier-Stokes equations for the incompressible flow in a thin liquid film and a surfactant mass balanace
    equation on the free surface, and using
    the lubrication theory together with certain assumptions on different dimensionless parameters, we have reached the final system of coupled PDEs. To solve these PDEs, we have applied a Galerkin--based isogeometric analyis and a generalized--$\alpha$ method for spatial and temporal discretizations, respectively. We have verified our isogeomtric code by considering three sets of examples which include comparisons of surfactant spreading rates with similarity solutions of a simplified \bl{form} of equations, studying the spreading of a surfactant drop on a thin liquid film over a smooth substrate by considering the full PDEs, and investigating the effect of finite--amplitude perturbations to the roughness of the substrate on the surfactant spreading. Next, we have focused on capturing Marangoni-driven fingering instabilities, initiated either by introducing perturbations to the initial film height or to the roughness of the substrate. Several examples, previously solved using mixed finite element methods \cite{LIU2019429}, have been solved in the IGA framework to demonstrate that IGA can easily capture these fingering patterns. Therefore,     
    through several numerical examples, we have shown that isogeometric analysis provides a natural finite element framework for solving high-order PDEs governing the problem of surfactant spreading over thin films on rough/smooth substrates.    
    \section*{Acknowledgements} \label{Acknowledgements} 
     \bl{The authors would like to thank the funding provided by the project H2020-MSCA-RISE-2016  BESTOFRAC (project number: 734370) from the European Union.}
     Navid Valizadeh and Timon Rabczuk would like to acknowledge the financial support of Deutsche Forschungsgemeinschaft (DFG) through the project RA1946/31-1 (project number: 405890576).
\bibliographystyle{model1-num-names}       
 \bibliography{bib_iga}          
\end{document}